\documentclass[a4paper,12pt,thmsa]{article}

\usepackage{amsfonts}
\usepackage{amssymb}
\usepackage{graphicx}
\usepackage{amsmath}
\usepackage{xr}
\usepackage{bezier}
\usepackage{tikz}
\usepackage[normalem]{ulem}

\makeatletter

\renewcommand{\paragraph}{\@startsection{paragraph}{4}{0mm}{-3mm}{-3mm} {\noindent \bf}}

\renewcommand{\subparagraph}{\@startsection{paragraph}{4}{3mm}{-3mm}{-3mm} {\noindent \bf}}

\newcommand{\RR}{{\mathbb{R}}}

\newcommand{\NN}{{\mathbb{N}}}

\newcommand{\CF}{{\mathcal{F}}}

\newcommand{\CV}{{\mathcal{V}}}

\newcommand{\CC}{{\mathcal{C}}}
\newcommand{\CI}{{\mathcal{I}}}

\newcommand{\CS}{{\mathcal{S}}}
\newcommand{\CL}{{\mathcal{L}}}
\newcommand{\CT}{{\mathcal{T}}}
\newcommand{\CB}{{\mathcal{B}}}

\newtheorem{theorem}{Theorem}
\newtheorem{corollary}{Corollary}
\newtheorem{definition}{Definition}
\newtheorem{lemma}{Lemma}
\newtheorem{proposition}{Proposition}

\let\inf\relax \DeclareMathOperator*\inf{\vphantom{p}inf}

\begin{document}

\hyphenchar\font=-1

\baselineskip=7mm

\vspace{3cm}

\author{Jean-Paul D\'ecamps\thanks{Toulouse School of Economics, University of Toulouse Capitole, Toulouse, France. E-mail: \tt{jean-paul.}  \tt{decamps@tse-fr.eu}.}~~~~~~~Fabien Gensbittel\thanks{Toulouse School of Economics, University of Toulouse Capitole, Toulouse, France. E-mail: \tt{fabien.} \tt{gensbittel@tse-fr.eu}.}~~~~~~~Thomas Mariotti\thanks{Toulouse School of Economics, CNRS, University of Toulouse Capitole, Toulouse, France, CEPR, and CESifo. Email: \tt{thomas.mariotti@tse-fr.eu}.}} \vspace{1cm}

\title{\textbf{Mixed-Strategy Equilibria in the War of Attrition under Uncertainty}\thanks{We thank Tiziano De Angelis, Bertrand Gobillard and Stéphane Villeneuve for very valuable feedback. We also thank seminar audiences at Collegio Carlo Alberto and Toulouse School of Economics, as well as conference participants at the 2021 Quimper Workshop on Dynamic Games for many useful discussions. This research has benefited from financial support of the ANR (Programmes d'Investissements d'Avenir CHESS ANR-17-EURE-0010 and ANITI ANR-19-PI3A-0004) and the research foundation TSE-Partnership.
}}

\vspace{1cm}

\maketitle \vspace*{6mm}

\begin{abstract}
We study a generic family of two-player continuous-time nonzero-sum stopping games modeling a war of attrition with symmetric information and stochastic payoffs that depend on an homogeneous linear diffusion. We first show that any Markovian mixed strategy for player $i$ can be represented by a pair $(\mu^i,S^i)$, where $\mu^i$ is a measure over the state space representing player $i$'s stopping intensity, and $S^i$ is a subset of the state space over which player $i$ stops with probability $1$. We then prove that, if players are asymmetric, then, in all mixed-strategy Markov-perfect equilibria, the measures $\mu^i$ have to be essentially discrete, and we characterize any such equilibrium through a variational system satisfied by the players' equilibrium value functions. This result contrasts with the literature, which focuses on pure-strategy equilibria, or, in the case of symmetric players, on mixed-strategy equilibria with absolutely continuous stopping intensities. We illustrate this result by revisiting the model of exit in a duopoly under uncertainty, and exhibit a mixed-strategy equilibrium in which attrition takes place on the equilibrium path though firms have different liquidation values.

\bigskip



\end{abstract}

\thispagestyle{empty}
\newpage
\setcounter{page}{1}

\section{Introduction}

The war of attrition is a workhorse to model situations in which, at each point of time, each player has to decide whether to hold fast or to concede and forfeit a prize to its opponent. Examples of such situations include animal conflict (Maynard Smith (1974)), public good provision (Bliss and Nalebuff (1984)), exit from a declining industry (Ghemawat and Nalebuff (1985), Fudenberg and Tirole (1986)), labor strikes (Kennan and Wilson (1989)), delays in agreement to stabilization policies (Alesina and Drazen (1991)), competition between technological standards (Bulow and Klemperer (1999)), bargaining (Abreu and Gul (2000)), and investment in the presence of informational externalities (D\'ecamps and Mariotti (2004)). Meanwhile, a growing literature attempts to test the theoretical predictions of these models and to estimate the welfare cost of delayed exit decisions (Takahashi (2015)).

Theoretical and empirical applications of war-of-attrition models face several challenges, however. The first is the multiplicity of equilibria, both in pure and mixed strategies, that characterize these models (Riley (1980), Hendricks, Weiss, and Wilson (1988)). It is therefore important to identify theoretical predictions of these models that are robust, in the sense that they hold in a large class of equilibria. The second is to account for observable asymmetries in players' characteristics, which many applied models disregard for simplicity. The third is to allow for stochastic payoffs, so as to capture uncertainty about the future evolution of, say, market conditions. The present paper is an attempt at resolving these issues in a unified framework. In so doing, it identifies a new class of equilibria in mixed strategies that have robust and novel empirical implications.

To this end, we study a generic model of the war of attrition with symmetric information, stochastic payoffs and potentially asymmetric players, which embeds the earlier models of Lambrecht (2001), Murto (2004), Steg (2015), and Georgiadis, Kim, and Kwon (2022). Two players initially present on a market face uncertainty about future market conditions---for instance, the future price of a relevant commodity, or the future state of market demand. Market conditions evolve according to an homogenous linear diffusion. Each player has the option to exit the market, which he may exert at any point in time. Specifically, both players continuously observe the evolution of market conditions; based on this information, each player then decides whether to remain in the market or to irreversibly exit, which terminates the game. In a Markovian way, the players' continuation payoffs when a player decides to exit the market only depend on current market conditions. Besides, there is a second-mover advantage in the sense that, if and when a player exits first, his continuation payoff is lower than the continuation payoff he would have obtained if the other player had exited first given the same market conditions. All payoff-relevant variables---the law of evolution of market conditions and the players' payoff functions---are assumed to be common knowledge. An example of this game is a war of attrition between two firms that may exit a market by liquidating their assets---say, because market demand deteriorates too much---but would meanwhile individually fare better as a monopolist than as a duopolist.

Given the payoff structure we postulate, it is natural to focus on Markov-perfect equilibria in which players' exit decisions at any point in time only depend on current market conditions (Maskin and Tirole (2001)). Our first contribution is to provide a precise definition of Markovian mixed strategies that allows for rich possibilities of randomization for the players. Specifically, our first main result, Theorem \ref{representation}, shows that a randomized stopping time for any player $i$, as defined by Touzi and Vieille (2002) by introducing an auxiliary randomizing device \`a la Aumann (1964), is Markovian if and only if it can be represented by a pair $(\mu^i,S^i)$, where $\mu^i$ is a measure over the state space of the diffusion representing player $i$'s stopping intensity, and $S^i$ is a subset of the state space over which player $i$ stops with probability $1$; the interpretation is that player $i$ exits the market with positive but finite intensity over the support of $\mu^i$, and with infinite intensity over $S^i$. Well-known examples of this representation include pure strategies---that is, stopping times---as in Lambrecht (2001) and Murto (2004), in which the intensity measure $\mu^i$ is degenerate, and mixed strategies in which $\mu^i$ is absolutely continuous with respect to Lebesgue measure, as in Steg (2015) and Georgiadis, Kim, and Kwon (2022). These authors characterize pure-strategy Markov-perfect equilibria, and, in the case of symmetric players, mixed-strategy Markov-perfect equilibria in which players exit the market with the same absolutely continuous intensity measure. In the latter case, attrition is maximal in the sense that each player obtains the payoff he would obtain when facing a stubborn opponent threatening him never to exit the market.

These examples, however, do not exhaust the range of possibilities made available by our formalization of Markovian mixed strategies. In particular, it is possible to conceive of such strategies in which the measure $\mu^i$ is singular with respect to Lebesgue measure. Such strategies need not be artificial nor exotic. For instance, $\mu^i$ may be a Dirac measure at a given point $x^i$ of the state space, weighted by some positive coefficient $a^i$. The interpretation is that, each time market conditions reach $x^i$, player $i$ exits the market with finite intensity $a^i$, a strategy that can easily be obtained as the limit of Markovian mixed strategies defined on discretized state spaces with increasingly fine mesh, or, alternatively, as the limit of Markovian strategies with absolutely continuous intensity measures with supports that degenerate to $x^i$.

Our second main result, Theorem \ref{CNmixed}, precisely shows that, if players are asymmetric---for instance, if firms in a duopoly have different liquidation values---then any mixed-strategy Markov-perfect equilibrium involves strategies with discrete intensity measures. At each point in the support of these measures, the corresponding player is indifferent between exiting and remaining in the market. This implies that the state space is partitioned into intervals in which players alternate between being in a dominated position (with a continuation payoff close to the value he could secure if facing a stubborn opponent) or in a dominant position (with a continuation payoff significantly above that value). Our third main result, Theorem \ref{CSmixed}, finally characterizes these mixed-strategy Markov-perfect equilibria through a variational system satisfied by the two players' continuation value functions. Solving for these equilibria then becomes a rather simple numerical task.

As an illustration, and to show that these necessary and sufficient conditions can be satisfied---in the sense that our variational system has a solution---we provide conditions under which the war of attrition between two duopolists with different liquidation values has a mixed-strategy Markov-perfect equilibrium. In this equilibrium, the firm with the lowest liquidation value randomizes between remaining in the market and exiting at the exit threshold for market conditions that would be optimal if its opponent were stubborn. By contrast, the firm with the highest liquidation value exits with probability 1 if market conditions fall below an even lower threshold, the value of which is determined precisely so as to meet its opponent's indifference condition.

It may rightly be objected that our construction does not solve the multiplicity problem that plagues standard models of the war of attrition: if anything, we exhibit additional equilibria that have been disregarded in the literature (Georgiadis, Kim, and Kwon (2022)). A first answer to this objection is that, in the case of asymmetric players, any mixed-strategy Markov-perfect equilibrium must feature the discrete intensity measures we highlight in Theorem \ref{CNmixed}. Yet there may be many equilibria of this form. It is thus important to identify robust implications of such equilibria that differentiate them from those the literature has focused upon. In that respect, a robust property of the novel equilibria we identify in this paper is that, at any point at which a player randomizes between exiting or remaining in the market, the equilibrium value function of its opponent exhibits a kink. In our duopoly example, the market value of its opponent reaches a peak at this kink, while the randomizing firm's market value goes down to its liquidation value. Moreover, a robust testable implication of such equilibria is that, along any path of the diffusion process modeling the evolution of market conditions, the two firms' market values fluctuate stochastically over the attrition region, moving in opposite directions as long as none of them exits the market. This negative comovement of firms' market values stands in contrast with the symmetric mixed-strategy Markov-perfect equilibrium with an absolutely continuous intensity measure that arises when firms have identical liquidation values, in which firms' market values are constant and equal to their common liquidation value over the attrition region.

\bigskip

\medskip

\noindent {\large \textbf{Related Literature}}

\medskip

\smallskip

\noindent  This paper belongs to the large literature on the war of attrition, starting with the seminal contribution of Maynard Smith (1974) on animal conflict. Ghemawat and Nalebuff (1985) study a war of attrition between duopolists who must decide when to exit from a declining industry. Hendricks, Weiss, and Wilson (1988) offer an exhaustive characterization of pure- and mixed-strategy equilibria in the war of attrition with symmetric information when players have potentially asymmetric payoffs that are deterministic functions of time. Riley (1980), Bliss and Nalebuff (1984), and Fudenberg and Tirole (1986) extend the analysis to asymmetric-information setups where, for instance, a firm is uncertain about its opponent's costs. In the same vein, D\'ecamps and Mariotti (2004) study an investment game that has the structure of a war of attrition because a firm's investment generates additional information for its opponent
about the return of a common-value project.

With the exception of the last paper---which, however, considers a very special Poisson information structure for signals---these papers confine their analysis to situations in which players' payoffs are deterministic. By contrast, a small literature, starting with Lambrecht (2001) and Murto (2004), examines the case where players in a war of attrition have symmetric information, but are uncertain about their future payoffs, which are driven by a diffusion process. Lambrecht (2001) analyzes the order in which firms go bankrupt within a given industry, and how this order is influenced by aggregate economic factors and firm- specific factors such as their financial structure. Murto (2004) studies a stochastic version of Ghemawat and Nalebuff (1985), and shows that a firm with a lower liquidation value may actually end up exiting the market first in equilibrium, despite being a priori more enduring than its opponent. Closest to the present paper in this literature is Georgiadis, Kim, and Kwon (2022). In a setting that extends Murto (2004), they show that, as soon as firms have different liquidation values, there exists no mixed-strategy Market-perfect equilibrium in which firms exit the market according to absolutely continuous intensity measures. They conclude that only pure-strategy Markov-perfect equilibria exist, and therefore that no attrition can actually take place in equilibrium. Our analysis shows that this conclusion is unwarranted once firms can exit the market according to Markovian randomized stopping times with singular intensity measures.

We have borrowed from Touzi and Vieille (2002) our concept of a randomized stopping time, which they introduced to show that continuous-time zero-sum Dynkin games admit a value. A technical contribution of the present paper is to provide a characterization of Markovian randomized stopping times in terms of an intensity measure and a stopping region. This characterization may prove useful for the study of general stochastic timing games.

\vskip 3mm

The paper is organized as follows. Section \ref{TheMod} describes the model. Section \ref{AMaFor} provides rigorous definitions of our strategy and equilibrium concepts, as well as preliminary properties of Markov-perfect equilibria. Section \ref{sec:benchmark} heuristically shows how to construct a mixed-strategy Markov-perfect equilibrium involving a singular intensity measure for one of the players. Section \ref{mainresults} states our main characterization results. Proofs not given in the main text are collected in Appendices A--C.

\section{The Model} \label{TheMod}

\subsection{A General Model of War of Attrition under Uncertainty}

We study a war of attrition with symmetric information between two players, 1 and 2, facing uncertainty about future market conditions. In what follows, $i$ (he) refers to an arbitrary player and $j$ (she) to his opponent. Time is continuous and indexed by $t \geq  0$. Both players observe the evolution of market conditions; based on this information, each player decides whether to \textit{hold fast}, that is, to remain in the market, or to \textit{concede}, that is, to exit the market, an irreversible decision that effectively terminates the game.

The evolution of market conditions is modeled as a one-dimensional time-homogeneous diffusion process $X \equiv (X_t)_{t \geq 0}$ defined over the canonical space $(\Omega,\CF,\mathbf P_x)$ of continuous trajectories with $X_0=x$ under $\mathbf P_x$, that is solution in law to the stochastic differential equation (SDE)
\begin{align} \label{eq1}
\mathrm dX_t = b(X_t) \, \mathrm dt  \, + \, \sigma(X_t) \, \mathrm dW_t, \quad t \geq 0,
\end{align}
driven by some Brownian motion $W \equiv (W_t)_{t \geq 0}$. The state space for $X$ is an interval $\mathcal I \equiv (\alpha, \beta)$, with $- \infty \leq \alpha < \beta \leq \infty$, and $b$ and $\sigma$ are continuous functions, with $\sigma>0$ over $\mathcal I$. We assume that $\alpha$ and $\beta$ are inaccessible (natural) endpoints for the diffusion. Therefore, $X$ is regular over $\mathcal I$ and the SDE \eqref{eq1} admits a weak solution that is unique in law.

Player 1 chooses a (random) time $\tau ^1$ and player 2 chooses a (random) time $\tau^ 2$. Both players discount future payoffs at a constant rate $r>0$. For each $i =1,2 $, the expected payoff of player $i$ is\footnote{By convention, we let $f(X_{ \tau}) \equiv 0$ over $\{\tau = \infty\}$ for any Borel function $f$ and any random time $\tau$.}
\begin{align}
\label{core} J^i(x,\tau^1,\tau^2) &= \mathbf E _x\! \left[1_{\{\tau^i\leq \tau^j\}} \,\mathrm e^{-r \tau^i}R^i (X_{\tau^i}) + 1_{\{\tau^i>\tau^j\}}\, \mathrm e^{-r \tau^j} G^i(X_{\tau^j}) \right]\hskip -1mm .
\end{align}
The payoff functions $R^i$ and $G^i$ in (\ref{core}) are continuous over their domain ${\cal I}$ and satisfy $G^i \geq  R^i$, with $G^i (x) >  R^i(x)$ for $x$ above some threshold $\alpha^i <\beta$.\footnote{Notice that one may have $\alpha^i \leq \alpha$. If $\alpha^i > \alpha$, then $G^i=R^i$ over $(\alpha, \alpha^i]$; this reflects that, for low values of $x$, it may be optimal for player $i$ to exit the market even as a monopolist.} Therefore, if player $i$ concedes at time $\tau^i \leq \tau^j$, he obtains a payoff $R^i(X_{\tau^i})$, whereas, if player $j$ concedes at time $\tau^j< \tau^i$ and $X_{\tau^j} > \alpha^i$, then player $i$ obtains a strictly higher payoff $G^i(\tau^j)$ than the payoff $R^i(\tau^j)$ he would have obtained by conceding at time $\tau^j$. The payoff functions $R^i$ and $G^i$, $i=1,2$, are assumed to be common knowledge among the players; hence the only primitive source of uncertainty in the model is the diffusion process (\ref{eq1}), whose realizations are observed by both players. We study the resulting war of attrition with symmetric information and uncertain payoffs under technical assumptions that we now present.

\subsection{Technical Assumptions}

We first recall useful properties of the solution $X$ to the SDE (\ref{eq1}). We next detail the assumptions on the payoff functions $R^i$ and $G^i$ and emphasize useful properties of the optimal stopping problem
\begin{align} \label{sa}
V_{R^i}(x) \equiv \sup_{\tau \in \mathcal T} \, \mathbf E_x \hskip 0.3mm[ \mathrm e^{-r\tau } R^i(X_\tau)]
\end{align}
faced by player $i$ when player $j$ is \textit{stubborn}, that is, plays $\tau^j = \infty$; here $\CT$ denotes the set of all stopping times of the usual augmentation $(\CF_t)_{t \geq 0}$ of the natural filtration generated by $X$ over the canonical space.\footnote{The definition of $(\CF_t)_{t \geq 0}$ is recalled in Appendix A.} We refer to (\ref{sa}) as player $i$'s \textit{stand-alone exit problem}, in which he cannot benefit from player $j$ conceding.

\paragraph{Properties of the Diffusion $X$}

The infinitesimal generator of the diffusion $X$ is defined for functions $u \in \mathcal C^2(\mathcal I)$ by
\begin{align}
\CL u (x) \equiv b(x) u'(x) + \frac{1}{2}\,\sigma^2(x) u''(x), \quad x \in \mathcal I. \label{defmathcalL}
\end{align}
That $\sigma >0$ over $\mathcal I$ ensures that the ordinary differential equation (ODE) $\CL u - ru=0$ admits a two-dimensional space of solutions in $\mathcal C^2(\mathcal I)$, which is spanned by two positive fundamental solutions $\psi$ and $\phi$, respectively strictly increasing and strictly decreasing, that are uniquely defined up to a linear transformation. By Abel's theorem, the ratio
\begin{align}
\gamma \equiv  \frac{\psi'(x)\phi(x) -\psi(x) \phi'(x) }{S'(x)} >0 \label{wronskian}
\end{align}
of the Wronskian of $\psi$ and $\phi$ and of the derivative of the scale function of the diffusion $X$, which is uniquely defined up to an affine transformation by
\begin{align}
S(x)\equiv \int_c^x \exp\left(-\int_c^y \frac{2 b(z)}{\sigma^2(z)}\, \mathrm dz \right)\mathrm dy, \quad x \in \mathcal I \label{scale}
\end{align}
for some fixed $c \in \mathcal I$, is a constant independent of $x$. Because the boundaries $\alpha $ and $\beta$ of $\mathcal I$ are natural, we know in particular that
\begin{align}
\lim_{x \to \alpha^+} \psi(x)=0, \quad \lim_{x \to \beta^-}\psi(x)=\infty, \quad \lim_{x \to \alpha^+} \phi(x)=\infty, \quad \lim_{x \to \beta^-} \phi(x)=0. \label{hboundaries}
\end{align}
Furthermore, letting $\tau_y \equiv \inf \hskip 0.5mm \{ t \geq 0 : X_t = y \}$ be the hitting time of $y \in \mathcal I$ from $X_0 = x$, we have that
\begin{align}
\label{laplace}  \mathbf E_x \hskip 0.3mm [\mathrm e^{-r\tau_y}] = \left\{ \begin {matrix} \frac{\psi(x)}{\psi(y)} & \text{if} & x \leq y, \\  \frac{\phi(x)}{\phi(y)} & \text{if} & x > y. \end{matrix} \right.
\end{align}

\paragraph{Assumptions on the Payoff Functions $R^i$ and $G^i$}

For each $i = 1, 2$, we assume that $R^i \in \mathcal C^2(\mathcal I)$, and that it satisfies
\begin{itemize}

\item[\bf A1]

For each $x \in \mathcal I$, $\mathbf E _x \hskip 0.3mm  [ \sup_{t \geq 0} \mathrm e^{-rt} | R^i(X_t) | ] < \infty$.

\item[\bf A2]

For each $x \in \mathcal I$, $\lim_{t \to \infty} \mathrm e^{-rt} R^i(X_t) = 0$, $\mathbf P_x$-almost surely.

\item[\bf A3]

There exists $x^i_0 \in \mathcal I$ such that ${\cal L}R^i - rR^i <0$ over $(\alpha, x^i_0)$ and ${\cal L}R^i - rR^i > 0$ over $( x^i_0, \beta)$.

\end{itemize}
A1 guarantees that the family $(\mathrm e^{-r  \tau} R^i(X_{ \tau}))_{\tau \in {\cal T}}$  is uniformly integrable. A1--A2 imply the useful growth property
\begin{align} \label{gp}
\lim_{x \to \alpha^+} \frac{R^i(x)}{\phi(x)} = \lim_{x \to \beta^-} \frac{R^i(x)}{\psi(x)} = 0.
\end{align}
A3 intuitively captures the idea that, as long as the market conditions remain in the portion $(\alpha, x^i_0)$ of the state space, the gains from staying in the market decline if no player has conceded yet. This guarantees that the optimal stopping region $\{x \in \mathcal I: V_{R^i}(x) = R^i(x)\}$ for problem (\ref{sa}) is of the form $(\alpha, x_{R^i}]$ for some threshold $x_{R^i} < x_0^i$, so that
\begin{align}
V_{R^i}(x) = \left\{ \begin {array} {lll} R^i(x) & \text{if} & x \leq x_{R^i}, \\   \frac{\phi(x)}{\phi(x_{R^i})} \, R^i(x_{R^i} )& \text{if} & x > x_{R^i}. \end{array} \right. \label{VR}
\end{align}
The smooth-fit property applies at $x_{R^i}$, that is, $R^{i\prime}(x_{R^i})=\frac{\phi'(x_{R^i})}{\phi(x_{R^i})}\,R^i(x_{R^i})$ (Peskir and Shiryaev (2006), Dayanik and Karatzas (2003, Corollary 7.1)). It follows from standard optimal stopping theory that $(\mathrm e^{-rt} V_{R^i}(X_t))_{t \geq 0}$ is a supermartingale and that ${\cal L}V_{R^i} - rV_{R^i} \leq 0 $ over $\mathcal I \setminus \{x_{R^i}\}$. The following lemma holds.

\begin{lemma} \label{lem_sign}
$V_{R^i}  >0$ over $\mathcal I$ and $R^i >0$ over $(\alpha, x_{R^i}]$.
\end{lemma}

We assume that $G^i \in \mathcal C^1( \mathcal I)$, that $G^i$ is piecewise $\mathcal C^2$ over $\mathcal I$, and that it satisfies

\begin{itemize}

\item[\bf A4]

For each $x \in \mathcal I$, $\mathbf E_x \hskip 0.3mm [ \sup_{t \geq 0} \mathrm e^{-rt} G^i(X_t)] < \infty $.

\item[\bf A5]

For each $x \in \mathcal I$, $\lim_{t \to \infty} \mathrm e^{-rt} G^i(X_t) = 0$, $\mathbf P_x$-almost surely.

\item[\bf A6]

$G^i \geq V_{R^i}$ over $\mathcal I$ and $G^i(x) > V_{R^i}(x)$ if and only if $x> \alpha^i$ for some $\alpha^i < x_{R^i}$.

\item[\bf A7]

${\cal L}G^i - rG^i \leq 0  $ everywhere $G^{i\prime \prime}$ is defined.

\end{itemize}
The interpretation of A7 is that player $i$ would rather obtain the payoff $G^i(X_t)$ sooner than later. This is the case, for instance, when $G^i$ is the value function of an ulterior optimal stopping problem faced by the winner of the war of attrition. From (\ref{sa}) and A6--A7, we have $G^i > R^i \vee 0$ over $\mathcal I$, so that, by Lemma \ref{lem_sign}, $G^i >0$ over $\mathcal I$; hence A4 guarantees that the family $(\mathrm e^{-r  \tau} G^i(X_{ \tau}))_{\tau \in {\cal T}}$ is uniformly integrable. A4--A5 imply the useful growth property
\begin{align} \label{gp'}
\lim_{x \to \alpha^+} \frac{G^i(x)}{\phi(x)} = \lim_{x \to \beta^-} \frac{G^i(x)}{\psi(x)} = 0.
\end{align}

\subsection{A Running Example: Exit in Duopoly} \label{running}

Consider the following model of exit in duopoly, in the spirit of Murto (2004) or Giorgiadis, Kim, and Kwon (2022). Two firms are initially present on the market. As long as both firms remain in the market, each earns a flow duopoly profit $X_t$, where $X$ follows a geometric Brownian motion with drift $b < r$ and volatility $\sigma$,
\begin{align*}
\mathrm dX_t= b X_t \, \mathrm dt + \sigma X_t \, \mathrm dW_t, \quad t\geq 0,
\end{align*}
over the state space $\mathcal I \equiv(0, \infty)$. If firm $i$ concedes at time $\tau^i$, then its assets are liquidated for a value $l^i >0$, while firm $j$ enjoys from time $\tau^i$ on a flow monopoly profit $mX_t$ for some $m >1$, until it in turn decides to exit the market and receive its liquidation value $l^j >0$. Thus the expected discounted profit of every firm $i$ for given exit times $\tau^i$ and $\tau^j$ is
\begin{align*}
F^i(x,\tau^1, \tau^2) \equiv \mathbf E_x \! \left[ \int_0^{\tau_1 \wedge \tau_2} \mathrm e^{-rt} X_t \, \mathrm dt + 1_{\{\tau^i \leq \tau^j\} } \, \mathrm e^{-r \tau^i} l^i
 +  1_{\{\tau^i > \tau^j\} } \, \mathrm e^{-r \tau_j} V^i_m(X_{\tau^j}) \right]\hskip -1mm ,
\end{align*}
where $V^i_m$ is firm $i$'s value function as a monopolist,
\begin{align*}
V^i_m(x) \equiv \sup _{\tau \in \mathcal T} \, \mathbf E_x \! \left[ \int_0^\tau \mathrm e^{-rt} mX_t \, \mathrm dt + \mathrm e^{-r \tau} l^i\right]\hskip -1mm .
\end{align*}
Letting $E(x) \equiv \mathbf E_x \!\left [\int_0^\infty \mathrm e^{-rt} X_t \, \mathrm dt  \right ] \!= \frac{x}{r - b}$, $R^i \equiv  l^i - E$, and $G^i \equiv V^i_m - E$, we obtain the expression (\ref{core}) for $J^i (\cdot ,\tau^1, \tau^2) \equiv F^i(\cdot ,\tau^1, \tau^2) - E$. Standard computations (see, for instance, Dixit and Pindyck (1994)) yield
\begin{align*}
x_{R^i}= {\rho^- \over \rho^- -1} \, (r-b) l^i \quad \mbox{and} \quad \alpha^i = {x_{R^i}\over m},
\end{align*}
where
\begin{align*}
\rho^- \equiv \frac{1}{2} - \frac{b}{\sigma^2} - \sqrt{\left( \frac{1}{2} - \frac{b}{\sigma^2}\right) ^2 + \frac{2r }{\sigma^2}}.
\end{align*}
Notice that $G^i(x) = R^i(x) = l^i - E(x)$ for all $x \in (\alpha, \alpha^i]$. It is easy to check that this specification satisfies A1--A7. We will use it in Section 4 to illustrate our results.

\section{Mixed Strategies and Equilibrium Concept} \label{AMaFor}

Our key methodological contribution is to allow players to play randomized stopping times. We first recall the definition and basic properties of randomized stopping times. Imposing a Markov restriction leads to our first main result, which is a representation theorem for Markov randomized stopping times. We then define the concept of Markov-perfect equilibrium and give some important properties of best replies.

\subsection{Randomized Stopping Times}

One classical definition of a randomized stopping time consists, following Aumann (1964), in enlarging the probability space; this compensates for the absence of a natural measurable structure over the space of stopping times. For every player $i= 1,2$, the corresponding enlarged probability space is $(\Omega ^i ,\mathcal F^i) \equiv (\Omega\times [0,1], \mathcal F \otimes \mathcal B([0,1]))$, endowed with the product probability $\mathbf P^i_x \equiv \mathbf P_x \otimes Leb$, where $Leb$ denotes Lebesgue measure. We borrow the following definition from Touzi and Vieille (2002).

\begin{definition}\label{def:randomized_stopping_time}
A {\rm randomized stopping time} for player $i=1,2$ is a $\CF \otimes \CB([0,1])$-measurable function $\gamma^i: \Omega^i \to \mathbb R_+$ such that$,$ for $Leb$-almost every $u^i\in [0,1],$ $\gamma^i(\cdot,u^i) \in \CT$. The process $\Gamma^i \equiv (\Gamma^i_t)_{t \geq 0}$ defined by
\begin{align} \label{ccdf'}
\Gamma^i_t (\omega) \equiv \int_{[0,1]} 1_{\{\gamma^i(\omega,u^i) \leq t\}} \, \mathrm du^i, \quad (\omega,t) \in \Omega \times \mathbb R_+,
\end{align}
is the {\rm conditional cumulative distribution function} (ccdf) of the randomized stopping time $\gamma^i$. The process $\Lambda^i \equiv (\Lambda^i_t)_{t \geq 0}$ defined by
\begin{align} \label{cst}
\Lambda^i_t(\omega) \equiv 1- \Gamma^i_t(\omega) , \quad (\omega,t) \in \Omega \times \mathbb R_+,
\end{align}
is the {\rm conditional survival function} (csf) of the randomized stopping time $\gamma^i$.
\end{definition}

It is immediate that the ccdf process $\Gamma^i$ defined by (\ref{ccdf'}) takes values in $[0,1]$ and has nondecreasing and right-continuous trajectories. The following lemma shows that the process $\Gamma^i$ is adapted and provides a useful representation.

\begin{lemma}\label{lemma_properties_lambda}
The ccdf process $\Gamma^i$ is $(\CF_t)_{t\geq 0}$-adapted and$,$ for $\mathbf P_x$-almost every $\omega \in \Omega,$
\begin{align}\label{ccdf2}
\Gamma^i_t (\omega)= \mathbf P^i_x \hskip 0.3mm [\gamma^i \leq t \! \mid \! \CF_t] (\omega)
\end{align}
for all $x\in \CI$ and $t \geq 0$.
\end{lemma}

We adopt the convention $\Gamma^i_{0-}\equiv0$. This allows us in what follows to interpret integrals of the form $\int_{[0, \tau)} \cdot \, \mathrm d \Gamma^i_t$ in the Stieltjes sense for any ccdf $\Gamma^i$.

If the players use randomized stopping times $\gamma^1$ and $\gamma^2$, then their expected payoffs are defined over the product probability space $\Omega\times [0,1] \times [0,1]$ with canonical element $(\omega,u^1,u^2)$ endowed with the product probability $\overline{\mathbf P}_x \equiv \mathbf P_x \otimes Leb \otimes Leb$ by
\begin{align}
J^i(x,\gamma^1,\gamma^2) \equiv \overline {\mathbf E }_x\! \left[1_{\{\gamma^i\leq \gamma^j\}} \,\mathrm e^{-r \gamma^i}R^i (X_{\gamma^i}) + 1_{\{\gamma^i>\gamma^j\}}\, \mathrm e^{-r \gamma^j} G^i(X_{\gamma^j}) \right]\hskip -1mm, \label{BDCOLL}
\end{align}
where $\gamma^1 \equiv\gamma^1(\omega,u^1)$ and $\gamma^2 \equiv \gamma^2(\omega,u^2)$, reflecting that player $1$ and player $2$ use the independent randomization devices $u^1$ and $u^2$, respectively.

Our next result shows that we may equivalently work with the family of ccdf processes $\Gamma^i$. As similar results appear elsewhere in the literature (Touzi and Vieille (2002), Riedel and Steg (2017)), its proof is relegated to Appendix A.

\begin{lemma} \label{ccdf}
If the players use randomized stopping times with ccdf $\Gamma^1$ and $\Gamma^2,$ then their expected payoffs write as
\begin{align}
J^i(x,\Gamma^1,\Gamma^2) =  \mathbf E _x \! \left[ \int_{[0,\infty)} \mathrm e^{-r t} R^i(X_t)\Lambda^j_{t-} \, \mathrm d\Gamma^i_t+ \int_{[0,\infty)} \mathrm e^{-r t}G^i(X_t)\Lambda^i_t \, \mathrm d\Gamma^j_t\right] \hskip -1mm. \label{forbpr}
\end{align}
Moreover$,$ any nondecreasing, right-continuous$,$ $\CF_t$-adapted$,$ $[0,1]$-valued process $\Gamma^i$ is the ccdf of the randomized stopping time $\hat\gamma^i$ defined by
\begin{align}\label{inverse_stopping}
\hat\gamma^i(u^i) \equiv \inf \hskip 0.5mm \{ t \geq 0 : \Gamma^i_t > u^i \}.
\end{align}
\end{lemma}

\subsection{Markovian Randomized Stopping Times} \label{MRST}

Our goal in this paper is to characterize equilibria in which players concede according to Markov randomized strategies that only depend on current market conditions. Notice that such strategies have to be defined for any initial market conditions $x \in \CI$. As a preliminary, we need the following standard definition (Revuz and Yor (1999, Chapter I, \S3)).

\begin{definition}
Let $Y \equiv (Y_t) _{t \geq 0}$ be the {\rm coordinate process} over the canonical space $\Omega,$ defined by $Y_t(\omega) \equiv \omega_t$ for all $\omega \in \Omega$ and $t \geq 0$. Then$,$ for each $t \geq 0,$ the {\rm shift operator} $\theta_t: \Omega \to \Omega$ is defined by $Y_s \circ \theta_t \equiv Y_{s +t} $ for all $s\geq 0$.
\end{definition}

In words, the effect of $\theta_t$ on a trajectory $\omega$ is to forget the part of the trajectory prior to time $t$ and to shift back the remaining part by $t$ units of time. We are now ready to define our notion of a Markovian randomized stopping time.

\begin{definition} \label{marksur}
A randomized stopping time for player $i=1,2$ with csf $\Lambda^i: \Omega \times \RR_+ \to [0,1]$ is {\rm Markovian} if$,$ for $\mathbf P_x$-almost every $\omega \in \Omega$ and for all $x \in \mathcal I,$ $\tau \in \CT,$ and $s\geq 0,$
\begin{align}\label{mark}
\mbox{$\tau(\omega) < \infty$ \, implies \, $\Lambda^i_{\tau(\omega)+s}(\omega)=\Lambda^i_{\tau(\omega)}(\omega)\Lambda^i_s( \theta_{\tau (\omega)} (\omega))$}
\end{align}
or$,$ more compactly$,$ $\Lambda^i_{\tau+s}=\Lambda^i_{\tau}\cdot \Lambda^i_s \circ \theta_{\tau}$ over the event $\{\tau<\infty\}$.
\end{definition}

Definition \ref{marksur} can be intuitively understood as follows. According to Definition \ref{def:randomized_stopping_time} and Lemma \ref{lemma_properties_lambda}, $\Lambda^i_{\tau+s}$ is the probability that player $i$ concedes after time $\tau+s$ conditionally on $\CF_{\tau+s}$. The Markov restriction then states that, at time $\tau$, and conditionally on the fact that player $i$ did not concede by then, the probability that he holds fast for at least $s$ additional units of time should not depend on the trajectory prior to time $\tau$. This probability is thus given by $\Lambda^i_s \circ \theta_{\tau}$, that is, the probability induced by the randomized strategy applied to the shifted trajectory. Formula (\ref{mark}) then follows from the standard formula for conditional probabilities.

Processes satisfying (\ref{mark}) are known as \textit{strongly Markovian multiplicative functionals} of the Markov process $X$ and are studied in the literature on general Markov processes (Blumenthal and Getoor (1968)). Combining a result by Sharpe (1971) with the classical representation result of additive functionals of regular diffusions (Borodin and Salminen (2002, Chapter II, Section 4, \S23)), we can deduce the following representation result for Markovian randomized stopping times.

\begin{theorem} \label{representation}
For any Markovian randomized stopping time for player $i =1,2$ with csf $\Lambda^i : \Omega \times \RR_ +\to [0,1],$ there exists a closed set $S^i \subset \CI$ and a Radon measure\footnote{Recall that a Radon measure over an open set $U \subset \RR$ is a nonnegative Borel measure that is locally finite in the sense that every point of $U$ has a neighborhood having finite $\mu$-measure.} $\mu^i$ over $\CI \setminus S^i$ such that$,$ for each $x \in \CI,$ and for $\mathbf P_x$-almost every $\omega \in \Omega$ and each $t \geq 0,$
\begin{align} \label{repre}
\Lambda_t^i (\omega)=  \mathrm e^{- \int_{\CI} L_t^y(\omega) \, \mu^i(\mathrm dy)} \, 1_{\{\tau_{S^i}(\omega) > t\}},
\end{align}
where
\begin{align}
L_t^y \equiv \lim_{\varepsilon \downarrow 0}\, {1 \over 2 \varepsilon} \int_0^t 1_{(y- \varepsilon, y + \varepsilon)}(X_s) \sigma^2( X_s) \, \mathrm ds, \label{deflocaltime}
\end{align}
is the local time of $X$ at $(y,t),$ and
\begin{align*}
\tau_{S^i}\equiv \inf \hskip 0.5mm \{t \geq 0: X_t \in S\}
\end{align*}
is the hitting time of $S^i$ by $X$. In particular, the mapping $t \mapsto \Lambda_t^i(\omega)$ is continuous over $[0,\tau_{S^i}(\omega))$ for $\mathbf P_x$-almost every $\omega \in \Omega$.
\end{theorem}

The interpretation of (\ref{repre}) is that player $i$ concedes with probability 1 over $S^i$, and with positive but finite intensity over $\mathrm {supp} \, \mu^i$. The relation (\ref{repre}) allows us in the following to indifferently refer to a Markov strategy for player $i$ as a ccdf $\Gamma^i$, a csf $\Lambda ^i$, or a pair $(\mu^i, S^i)$; we shall use these notations interchangeably in the definition of players' payoffs. Three cases of the representation (\ref{repre}) are worth mentioning.

\paragraph{The Pure Stopping Case}

If $\mu^i \equiv 0$, then the Markov strategy $(0,S^i)$ is just the pure stopping time $\tau_{S^i}$. This is the class of Markov strategies considered by Murto (2004).

\paragraph{The Absolutely Continuous Case}

If $\mu^i \equiv g^i \cdot Leb$ is absolutely continuous with density $g^i$ with respect to Lebesgue measure, then, using the occupation time formula (Revuz and Yor (1999, Chapter VI, \S1, Corollary 1.6)), the corresponding csf writes as
\begin{align}
\Lambda_t ^i=\mathrm e^{- \int_{\CI} L_t^y g^i(y) \, \mathrm dy} \, 1_{\{\tau_{S^i}>t\}}=\mathrm e^{- \int_0^t g^i(X_s)\sigma^2(X_s)\, \mathrm ds}\,1_{ \{\tau_{S^i}>t\}}. \label{Brown}
\end{align}
Outside $S^i$, this strategy consists for player $i$ in conceding according to a Poisson process with stochastic intensity $\lambda^i(X_t) \equiv g^i(X_t)\sigma^2(X_t)$; that is, during a short time interval $[t, t+ \mathrm dt)$, he concedes with probability $1$ if $X_t \in S^i$ and with probability $\lambda^i(X_t) \, \mathrm dt $ otherwise. This is the class of Markov strategies considered by Giorgiadis, Kim, and Kwon (2022).

\paragraph{The Singular Case}

If $\mu^i \equiv a^i \delta_{x^i}$, where $a^i >0$ and $\delta_{x^i}$ is the Dirac mass at $x^i \in \CI\setminus S$, then the corresponding csf writes as
\begin{align}
\Lambda^i_t= \mathrm e ^{-a^i L^{x^i}_t}  1_{\{\tau_{S^i}>t\}}. \label{Providence}
\end{align}
In particular, the mapping $t \mapsto \Lambda_t^i(\omega)$ is singular over $[0,\tau_{S^i}(\omega))$ for $\mathbf P_x$-almost every $\omega \in \Omega$ such that the trajectory of $X$ crosses $x^i$; that is, its derivative is zero for $Leb$-almost every $t \in [0,\tau_{S^i}(\omega))$, though $\Lambda_t^i(\omega)$ is not constant as it decreases each time $X$ crosses $x^i$. To the best of our knowledge, Markov strategies with singular csf have not been considered in the literature. Yet there is no reason to discard such strategies, as they naturally emerge as limits of more familiar ones. Here are two illustrations:
\begin{itemize}

\item[(i)]

First, discretize the state space (and possibly the time space) and consider Markov strategies for player $i$ prescribing him to concede with positive intensity when the current state is $x^i$. Then, with appropriate normalizations, the natural limit of such strategies when the mesh of the discretization goes to 0 corresponds to a distribution with hazard rate proportional to the local time of the diffusion at $x^i$.

\item[(ii)]

Second, consider the Markov strategy that, outside $S^i$, consists for player $i$ in conceding according to a Poisson process with stochastic intensity $\lambda^i_\varepsilon (X_t) \equiv \frac{a^i}{2\varepsilon} \, \sigma^2(X_t) \, 1_{(x^i- \varepsilon ,x^i+\varepsilon)}$ for $a^i >0$ and some small $\varepsilon>0$. By (\ref{Brown}), the corresponding csf writes as
\begin{align*}
\Lambda^i_{\varepsilon,t}=\mathrm e^{- \int_0^t \lambda^i_\varepsilon(X_s) \, \mathrm ds} \,1_{ \{\tau_{S^i}>t\}}.
\end{align*}
From the definition \eqref{deflocaltime} of the local time $L^{x^i}_t$ of $X$ at $(x^i,t)$, we deduce that, for each $t \geq 0$, $\Lambda^i _{\varepsilon,t}$ converges $\mathbf P_x$-almost surely to $\Lambda^i_t$ in (\ref{Providence}) as $\varepsilon$ goes to 0.

\end{itemize}
Let us finally mention an important property of a Markov strategy, such as (\ref{Providence}), associated to a singular measure with an atom at $x^i$. Using the properties of the local time, one can check\footnote{This be obtained for example by adapting the method used in Lemma 15 in Peskir (2019).} that the total probability of conceding before time $t$ starting from $x^i$ is of order $\sqrt{t}$, whereas the same quantity is of order $t$ for a Markov strategy, such as (\ref{Brown}), associated to an absolutely continuous measure. As we will see in Sections 4--5, this particular singular behavior will create points of nondifferentiability in the players' equilibrium value functions.

\subsection{Markov-Perfect Equilibrium and Properties of Best Replies}

We are now ready to define our equilibrium concept and to provide some basic properties of best replies. Our first result, which we will repeatedly use in what follows, illustrates the standard fact that a player, given the behavior of his opponent, cannot improve his payoff merely by randomizing over pure strategies.

\begin{lemma} \label{usepure}
For each $x \in {\cal I}$ and for any pair of randomized stopping times with ccdf $(\Gamma^1,\Gamma^2),$
\begin{align*}
J^1 (x,  \Gamma^1, \Gamma^2) & \leq \sup_{\tau^1 \in \CT} \,J^1(x, \tau^1,\Gamma^2),
\\
J^2 (x,  \Gamma^1, \Gamma^2) & \leq \sup_{\tau^2 \in \CT} \,J^2(x, \Gamma^1, \tau^2).
\end{align*}
\end{lemma}

\vskip 3mm

This motivates the following definition.

\begin{definition}
A {\rm Markov-Perfect Equilibrium} (MPE) is a profile of Markov strategies \linebreak $((\mu^1 ,S^1),(\mu^2,S^2))$ such that$,$ for each $x \in {\cal I},$
\begin{align*}
J^1(x, (\mu^1, S^1), (\mu^2, S^2)) &= \bar J^1(x, (\mu^2, S^2)) \equiv \sup_{\tau^1 \in {\cal T}} \, J^1(x, \tau^1, (\mu^2, S^2)),
\\
J^2(x, (\mu^1, S^1), (\mu^2, S^2)) &=\bar J^2(x,  (\mu^1, S^1)) \equiv \sup_{\tau^2 \in {\cal T}} \,J^2(x, (\mu^1, S^1), \tau^2).
\end{align*}
That is$,$ for each $i=1,2,$ $(\mu^i,S^i)$ is a {\rm perfect best reply} (pbr) for player $i$ to $(\mu^j,S^j),$ and $\bar J^i(\cdot, (\mu^j, S^j))$ is player $i$'s {\rm best-reply value function} (brvf) to $(\mu^j, S^j)$.
\end{definition}

When no confusion can arise as to the strategy of player $j$, we write $\bar J^i$ instead of $\bar J^i (\cdot, (\mu^j, S^j))$. The next proposition  provides useful general properties of pbr and brvf, and is key to establish our main results.

\begin{proposition} \label{geneprop}
If $(\mu^i,S^i)$ is a pbr to $(\mu^j,S^j)$ with associated brvf $\bar J^i,$ then $V_{R^i}\leq \bar J^i \leq G^i$. Furthermore$,$
\begin{itemize}

\item[(i)]

$S^1\cap S^2\cap (\alpha^i,\beta)=\emptyset;$

\item[(ii)]

$S^i \subset C^i\equiv \{ x \in \CI : \bar {J}^i(x)=R^i(x)\};$

\item[(iii)]

$\mathrm {supp}\,\mu^i \subset C^i \cup S^j;$

\item[(iv)]

$\mathrm {supp}\,\mu^i\cup S^i \subset (\alpha,x_{R^i}];$

\item[(v)]

$(0,S^i)$ is also a pbr to $(\mu^j,S^j)$ and more generally $(\tilde \mu^i,S^i)$ is a pbr to $(\mu^j,S^j)$  for any $\tilde\mu^i$ such that $\mathrm {supp} \, \tilde\mu^i\subset C^i \cup S^j;$

\end{itemize}
\end{proposition}

Property (i) intuitively states that player $i$ should never concede when market conditions $x$ are such that player $j$ concedes with probability 1 and player $i$'s payoff from conceding is strictly less than the payoff from letting player $j$ concede, that is, $x \in S^j$ and $G^i(x) > V_{R^i}(x)$. Property (ii) simply expresses the fact that player $i$'s brvf coincides with $R^i$ over the portion $S^i$ of the state space over which he concedes with probability 1. Property (iii) states that player $i$'s payoff is $R^i$ when he concedes with positive intensity outside of player $j$'s stopping region $S^j$. Property (iv) reflects that player $i$ should never concede when market conditions are above the optimal threshold $x_{R^i}$ for his stand-alone exit problem; intuitively, this is because waiting for $X$ to drop down to $x_{R^i}$ before conceding is player $i$'s optimal strategy even in the worst-case scenario in which player $j$ is stubborn, that is, $(\mu^j,S^j) = (0, \emptyset)$. Finally, property (v) states that, when conceding with positive intensity outside of $S^i$, player $i$ should be indifferent between holding fast and conceding.

\paragraph{Remark}

Some authors (see, for instance, Murto (2004)) include, as a refinement in the definition of an MPE, the requirement that $(\alpha, \alpha^i] \subset S^i$ for all $i$. The rationale for this assumption is that, because $G^i=V_{R^i}=R^i$ over $(\alpha, \alpha^i]$, holding fast further below $\alpha^i$ would be weakly dominated for player $i$ by conceding with probability 1 over this interval. For instance, being stubborn is a best reply for player $i$ over $(\alpha,\alpha^i)$ only if player $j$ concedes with probability 1 over this interval, except perhaps over a set of Lebesgue measure 0. This behavior is not per se inconsistent with an MPE, but it is not consistent with trembling-hand perfection in the spirit of Selten (1975), see Ghemawat and Nalebuff (1985) for a discussion of a similar point in a deterministic model. Hereafter, we do not systematically impose this refinement, especially in Section 5 where this allows to simplify notation; however, we will indicate which MPEs can be modified so as to satisfy it.

\vskip 3mm

We close this section with an important global regularity result.

\begin{proposition}\label{lem:continuity}
If $((\mu^1,S^1),(\mu^2,S^2))$ is an MPE$,$ then$,$ for each $i=1,2$, player $i$'s brvf $\bar J^i$ is continuous over $\CI$.
\end{proposition}

\section{MPEs with Singular Strategies: Heuristics}\label{sec:benchmark}

We first recall within our general framework two standard MPEs, respectively in pure and mixed strategies, that have been emphasized in the literature. Based on these examples and on our representation theorem for Markovian randomized stopping times, we next describe a novel type of MPE involving a singular strategy for one of the two players. Our heuristic presentation leads to a variational system that turns out to fully characterize the candidate equilibrium. We use the running example of Section \ref{running} to illustrate our findings.

\subsection{A Pure-Strategy MPE} \label{08:13}

We say that player 1 is \textit{as least as enduring} as player 2 if $\alpha^1 \leq\alpha^2$ and $x_{R^1}\leq x_{R^2}$; intuitively, player 1 is at least as willing to hold fast as player 2. Suppose then that player 1 threatens to hold fast maximally and concede only at $\tau^1 = \inf \hskip 0.5mm \{t \geq 0 : X_t \leq \alpha^1\}$. Then, because $\alpha^1 \leq \alpha^2$, we have $G^2(X_{\tau^1}) =R^2(X_{\tau^1})$ by definition of $\alpha^2$. In light of (\ref{core})--(\ref{sa}), this implies that, for all $x \in \mathcal I$ and $\tau^2 \in \CT$,
\begin{align*}
J^2(x,\tau^1,\tau^2) = \mathbf E _x\! \left[\mathrm e^{-r \tau^1 \wedge \tau^2}R^2 (X_{\tau^1 \wedge \tau^2}) \right] \! \leq V_{R^2}(x).
\end{align*}
Thus a pbr for player 2 to $\tau^1$ is to concede at $\tau^2 = \inf \hskip 0.5mm \{t \geq 0 : X_t \leq x_{R^2}\}$. As for player 1, if player 2 concedes at $\tau^2$, then, for each $x\in \mathcal I$,
\begin{align*}
\mathbf E_x \! \left[\mathrm e^{-r \tau^2}G^1(X_{\tau^2}) \right] \geq R^1(x).
\end{align*}
For $x \leq x_{R^2}$, this follows from the fact that $G^1 (x) \geq R^1(x)$ by A6, with a strict inequality if $x > \alpha^1$. For $x > x_{R^2}$, this follows from A6 again along with the fact that the process $(\mathrm e^{-rt}V_{R^1}(X_t))_{t\geq 0}$ is a martingale up to $\tau_{x_{R^1}}$, the hitting time of $ x_{R^1}$, which is no less than $\tau^2$ because $x_{R^1} \leq x_{R^2}$ by assumption. Thus a pbr for player 1 to $\tau^2$ is to concede at $\tau^1$. This implies the following result, which has many counterparts in the literature (see, for instance, Ghemawat and Nalebuff (1985), D\'ecamps and Mariotti (2004), Murto (2004), Giorgiadis, Kim, and Kwon (2022)).

\begin{proposition} \label{pure}
If player 1 is at least as enduring as player 2$,$ then $((0, (\alpha, \alpha^1]), (0,(\alpha, x_{R^2}]) )$ is a pure-strategy MPE.
\end{proposition}

In the case where the asymmetry between the players is small, $((0, (\alpha, x_{R^1}]), (0, \emptyset) )$ is also an MPE in which the more enduring player 1 follows his stand-alone optimal strategy because the less enduring player 2 is stubborn (Giorgiadis, Kim, and Kwon (2022)). However, this MPE does not satisfy Murto's (2004) trembling-hand-perfection refinement, because, for $x \in (\alpha^1,\alpha^2)$, player 2's strategy is no longer a best response when player 1 does not concede with probability 1 in any small enough neighborhood of $x$. Nevertheless, Murto (2004) shows that, when we allow player 1's stopping set $S^1$ to exhibit a gap, there may exist an MPE satisfying this refinement in which, when $x > x_{R^1}$, player 1 exits first when $X$ reaches $x_{R^1}$.

\subsection{A Mixed-Strategy MPE in the Symmetric Case}

Suppose now that players are symmetric, in the weak sense that they are as enduring as each other, $\alpha^1 =\alpha^2\equiv \alpha^*$ and $x_{R^1}= x_{R^2}\equiv x^*$. This is of course the case when the players have identical payoff functions, $R^1=R^2$ and $G^1 = G^2$. The following result, which restates in our framework earlier results in the literature (Steg (2015), Georgiadis, Kim, and Kwon (2022))\footnote{ A related construction also appears in Kwon and Palczewski (2022). There, a symmetric Bayesian equilibrium is constructed in a model with asymmetric information and a continuum of types. The pure strategies, seen as randomized strategies assimilating the types as randomization devices, use absolutely continuous intensities depending on $X$ and on an auxiliary belief process.}, characterizes a mixed-strategy MPE in which the players concede with absolutely continuous intensities over the interval $(\alpha^*, x^*]$.

\begin{proposition} \label{sym}
If the players are as enduring as each other$,$ then the strategy profile
\begin{align*}
((\lambda^1(x)\sigma^{-2}(x)\, \mathrm dx, (\alpha, \alpha^*]), (\lambda^2(x)\sigma^{-2}(x)\, \mathrm dx, (\alpha, \alpha^*]))
\end{align*}
defined$,$ for each $i=1,2,$ by
\begin{align}\label{def_density}
\lambda^i(x)  \equiv \frac{rR^j(x)  - {\cal L} R^j(x)}{G^j(x) - R^j(x)} \, 1_{\{\alpha^*<x \leq x^*\}},
\end{align}
is a mixed-strategy MPE.
\end{proposition}

We know from Proposition \ref{geneprop}(iv) that $\mathrm {supp}\,\mu^i \subset (\alpha^*, x^*]$. Following Theorem 5.1 in Steg (2015), the MPE constructed in Proposition \ref{sym}  is such that each player stops with an intensity function $\lambda^i$ with support $(\alpha^*, x^*]$. This intensity is constructed so that, at each point of this interval, each player is indifferent between holding fast and conceding. In equilibrium, the value function of each player coincides with the value function of his stand-alone exit problem (\ref{sa}). Thus, in expectation, the war of attrition yields no benefit to either player.

\subsection{A Singular Mixed-Strategy MPE} \label{EXAMPLE}

When players are asymmetric but there is no uncertainty about their future payoffs, the war of attrition may admit mixed-strategy equilibria in which players strategies are described, over some interval of exit times, by absolutely continuous distributions (see, for instance, Hendricks, Weiss, and Wilson (1988)); this is certainly the case in the limiting case of our model where $\mu = \sigma = 0$, so that market conditions are constant. This result has no counterpart under Brownian uncertainty. Indeed, Georgiadis, Kim, and Kwon (2022) have shown that the construction of Proposition \ref{sym} does not extend to the case of asymmetric players: specifically, when players are not as enduring as each other, there exists no mixed- strategy MPE in which the players concede with absolutely continuous intensities. For all that, it would be incorrect to conclude that only pure-strategy MPEs exist, and thus that attrition cannot take place in an MPE of our model. This section argues for this claim by describing an MPE involving a singular strategy for one of the two players. For the sake of simplicity, the analysis below remains at a heuristical level. A full justification of our arguments is provided in Section 5.

From now on, assume as in Section \ref{08:13} that player 1 is at least as enduring as player 2, that is, $\alpha^1 \leq \alpha^2$ and $x_{R^1} \leq x_{R^2}$. Consider then the following equation in $x$:
\begin{align} \label{barb1}
R^1(x_{R^1}) = \frac{ \phi (x_{R^1})}{\phi(x)}\,G^1(x).
\end{align}
We show in Appendix B that (\ref{barb1}) admits a unique solution $\underline x^2 \in(\alpha^1, x_{R^1})$. In words, the threshold $\underline x^2$ is such that, if player 2 threatens to concede only at $\tau^2 = \inf \hskip 0.5mm \{t \geq 0 : X_t \leq \underline x^2\}$, then, at $x_{R^1}$, player 1 is indifferent between conceding and obtaining $R^1(x_{R^1})$ immediately and waiting for player 2 to concede at $\tau^2$ and obtaining $G^1(\underline x^2)$ only then. We claim that, if $x_{R^2}$ is close enough to $x_{R^1}$, then there exists an MPE in which player 1 randomizes between holding fast and conceding at $x_{R^1}$ and player 2 concedes only at $\tau^2$. Using the representation for randomized stopping times provided in Theorem \ref{representation}, this amounts to the existence of a constant $a^1>0$ such that the profile of strategy $((a^1 \delta_{x_{R^1}}, (\alpha, \alpha^1)), (0, (\alpha, \underline x^2))$ is an MPE in which player 1 concedes, with positive but finite intensity, only at $x_{R^1}$.

\subsubsection{Necessary Conditions} \label{NCHEUR}

To establish this claim, we first assume that such an MPE exists, and we derive necessary conditions for the brvf $\bar J^1$ and $\bar J^2$. An obvious preliminary observation is that $\bar J^1 \geq R^1$ and $\bar J^2 \geq R^2$, because, when current market conditions are $x$, every player $i$ can always guarantee himself the payoff $R^i(x)$ by exiting the market immediately.

\paragraph{Player 1}

Player 1, whose strategy involves randomization at $x_{R^1}$, should be indifferent at $x_{R_1}$ between conceding and holding fast until $\tau^2$. This implies that his brfv $\bar J^1$ must be $\mathcal C^2$ over $(\underline x^2, \beta)$, with $\bar J^1 (x_{R^1}) = R^1(x_{R^1})$ (value-matching). Because $\bar J^1 \geq R^1$, it follows in turn that $\bar J^{1 \prime} (x_{R^1}) = R^{1 \prime}(x_{R^1})$ as well. Moreover, by standard dynamic-programming arguments, $\bar J^1$ must satisfy the ODE $\mathcal L \bar J_1 -r\bar J_1= 0$ over $(\underline x^2, \beta)$ (see, for instance, Dixit and Pindyck (1994)). This leads to
\begin{align}
\bar J^1 (x) = \frac{\phi(x) }{\phi (x_{R^1})}\,R^1(x_{R^1}), \quad x \in (\underline x^2, \beta). \label{meuble}
\end{align}
In particular, $\bar J^1 = V_{R^1}$ over $[x_{R^1}, \beta)$: player 1 does not benefit from the war of attrition over $[x_{R^1}, \beta)$. By contrast, $\bar J^1 > V_{R^1}$ over $[\underline x^2, x_{R^1})$, reflecting that player 1 can hope that player 2 may concede at $\underline x^2$ before he himself concedes at $x_{R^1}$.

\paragraph{Player 2}

Player 2 plays a pure strategy and hopes to benefit from player 1 conceding at $x_{R^1}$. We guess that $\bar J^2$ is $\mathcal C^2$ over $(\underline x^2, x_{R^1}) \cup (x_{R^1},\beta)$, with $\bar J^2(\underline x^2)= R^2(\underline x^2)$ (value-matching) and $\bar J^{2\prime}(\underline x^2)= R^{2 \prime}(\underline x^2)$ (smooth pasting), and that it satisfies the ODE ${\cal L}\bar J^2 - r\bar J^2 =0 $ over that region. There remains to characterize the behavior of $\bar J^2$ at $x_{R^1}$. Because player 1 randomizes at $x_{R^1}$ between holding fast and conceding, we expect that $G^2( x_{R^1}) > \bar J^2( x_{R^1}) >R^2(x_{R^1})$. This, along with the properties of the local time highlighted in Section \ref{MRST}, implies that $\bar J^2$ is not differentiable at $ x_{R^1}$. Indeed, starting from $x_{R^1}$, player 1 concedes in a small time interval of length $\mathrm dt$ with probability $\mathbf E_{x_{R^1}}[\Gamma_{\mathrm dt}] = a^1 c \sqrt{\mathrm dt} + o (\sqrt{\mathrm dt})$, where $\Gamma_{\mathrm dt} = 1 - \mathrm e^{-a^1 L_{\mathrm dt}^{ x_{R^1}}}$ and $c$ is a positive constant. If player 1 concedes, then player 2 benefits from the follower payoff $G^2( x_{R^1})$, while if player 1 holds fast, player 2 achieves the value $\bar J^2(X_{\mathrm dt})$. Thus
\begin{align} \label{heur2}
\bar J^2(x_{R^1}) =a^1 c \sqrt{\mathrm dt} \,G^2( x_{R^1}) +(1 - a^1c \sqrt{\mathrm dt})\, \mathbf E_{ x_{R^1}} [\mathrm e^{-r \mathrm dt} \bar J^2(X_{\mathrm \mathrm dt}) ] + o (\sqrt{\mathrm dt}).
\end{align}
Now, suppose, by way of contradiction, that $J^2$ is $\mathcal C^2$ in a neighborhood of $x_{R^1}$. Then, from It\^o's formula,
\begin{align}
\mathbf E_{ x_{R^1}} [\mathrm e^{-r \mathrm dt} \bar J^2(X_{\mathrm dt}) ] = \bar J^2(x_{R^1}) +  ({\cal L} \bar J^2 - r \bar J^2)(x_{R^1})\, \mathrm dt + o(\mathrm dt). \label{EpBa}
\end{align}
Plugging (\ref{EpBa}) into (\ref{heur2}) yields $a^1 c[G^2( x_{R^1}) - \bar J^2( x_{R^1})] \sqrt{ dt}+ o (\sqrt{\mathrm dt})=0$, a contradiction as $G^2( x_{R^1})> \bar J^2( x_{R^1})$ and $a^1$ and $c$ are positive constants. This is an indication that $\bar J^2$ is not differentiable at $x_{R^1}$; let us denote by $ \Delta \bar J^{2 \prime}(x_{R^1}) \equiv \bar J^{2\prime +}(x_{R^1})  - \bar J^{2 \prime -}(x_{R^1})$ the corresponding derivative jump. From the It\^o--Tanaka--Meyer formula, which generalizes It\^o's formula to functions, such as $\bar J^2$, that can be written as the difference of two convex functions (Karatzas and Shreve (1991, Theorem 3.7.1 and Problem 3.6.24)), we have
\begin{align}
\mathbf E_{ x_{R^1}} [\mathrm e^{-r\mathrm dt} \bar J^2(X_{\mathrm dt}) ] &= \bar J^2( x_{R^1}) + \mathbf E_{ x_{R^1}} \! \left[\int_0^t \mathrm e^{-rs} ({\cal L} \bar J^2 - r \bar J^2) (X_s) \, \mathrm ds \right. \notag
\\
& \hskip 3.55cm  + \left. \int_0^t \mathrm e^{-rs} \bar J^{2 \prime -}(X_s) \sigma(X_s) \, \mathrm dW_s + \frac{1}{2}\, \Delta \bar J^{2\prime}(x_{R^1}) L_{\mathrm dt}^{x_{R^1}}\right] \notag
\\
& = \bar J^2( x_{R^1}) + \frac{1}{2}\, \Delta \bar J^{2\prime}(x_{R^1}) c\sqrt{dt} + o (\sqrt{\mathrm dt}), \label{EpBaITM}
\end{align}
where the second equality follows from the fact that ${\cal L}\bar J^2 - r\bar J^2 =0 $ over $(\underline x^2, x_{R^1}) \cup (x_{R^1},\beta)$ and from the properties of local time. Plugging (\ref{EpBaITM}) into (\ref{heur2}) yields
\begin{align}
a ^1 [G^2( x_{R^1}) - \bar J^2( x_{R^1})]  + \frac{1}{2}\,\Delta \bar J^{2\prime}(x_{R^1}) =0. \label{Choli}
\end{align}
Notice from $G^2( x_{R^1}) > J^2( x_{R^1})$ and (\ref{Choli}) that $\Delta \bar J^{2\prime}(x_{R^1}) <0$. Intuitively, player 2 gets more and more optimistic as $X $ approaches $x_{R^1}$, but is disappointed if $X$ crosses $x_{R^1}$ but player 1 does not concede at $x_{R^1} $.

\paragraph{The Variational System}

Our discussion so far leads to the following variational system: find a constant $a>0$, and two functions $w^1 \in {\cal C}^0((\alpha, \beta))\cap {\cal C}^2 ((\alpha, \beta) \setminus \{  \underline x^1 \})$ and $w^2 \in {\cal C}^0((\alpha, \beta)) \cap {\cal C}^2 ((\alpha,\beta) \setminus \{\underline x^1,  x_{R^1} \})$ such that
\begin{align}
w^1 &\geq R^1 \text{ over } (\alpha, \beta), \label{Gonsolin}
\\
\CL w^1 - rw^1  &= 0 \text{ over } (\underline x^2, \beta), \label{var11e}
\\
w^1&= G^1\text{ over } (\alpha, \underline x^2], \label{var12e}
\\
w^1(x_{R^1})&= R^1(x_{R^1}),  \label{var13e}
\\
w^1(\beta^-)&= 0, \label{var15e}
\\
\notag
\\
w^2 &\geq R^2 \text{ over } (\alpha, \beta), \label{Ermet}
\\
\CL w^2 - rw^2 &= 0 \text{ over } (\underline x^2, x_{R^1}) \cup (x_{R^1}, \beta), \label{var1'e}
\\
w^2(\underline x^2 ) &= R^2(\underline x^2), \label{var2'e}
\\
w^{2\prime}(\underline x^2)&= R^{2 \prime}(\underline x^2), \label{var3'e}
\\
a[ G^2(x_{R^1}) - w^2(x_{R^1})] + \frac{1}{2} \,\Delta {w^{2\prime}}(x_{R^1})&=0,  \label{var4'e}
\\
w^2(\beta^-) &=0. \label{var5'e}
\end{align}

\subsubsection{Sufficient Conditions}

It is an implication of our main characterization result, Theorem 4, that, if $(a^1, \bar J^1,\bar J^2)$ is a solution to the variational system (\ref{Gonsolin})--(\ref{var5'e}), then $\bar J^1$ is the brfv to $(0,(\alpha, \underline x^2))$ and $w^2$ is the brfv to $(a^1 \delta_{x _{R^1}}, (\alpha, \alpha^1))$, so that, according to the construction in Section \ref{NCHEUR}, $((a^1 \delta_{x_{R^1}}, (\alpha, \alpha^1)), (0, (\alpha, \underline x^2))$ is an MPE. As for $\bar J^1$, we have already seen that (\ref{Gonsolin})--(\ref{var11e}) and (\ref{var13e}) pin down a unique solution, given by (\ref{meuble}), which satisfies $\bar J^1 (\underline x^2) = G^1 (\underline x^2)$ by definition of $\underline x^2$. As for $\bar J^2$, the analysis is a bit more delicate due to the presence of the derivative jump $\Delta {\bar J^{2\prime}}(x_{R^1})$ at $x_{R^1}$, which, by (\ref{var4'e}), is pinned down by the intensity $a^1$ with which player 1 exits at $x_{R^1}$. In our running example, it can be shown that, as long as the asymmetry between the players is small, assuming that $b>0$ and $m$ is sufficiently large, one can indeed find a positive value for $a^1$ such that (\ref{Ermet})--(\ref{var5'e}) holds. The following result then holds\footnote{A numerical study that we performed suggests that, for firms' liquidation values $l^1 \leq l^2$  close to each other, the variational systems (\ref{Gonsolin})-(\ref{var15e}) and (\ref{Ermet})-(\ref{var5'e}) admit  a solution whatever the parameter values of the model if $b>0$, and for $m$  in some compact interval $[1,C]$ if $b<0$ for some constant $C$ which increases with $\sigma$.}. 

\begin{proposition} \label{armchair}
In the running example$,$ if the firms' liquidation values $l^1 \leq l^2$ are close to each other$,$ $b>0$ and $m$ is sufficiently large, then there exists a mixed-strategy MPE $((a^1 \delta_{x_{R^1}}, (\alpha, \alpha^1)), (0, (\alpha, \underline x^2))$ in which the more enduring firm 1 randomizes between holding fast and conceding at $x_{R^1}$ while the less enduring firm 2 exits with probability 1 as soon as market conditions fall below $\underline x^2 < x_{R^1}$.
\end{proposition}

The MPE constructed in Proposition \ref{armchair} differs from the pure-strategy MPE of Proposition \ref{pure} in that, for $x \geq x_{R^1}$, it is the more enduring enduring player 1 who does not benefit from the war of attrition, in the sense that $\bar J^1 = V_{R^1}$ over $[x_{R^1},\beta)$; by contrast, $\bar J^2 > V_{R^2}$ over this portion of the state space. The reason is that player 2 adopts a tougher stance by threatening to exit the market only at $\underline x^2 < x_{R^1}< x_{R^2}$, which makes player 1 indifferent between holding fast and conceding at $x_{R^1}$. By construction, this MPE satisfies the requirement that $(\alpha , \alpha^i] \in S^i$ for every player $i$, as in Murto (2004).

It should also be noted that, in this singular mixed-strategy MPE, we have $\max \hskip 0.5mm S^1 \vee \max \hskip 0.5mm S^2 = \underline x^2 < x_{R^1} <x_{R^2}$. This contrasts with pure-strategy MPEs, in which one always have $\max \hskip 0.5mm S^1 \vee \max \hskip 0.5mm S^2 \in \{ x_{R^1} ,x_{R^2}\}$. Thus mixing by player 1 delays the time at which a firm must necessarily exit the market. In particular, the difference with the pure-strategy MPE characterized by Murto (2004), in which the stopping set $S^1$ of player 1 exhibits a gap below $\max \hskip 0.5mm S^1 = x_{R^1}$ and player 1 is the first to exit the market at $x_{R^1}$ when $x > x_{R^1}$, is that player 1 does not exit with probability 1 at $x_{R^1}$. This leads to a richer dynamics, whereby, on the equilibrium path, every player $i$ can alternate between being in a dominated position (with a value close to $V_{R^i}$) or in a dominant position (with a value significantly above $V_{R^i}$); specifically, player 1 is in a dominant position when $X$ is close to $\underline x^2$, while player 2 is in a dominant position when $X$ is close to $x_{R^1}$. As we show in Theorem 3, this alternation phenomenon is a robust feature of any mixed-strategy MPE.

Figure 1 illustrates this point in the running example by plotting the firms' market-value functions $F^i \equiv \bar J^i +E$ in the singular mixed-strategy MPE constructed in Proposition \ref{armchair}. Notice that firm 1's value function coincides with its monopolist value function $V^1_m$ over $(\alpha, \underline x_2]$, because player 2 exits the market with probability 1 at any point of this interval. It can also be checked that $F^{2 \prime - }(x_{R^1}) >0 > F^{2 \prime +} (x_{R^1})$, reflecting that player 2's market value reaches a local maximum when $X_t = x_{R^1}$. Finally, a testable implication of this MPE, which is also apparent from Figure 1, is that, along any path of $X$, the firms' market values move in opposite directions as long as none of them exits the market and market conditions do not wander too much above $x_{R^1}$. Again, this negative comovement of firms' market values is a robust feature of any mixed-strategy MPE.
\noindent \setcounter{figure}{0}
\begin{figure}
\noindent {\normalsize \unitlength=1.00mm \special{em:linewidth
0.4pt} \linethickness{0.4pt}
\begin{picture}(165.00,80)

\put(10,10){\vector(1,0){151}}
\put(161,6.25){\makebox(0,0)[cc]{$x$}}

\put(10,10){\vector(0,1){80}}
\put(2,90){\makebox(0,0)[cc]{Values}}

\bezier{1000}(10,10)(85.5,22)(161,34)
\put(125,22){\makebox(0,0)[cc]{$\displaystyle  {x \over r-b}$}}

\color{purple}
\bezier{1000}(10,10)(30,49) (50,88)
\put(53,80){\makebox(0,0)[cc]{$\displaystyle  {m x \over r-b}$}}

\color{blue}
\bezier{141}(10,30)(85.5,30)(161,30)
\bezier{1000}(40,70)(50,16)(161,35)
\bezier{1000}(12,30)(18,30) (40,70)
\bezier{1000}(10,30)(11,30)(12,30)
\bezier{50}(110,30)(110,20)(110,10)
\bezier{50}(40,40)(40,55)(40,70)
\bezier{50}(12,10)(12,20)(12,30)
\put(51,60){\makebox(0,0)[cc]{$F^1(x)$}}
\put(5.8,30.6){\makebox(0,0)[cc]{$l^1$}}
\put(110,6.25){\makebox(0,0)[cc]{$x_{R^1}$}}
\put(12,6.25){\makebox(0,0)[cc]{$\alpha^1$}}

\color{red}
\bezier{141}(10,40)(85.5,40)(161,40)
\bezier{1000}(10,40)(25,40)(40,40)
\bezier{1000}(40,40)(90,40) (110,65)
\bezier{1000}(110,65)(130,43) (161,45)
\bezier{50}(110,30)(110,47.5)(110,65)
\bezier{50}(40,10)(40,25)(40,40)
\put(5.8,40.6){\makebox(0,0)[cc]{$l^2$}}
\put(125,60){\makebox(0,0)[cc]{$F^2(x)$}}
\put(40,6.25){\makebox(0,0)[cc]{$\underline x^2$}}

\end{picture}}
\vskip -2mm
\caption{The value of never exiting the market in a duopoly (in black), the value of never exiting the market in a monopoly (in purple), and the more enduring firm 1's value (in blue) and the less ensuring firm 2's value (in red) in the singular mixed-strategy MPE.}
\end{figure}
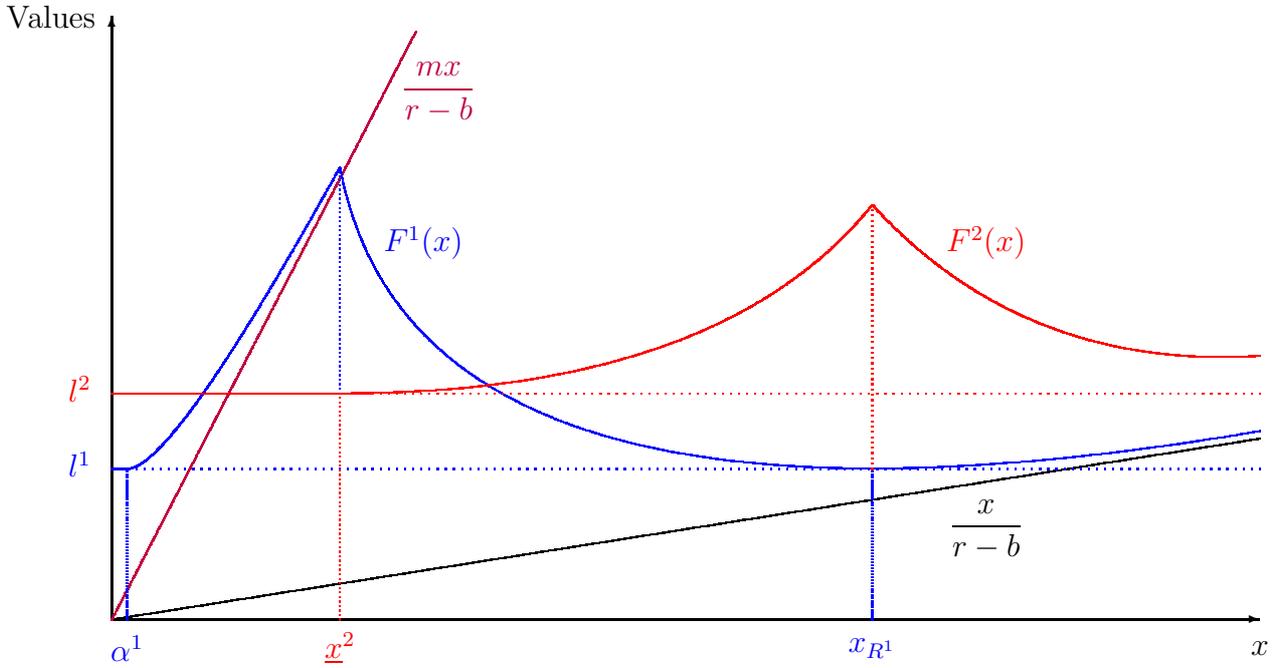

Our analysis also has novel welfare implications. As already noted, because firm 1 does not exit with probability 1 at $x_{R^1}$, firm 2's market value at $x_{R^1}$, $F^2(x_{R^1})$, must be less than its value as a monopolist, $V^2_m(x_{R^1})$. Because firm 1's market value satisfies $F^1 = V_{R^1} + E$ over $[x_{R^1}, \infty)$ and thus coincides with his stand-alone market value, this implies that, for $x > x_{R^1}$, the mixed-strategy MPE characterized in Proposition \ref{armchair} is ex-ante Pareto dominated by the pure-strategy MPE $((0, (\alpha, x_{R^1}]), (0, \emptyset) )$. Thus, even when firms have asymmetric liquidation values, wasteful attrition takes place with positive probability on the equilibrium path, in contrast with the conclusion of Georgiadis, Kim, and Kwon (2020).

\section{Main Results} \label{mainresults}

This section presents our main results, which extend Proposition \ref{armchair} to our general model. We first provide a necessary condition for mixed-strategy MPEs in the case of asymmetric players. We then characterize such equilibria by a variational system satisfied by the two players' value functions.

\subsection{The Alternating Structure of Mixed-Strategy MPEs}

The proofs of our main results require an additional regularity assumption, which we will throughout maintain in what follows.

\begin{itemize}

\item[\bf A8]

The functions $b$, $\sigma$, and $R^{i\prime\prime}$ are locally Lipschitz.

\end{itemize}

By convention, let us set $\max \hskip 0.5mm \emptyset \equiv \alpha$ and, for any MPE $((\mu^1, S^1), (\mu^2, S^2))$, let $s^i \equiv \max \hskip 0.5mm S^i$. The following theorem then holds.

\begin{theorem} \label{CNmixed}
Suppose that $x_{R^1} \neq  x_{R^2}$. Then$,$ for any mixed-strategy MPE $((\mu^1, S^1), (\mu^2, S^2)),$ the restrictions of the intensity measures $\mu^1$ and $\mu^2$ to $(s^1 \vee s^2, \beta)$ are purely atomic.
\end{theorem}

Theorem \ref{CNmixed} confirms the basic insight of Georgiadis, Kim, and Kwon (2022), according to which there exists no mixed-strategy MPE with absolutely continuous intensity measures when $x_{R^1} \neq  x_{R^2}$. Thus, if a mixed-strategy MPE exists at all, it must feature intensity measures that are singular with respect to Lebesgue measure. The additional information provided by Theorem \ref{CNmixed} is that these measures must be discrete, which, for instance, rules out intensity measures with Cantor-set supports. The proof can be sketched as follows. Let us consider a mixed-strategy MPE $((\mu^1, S^1), (\mu^2, S^2))$, supposing one exists. Property (v) of Proposition \ref{geneprop} and dynamic-programming arguments imply that the brvf $\bar J^i$ of every player $i$ satisfies the ODE $\CL u - ru =0$ over any interval where player $j$ does not concede. It also follows from Proposition \ref{geneprop} that $\bar J^i \geq V_{R^i}$ and that $\bar J^i(q^i) = R^i(q^i)$ for all $q^i \in \mathrm{supp} \, \mu^i$. These observations imply that any point $q^i>s^1 \vee s^2$ in the support of $\mu^i$ is isolated in this set. Thus the only potential accumulation point of the support of $\mu^i$ is $s^1 \vee s^2$; it turns out that in this case $s^1 \vee s^2 = \alpha$, so that $S^1=S^2= \emptyset$. (This is where A8 is needed.) In turn, over an open interval around any point $q^i>s^1 \vee s^2$ in the support of $\mu^i$, $\bar J^i$ coincides with the solution $T^i_{q^i}$ to the ODE ${\cal L } u - r u = 0$ that is tangent to $R^i$ at $q^i$. At any such point $q^i$, player $i$ is indifferent between holding fast and conceding.

Because the brvf of both players must satisfy these properties, a lot of structure is thereby induced on the supports of $\mu^1$ and $\mu^2$. First, it is easy to see that, because $\bar J^i \geq V_{R^i}$, we must have $\max \hskip 0.5mm \mathrm {supp} \, \mu^i \leq x_{R^i}$ for every player $i$. To fix ideas, let us suppose that $q^1_1 \equiv \max \hskip 0.5mm \mathrm {supp} \, \mu^1 \geq \max \hskip 0.5mm \mathrm {supp} \, \mu^2 \equiv q^2_1$. We can then show that in fact $x_{R^1} = q^1_1 > q^2_1$. It also follows that, for each $i=1,2$, and for any two consecutive points $q^i_n > q^i_{n+1} > s^1 \vee s^2$ in the support of $\mu^i$, there must exists a single point $q^j \in (q^i_{n+1},q^i_n)$ in the support of $\mu^j$ at which player $j$ is indifferent between conceding or holding fast; as discussed above, this implies that $\bar J^j = T^j_{q^{j }}$ over $(q^i_{n+1},q^i_n)$. We thus obtain two decreasing sequences of randomization thresholds $(q^{1}_n)_{n=1}^{N^1}$ and $(q^{2}_n)_{n=1}^{N^2}$, with either $N^1= N^2 = \infty$ or $0 \leq N^1-N^2 \leq 1$, which are intertwined in the sense that $q^1_1 > q^2_1 >q^1_2 > q^2_2 >\ldots$ as long as these thresholds are defined, and which satisfy $\lim _{n \to \infty} q^{1}_n=\lim _{n \to \infty} q^{2}_n = \alpha$ if $N^1= N^2 = \infty$. These two sequences characterize the restrictions of $\mu^1$ and $\mu^2$ to $(s^1 \vee s^2, \beta)$. Thus any mixed-strategy MPE must fall into one of three categories.

\begin{corollary} \label{coromix}
Suppose that $x_{R^1} \neq  x_{R^2},$ and let $((\mu^1, S^1), (\mu^2, S^2))$ be a mixed-strategy MPE such that$,$ for every player $i,$ $\mathrm {supp} \, \mu^i \cap (s^1 \vee s^2,\beta) = \{q^{i}_n : n=1,\ldots ,N^i\}$ for intertwined decreasing sequences of randomization thresholds $(q^1_n)_ {n=1}^{N^1}$ and $(q^2_n)_{n=1}^{N^2}$ satisfying$,$ with no loss of generality$,$ $q^1_1 > q^2_1$. Then $q^{ 1}_1 = x_{R^1}$ and
\begin{enumerate}

\item

if $N^1=N^2\equiv N \in  \mathbb N \setminus \{0\},$ then $q^{ 1}_{N}> q^{ 2}_N> s^1 >s^2;$

\item

if $N^1=N^2+1\equiv N \in  \mathbb N \setminus \{0\},$ then $q^2_{N-1} > q^{ 1}_{N}> s^2 >s^1, $ with $q^2_0 \equiv \beta$ by convention$;$

\item

if $N^1=N^2 = \infty,$ then $\lim _{n \to \infty} q^{1}_n=\lim _{n \to \infty} q^{2}_n = s^1=s^2= \alpha,$ so that $S^1=S^2=\emptyset$.
\end{enumerate}
\end{corollary}

In an MPE of type 1, player 1 exits the market with probability 1 at $s^1$, and player 2 has the lowest randomization threshold. In an MPE of type 2, player 1 has the lowest randomization threshold, and player 2 exits the market with probability 1 at $s^2$---the example of Section \ref{EXAMPLE} is a case in point, with $N^1=1$ and $N^2=0$. In an MPE of type 3, neither player exits the market with probability 1 at any point of the state space, and players keep randomizing all the way down to $\alpha$.  It should be noted that an equilibrium of type 3 can exist only if $\alpha^1 = \alpha^2 = \alpha$; indeed, every player $i$ such that $\alpha^i > \alpha$ would not be willing to delay exiting the market over $(\alpha, \alpha^i)$ if his opponent were to do the same. 

The upshot from Theorem \ref{CNmixed} and Corollary \ref{coromix} is that, when players have different stand- alone optimal exit thresholds, alternation is a robust feature of any mixed-strategy MPE, which generalizes the insights from Section \ref{EXAMPLE}. In the attrition region, players randomize between conceding and holding fast at isolated thresholds. In an equilibrium of type 1 and type 2, this process may persist until one player eventually reaches his stopping region and exits the market with probability 1. By contrast, in an equilibrium of type 3, exit must take place at some randomization threshold.

Corollary \ref{coromix} fully characterizes equilibrium outcomes for an MPE of type 3, because any market conditions in $\mathcal I$ can be reached with positive probability from any initial market conditions $x \in \mathcal I$. The same holds true for MPEs of types 1 and 2, provided $x > x_{R^1}$ (with the convention that $q_1^1 > q^2_1$). Indeed, for any such MPE $((\mu^1, S^1), (\mu^2, S^2))$ and for each $x > x_{R^1}$, there exists an outcome-equivalent MPE $((\tilde \mu^1, \tilde S^1), (\tilde \mu^2, \tilde S^2))$ such that $\mathrm {supp} \, \tilde \mu^i = \mathrm {supp} \, \mu^i \cap (s^1 \vee s^2, \beta)$ for every player $i$ and $\tilde S^1 = (\alpha, s^1)$ and $\tilde S^2 = \emptyset$ (for an equilibrium of type 1), or $\tilde S^1 = \emptyset$ and $\tilde S^2 = (\alpha,s^2)$ (for an equilibrium of type 2).

By contrast, Corollary \ref{coromix} does not determine equilibrium outcomes of MPEs of types 1 and 2 for lower initial market conditions. First, as in Murto (2004), it may be possible to construct MPEs in which the stopping regions $S^1$ and $S^2$ exhibit gaps. Second, these gaps may themselves include randomization thresholds at which players exit the market with positive but finite intensity.

\subsection{The Characterization Result}

The next Theorem gives a necessary and sufficient condition for the existence  of an equilibrium of type 2. Similar results can be written for equilibria of types 1 and 3.

\begin{theorem} \label{CSmixed}
Let us consider two sequences of positive real numbers $ (q^{ 1}_i)_{1 \leq i \leq n}$ and  $ (q^{ 2}_i)_{1 \leq i \leq n-1}$ and $s^2\in \CI$  with  $q^{ 1}_1 = x_{R^1} > q^{ 2}_1 > q^{ 1}_2> ... > q^{ 2}_{n-1}>q^{ 1}_{n}> s^{ 2} $ and two sequences of positive real numbers
$ (a_i)_{1 \leq i \leq n}$ and $ (b_i)_{1 \leq i \leq n - 1}$.
 Then, the strategy profile $((\mu^1,S^1), (\mu^2, S^2)) \equiv ((\sum_{i=1}^{n} a_i \delta_{q^{ 1}_{i}}, \emptyset)$, $(\sum_{i=1}^{n-1} b_i \delta_{q^{ 2}_{i}}, (\alpha,s^{ 2}])))$ is a mixed-strategy MPE if and only if there exists a pair $(w^1,w^2)$ solution to the variational systems $\CV\CS^1$ and $\CV\CS^2$ below.\\

\noindent
Precisely:
\begin{itemize}
\item[(i)] $w^1 \in{\cal C}^0({\cal I}) \cap {\cal C}^2({\cal I} \setminus \{(q^{ 2}_i)_{1 \leq i \leq n-1 }, s^2 \})$ is solution to $\CV\CS^1$ if:
\begin{eqnarray}
\CL w^1 - rw^1 &=& 0   \;\text{ on } (s^{ 2},\beta)\setminus \{ (q^{ 2}_i)_{1 \leq i \leq n-1}\}
\label{var11n}\\
w^1(x) &=& G^1(x)  \;\text{ on } x \leq  s^{ 2},
 \label{var12n}\\
w^1(q^{ 1}_i) &=& R^1(q^{ 1}_i),   \; 1 \leq i \leq n \label{var13n}    \\
b_i [ G^1(q^{ 2}_i) - w^1(q^{ 2}_i)] + \frac{1}{2} \Delta (w^1)'(q^{ 2}_i)&=&0,  \; 1 \leq i \leq n-1  \label{var14n'}\\
\lim_{x\rightarrow \beta} w^1(x)&=& 0, \label{var15n} \\
w^1(x) &\geq& R^1(x) \;\text{ on } {\cal I}. \label{var16n}
\end{eqnarray}
\item[(ii)] $w^2\in {\cal C}^0({\cal I})  \cap {\cal C}^2({\cal I} \setminus (\{(q^{ 1}_i)_{1 \leq i \leq n}\} \cup  \{s^{ 2} \} ))$ is solution to $\CV\CS^2$ if:
\begin{eqnarray}
\CL w^2 - rw^2  &=& 0 \; \text{ on }  (s^{ 2},\beta)\setminus \{ (q^{ 1}_i)_{1 \leq i \leq n}\}
 \label{var11g}\\
w^2(x) &=& R^2(x) \;\text{ on } x \leq  s^{ 2}, \label{var11g'}\\
w^2(q^{ 2}_i) &=& R^2(q^{ 2}_i), \;    1 \leq i \leq n-1 \label{var13g}    \\
(w^2)'(s^{ 2}) &=& (R^2)'(s^{ 2}),    \label{var14g}\\
a_i [ G^2(q^{ 1}_i) - w^2(q^{ 1}_i)] + \frac{1}{2} \Delta (w^2)'(q^{ 1}_i)&=&0,  \; 1 \leq i \leq n   \label{var14g'}\\
\lim_{x \rightarrow \beta} w^2(x)&=& 0, \label{var15g} \\
w^2 &\geq& R^2 \;\text{ on } \CI.\label{var16g}
\end{eqnarray}
\end{itemize}
\end{theorem}

The proof is relegated to Appendix C and is based on the properties obtained from the proof of Theorem \ref{CNmixed} together with classical methods employed in optimal stopping and stopping games. In particular, the conditions \eqref{var14n'} and \eqref{var14g'} are obtained by applying the It\^{o}-Tanaka-Meyer formula.

\newpage

\renewcommand{\thesection}{Appendix A: Proofs for Section 3}

\renewcommand{\theequation}{A.\arabic{equation}}
\setcounter{equation}{0}

\section{}\label{appendix_on_section3}

We start with some definitions (Revuz and Yor (1999, Chapter I, \S4)). The process $X$ is defined over the canonical space $(\Omega,\mathcal F)$ of continuous trajectories, and $\mathbf P_\mu$ denotes the law of the process $X$ given an initial distribution $\mu \in \Delta(\mathcal I)$ over $\mathcal B( \mathcal I)$. We denote by $(\CF_t^0)_{t \geq 0}$ the natural filtration $(\sigma( X_s; s \leq t))_{t \geq 0}$ generated by $X$, and we let $\CF_\infty^0 \equiv \sigma(\bigcup_{t \geq 0} \mathcal F_t^0)$. For each $\mu\in \Delta(\mathcal I)$, we denote by $\CF_\infty^\mu$ the completion of $\CF _\infty^0$ with respect to $\mathbf P_\mu$, and, for each $t\geq 0$, we let $\CF_t^\mu$ be the augmentation of $\CF_t^0$ by the $\mathbf P_\mu$-null, $\CF_\infty^\mu$-measurable sets. The usual augmented filtration $(\CF_t)_{t \geq 0}$ is then defined by $\CF_t \equiv \bigcap_{\mu\in \Delta(\mathcal I)} \CF_t^\mu$ for all $t\geq 0$. Because the process $X$ is a Feller process in the sense of Revuz and Yor (1999, Chapter III, \S2, Definition 2.5) and a standard process in the sense of Blumenthal and Getoor (1968, Chapter I, Definition 9.2), the filtration $(\CF_t)_{t \geq 0}$ is actually right-continuous. As usual in this literature, we say that a property of the trajectories $\omega \in \Omega$ is satisfied \textit{almost surely} if it is satisfied $\mathbf P_\mu$-almost surely for all $\mu\in \CI$ or, equivalently, $\mathbf P_x$-almost surely for all $x\in \CI$.

\bigskip

\noindent \textbf{Proof of Lemma \ref{lem_sign}.} The proof proceeds along the same lines as that of Lemma 1 in D\'ecamps, Gensbittel, and Mariotti (2021). The result follows. \hfill $\blacksquare$

\bigskip

\noindent \textbf{Proof of Lemma \ref{lemma_properties_lambda}.} For each $\mu \in \Delta(\mathcal I)$, $\omega$ and $u^i$ are independent under $\mathbf P^i_\mu\equiv \mathbf P_\mu \otimes Leb $, and hence
\begin{align*}
\Gamma^i_t(\omega)= \mathbf P^i_\mu \hskip 0.3mm [\gamma^i \leq t \! \mid \! \CF](\omega)
\end{align*}
almost surely for all $t \geq 0$. We may assume that $\gamma(\cdot,u^i)\in \CT$ for all $u^i$, as we can  replace $\gamma^i$ by the constant stopping time $0$ for all $u^i$ in a Borel set of zero Lebesgue measure without  modifying the process $\Gamma^i$. Therefore, for all $u^i \in [0,1]$ and $t\geq 0$, we have $\{\omega \in \Omega : \gamma^i(\omega,u^i)\leq t\} \in \CF_t$ as $\gamma(\cdot,u^i)\in \CT$. Using Corollary 2 in Solan, Tsirelson, and Vieille (2012), this implies that $\Gamma^i_t$ is measurable with respect to the augmentation of $\CF_t$ by the $\mathbf P_\mu$-null, $\CF_\infty^\mu$-measurable sets, which coincides with $\CF_t^\mu$. As this is true for all $\mu \in \Delta(\mathcal I)$, we deduce that $\Gamma^i$ is adapted with respect to $\CF_t$. In particular, letting $\mu \equiv \delta_x$ yields
\begin{align*}
\Gamma^i_ t(\omega)=\mathbf P^i_x \hskip 0.3mm [\gamma^i \leq t \! \mid \! \CF_t](\omega)
\end{align*}
almost surely for all $t \geq 0$ by the law of iterated expectations. The result follows. \hfill $\blacksquare$

\bigskip

\noindent \textbf{Proof of Lemma \ref{ccdf}.} Assume that, for each $i=1,2$, $\gamma^i$ is a randomized stopping time with ccdf $\Gamma^i$. We have
\begin{align*}
\overline{\mathbf E}_x \!\left[1_{\{\gamma^i\leq \gamma^j\}}\, \mathrm e^{-r \gamma^i}R^i (X_{\gamma^i}) \right] \! &= \int_0^1 \!\int_0^1 \mathbf E _x \! \left[1_{\{\gamma^i(u^i)\leq \gamma^j(u^j)\}} \, \mathrm e^{-r \gamma^i(u^i)}R^i (X_{\gamma^i(u^i)})\right] \!\mathrm du^j \, \mathrm du^i
\\
&= \int_0^1 \mathbf E_x \! \left[\mathrm e^{-r \gamma^i(u^i)}R^i (X_{\gamma^i(u^i)}) \int_0^1 1_{\{ \gamma^i(u^i)\leq \gamma^j (u^j)\}}\, \mathrm du^j \right] \! \mathrm du^i \allowdisplaybreaks
\\
&= \int_0^1 \mathbf E_x \! \left[ \mathrm e^{-r \gamma^i(u^i)}R^i (X_{\gamma^i(u^i)}) \Lambda^j_{\gamma^i(u^i)-} \right] \! \mathrm du^i
\\
&=  \mathbf E_x \!\left[ \int_0^1 \mathrm e^{-r \gamma^i(u^i)}R^i (X_{\gamma^i(u^i)}) \Lambda^j_{\gamma^i(u^i)-} \, \mathrm du^i \right]
\\
&=  \mathbf E_x \! \left[ \int_{[0,\infty)} \mathrm e^{-r t}R^i (X_t) \Lambda^j_{t-} \, \mathrm d\Gamma_t^i \right] \hskip -1mm,
\end{align*}
where the second and fourth equalities follow from Fubini's theorem, and the third equality follows from the definition of $\Lambda^j$. The last equality follows from observing that, for $\mathbf P_x$-almost every $\omega \in \Omega$, $t \mapsto \Gamma ^i_t(\omega)$ is the cdf of the random variable $\gamma^i(\omega,\cdot)$ defined on the probability space $([0,1], \mathcal B([0,1]),Leb)$ and taking values in $[0,\infty]$, where $\Gamma^i_\infty (\omega) \equiv 1$ by convention; Fubini's theorem then implies that the random variable $u^i \mapsto \mathrm e^{-r \gamma^i(\omega,u^i)}R^i (X_{\gamma^i(\omega,u^i)}) \Lambda^j_{\gamma^i(\omega,u^i)-}$ is Lebesgue integrable over $[0,1]$ for $\mathbf P_x$-almost every $\omega \in \Omega$,\footnote{Recall that, by convention, this random variable is equal to 0 if $\gamma^i(\omega,u^i)=\infty$.} and we can thus apply the usual formula for the expectation. The proof for the second term appearing in (\ref{BDCOLL}) and (\ref{forbpr}) is similar and thus omitted.

Let us then verify that \eqref{inverse_stopping} defines a randomized stopping time in the sense of Definition \ref{def:randomized_stopping_time}. That $\hat\gamma^i(u^i)\in \CT$ for $Leb$-almost every $u^i \in [0,1]$ is standard (Jacod and Shiryaev (2003, Proposition I.1.28)). The random variable $(\omega,u^i)\mapsto \hat\gamma^i(u^i)(\omega)$ is $\CF_\infty\otimes \CB([0,1])$-measurable as it is nondecreasing and right-continuous with respect to $u^i$. That the ccdf associated with $\hat\gamma^i$ is $\Gamma^i$ is proven in De Angelis, Ferrari, and Moriarty (2018, Lemma 4.1), who use this representation as the definition of a randomized stopping time. The result follows. \hfill $\blacksquare$

\bigskip

\noindent
{\bf Proof of Theorem \ref{representation}}
If $\Gamma$ is a Markovian strategy, then $\Lambda=1-\Gamma$ satisfies almost surely for all $\tau \in \CT$, for all $s\geq 0$ 
\begin{equation}\label{markov2}
\Lambda_{\tau+s}=\Lambda_{\tau}(\Lambda_s \circ \theta_{\tau}), 
\end{equation}
on the event $\{\tau<\infty\}$, where $\theta_.$ denotes the shift operator on $\Omega$. 
In particular, this property applied at $\tau=s=0$ implies that $\Lambda_0=(\Lambda_0)^2$ and thus $\Lambda_0 \in \{0,1\}$ a.s.
In the terminology of Blumenthal and Getoor (1968), $\Lambda$ is a strong right-continuous multiplicative functional of $X$ adapted to $(\CF_t)$. 
The set $E=\{x | \mathbf{P}_x(\Lambda_0=1)=1\}$ is called the set of permanent points of $\Lambda$. Using Blumenthal's $0-1$ law and the fact that $\Lambda_0 \in \{0,1\}$, we have $\CI\setminus E=\{ x| \mathbf{P}_x(\Lambda_0=0)=1\}$. 
The stopping time $\tau=\inf \{ t > 0 | \Lambda_t=0\}\in \CT$ is called the lifetime of $\Lambda$. 
In order to apply the main result of Sharpe (1971), we need to prove that $\Lambda$ is an exact multiplicative functional, in the sense of Definition III.4.13  in Blumenthal and Getoor (1968). It is sufficient to prove that (see Proposition III.5.9 in Blumenthal and Getoor (1968)) 
\[ \forall x \notin E, \forall t>0, \; \lim_{u \downarrow 0}\mathbf E_x[ \Lambda_{t-u}\circ \theta_u]=0 .\]
Note that 
\[1_{\{t-u \geq \tau_{x}\circ \theta_u\}}\Lambda_{t-u}\circ \theta_{u} =0 \; \mathbf{P}_x-a.s.\]
so that for $u$ sufficiently small: 
\[ \mathbf E_x[ \Lambda_{t-u}\circ \theta_u] \leq \mathbf{P}_x( t-u < \tau_{x}\circ \theta_u) = \mathbf E_x[ \mathbf{P}_{X_u}(t-u < \tau_{x})]\leq \mathbf E_x[ \mathbf{P}_{X_u}(t/2 < \tau_{x})].\]
The map $y \rightarrow \mathbf{P}_y(\tau_x > t/2)$ is (universally) measurable and bounded, and has limit $0$ when $y\rightarrow x$ since $X$ is a regular diffusion. We conclude therefore by bounded convergence.
This implies that $E$ is open and thus that $\CI\setminus E$ is closed (see Blumenthal and Getoor (1968) p.126 last paragraph, using that the fine topology coincides with the usual topology in our case).

We can therefore apply Theorem 7.1 in Sharpe (1971). 
The stopping time appearing in the general expression (7.1) in Sharpe (1971) is called the exact regularization of $\tau$ and is equal in our particular case to the hitting time of a Borel set $B$ since the lifetime of $X$ is $+\infty$ and $X$ is continuous. Moreover, using that $X$ is a regular diffusion, this stopping time is a.s. equal to the hitting time of the closure of $B$, which we denote by $S$.
The product term in expression (7.1) of Sharpe (1971) is identically equal to $1$ and thus disappears (see Theorems 5.1 and 5.2 in Sharpe (1971), having in mind that all semipolar sets for $X$ are empty since $X$ is a regular diffusion, and that $X$ has continuous trajectories). Thanks to these remarks, there exists a  continuous additive functional $A$ with bounded $1$-potential (which means that $x \rightarrow \mathbf E_x[ \int_0^\infty e^{-t}dA_t]$ is bounded) and a positive Borel map $f$ such that
\[ \Lambda^i_t = 1_{t < \tau_S} e^{- \int_0^t f(X_s)dA_s} \; a.s. \]
Moreover, as stated in  Theorem 7.1 in Sharpe (1971), the integral $\int_0^t f(X_s)dA_s$ is $\mathbf{P}_x$ a.s. finite for all $t<\tau_S$ except maybe  for $x$ in a $\tau_S$-polar set. From the definition of a $\tau_S$ polar set (see the definition and Lemma 2.1 page 29 in Sharpe (1971)) we see that a $\tau_S$-polar set in our case must be a subset of $S$.

Using the classical representation result (see Borodin and Salminen (2002) chapter I.2 paragraph 23), there exists a positive Radon measure $\mu$ (locally finite) on $\CI$ such that $A_t=\int_{\CI} L^x_t d\mu(x)$ a.s.. Therefore, we have for $t < \tau$ 
\[ \tilde A_t=\int_0^t f(X_s)dA_s= \int_0^t \int_\CI f(X_s)dL^x_s d\mu(x)=\int_{\CI} L^x_t f(x) d\mu(x)   .\]
It follows that the measure $f(x)d\mu(x)$ on $\CI\setminus S$ is a Radon measure, i.e. a non-negative Borel measure which is locally finite. 
Indeed, assume by contradiction that it is not locally finite, then there exists $x \in \CI \setminus S$ such that for every $\varepsilon>0$ such that $(x-\varepsilon,x+\varepsilon) \subset \CI \setminus S$, we have  $\int_{(x-\varepsilon,x+\varepsilon)} f(x)d\mu(x)=\infty$. 
For all $t>0$, we have for all $\omega$ in a set of $\mathbf{P}_x$-probability $1$ that $L^x_t(\omega)>0$, and therefore since the local time of $X$ is a.s. jointly continuous, that $L^y_t(\omega)>0$ for $y\in [x-\varepsilon(\omega),x+\varepsilon(\omega)]$ for some $\varepsilon(\omega)>0$. This implies that if $0<t<\tau_S(\omega)$
\[   \tilde A_t(\omega)=\int_{\CI} L^x_t(\omega) f(x) d\mu(x) \geq \min_{y \in [x-\varepsilon(\omega),x+\varepsilon(\omega)]} L_t^y(\omega) \int_{(x-\varepsilon(\omega),x+\varepsilon(\omega))} f(x)d\mu(x)=+\infty .\]
Since $\mathbf{P}_x(\tau_S>0)=1$, this contradicts the statement of Theorem 7.1 in Sharpe (1971).

Reciprocally, if $S$ is a closed subset of $\CI$ and $\mu$ is a Radon measure on $\CI\setminus S$, then the process 
 \[ \Lambda_t = 1_{t < \tau_S} e^{- \int_\CI L^x_t d\mu(x) } ,\]
is well-defined and is a strong right-continuous multiplicative functional, and in particular satisfies \eqref{markov2}.
It can be checked directly that $\Gamma=1-\Lambda$ satisfies the assumptions of Lemma \ref{ccdf} and thus is the c.c.d.f. a randomized stopping time. \hfill $\blacksquare$

\bigskip

\noindent \textbf{Proof of Lemma \ref{usepure}.} We focus on player 1, the proof for player 2 being symmetrical. Observe from (\ref{forbpr}) that, for each $\tau^1 \in \mathcal T$, player 1's payoff from playing $\tau^1$ against $\Gamma^2$ is
\begin{align}
J^1(x,\tau^1,\Gamma^2) =  \mathbf E_x \! \left[ \mathrm e^{-r \tau^1}R^1(X_{\tau^1})\Lambda^2_{\tau^1-}+ \int_{[0,\tau^1)} \mathrm e^{-r t} G^1(X_t) \, \mathrm d\Gamma^2_t\right] \hskip -1mm. \label{forpbr1'}
\end{align}
Letting $\hat\gamma^1$ be the randomized stopping time associated to the ccdf $\Gamma^1$ by \eqref{inverse_stopping}, we have
\begin{align*}
J^1(x,\Gamma^1,\Gamma^2) &= \int_0^1 \mathbf E_x \! \left[ \mathrm e^{-r \hat \gamma^1(u^1)}R^1 (X_{\gamma^1(u^1)}) \Lambda^ 2_{\gamma^ 1(u^1)-}+ \int_{[0,\hat \gamma^1(u^1))} \mathrm e^{-r t} G^1(X_t) \, \mathrm d\Gamma^2_t\right] \!\mathrm d u^1 \allowdisplaybreaks
\\
&= \int_0^1 J^1(x,\hat \gamma^1(u^1),\Gamma^2)\, \mathrm du^1
\\
&\leq \sup_{u^1 \in [0, 1]} \,J^1(x,\hat\gamma^1(u^1) , \Gamma^2)
\\
&\leq \sup_{\tau^1 \in \CT} \,J^1(x, \tau^1, \Gamma^2).
\end{align*}
where the first equality follows along the same steps as in the proof of Lemma \ref{ccdf}, and the second equality follows from (\ref{forpbr1'}).
The result follows. \hfill $\blacksquare$

\bigskip
\noindent
{\bf Proof of Proposition \ref{geneprop}}
Before proving Proposition \ref{geneprop}, we need a few technical results related to continuous additive functional of diffusions. More precisely, we consider processes $A$ defined by $A_t=\int_{\CI}L^x_t d\mu(x)$ where $\mu$ is a Radon measure on $\CI\setminus S$ for some closed set $S$. \\

\noindent
If $(a,b) \subset \CI\setminus S$, the restriction of $\mu$ to $(a,b)$ is a Radon measure, and $A$ defines a continuous additive functional of the diffusion $X$ killed at $a$ and $b$, with state space $(a,b)$. \\

\noindent
If $\tau$ is the first exit time of $(a,b)$ with $[a,b] \subset \CI$, then
\begin{equation}\label{eq:formula_E_A}
\mathbf E_x[A_\tau]= \int_{\CI} \mathbf E_x[ L^y_{\tau} ] d\mu(y)=\int_{(a,b)} \mathbf E_x[ L^y_{\tau} ] d\mu(y)  =\int_{(a,b)} 2 (S'(y))^{-1} \Phi_{a,b}(x,y)d\mu(y), 
\end{equation}
where  
\[ \Phi_{a,b}(x,y)= \frac{(S(x\wedge y) - S(a))(S(b)-S(x\vee y))}{S(b)- S(a)}\]
denotes the Green function (potential density) of the diffusion $X$ killed at the boundaries $a$ and $b$. (see Borodin and Salminen (2002) page 20-21). It is easy to check that $\mathbf E_x[A_\tau]$ is finite if and only if for some $x\in(a,b)$:
\[ \int_a^{x} (S(y)-S(a))d\mu(y)< \infty \; \text{ and } \; \int_x^{b} (S(b)-S(y))d\mu(y)< \infty.\]
A more precise result can be shown (see Theorem 2.1 in Cetin (2018)). Precisely:
\begin{equation}\label{equiv_finite_A}
 \left\{ \begin{matrix} 
 A_{\tau_a}1_{\tau_a<\tau_b} =+\infty \, a.s. & \text{if} & \int_a^{x} (S(y)-S(a))d\mu(y)=+\infty \text{ for some $x\in(a,b)$} \\ A_{\tau_a}1_{\tau_a<\tau_b} <+\infty \, a.s. & \text{otherwise} &  \end{matrix} \right. 
\end{equation}
A symmetric result holds for $b$.   \\

\noindent
We will use the following lemma in the proof of Proposition \ref{lem:continuity}.
\noindent
\begin{lemma}\label{tech-lemma}

\noindent
Let $A_t=\int_{(a,b)} L^y_t d\mu(y)$ for some Radon measure $\mu$ on $(a,b)$.
Define the function $h$ by 
\[ \forall x\in(a,b), \; h(x)=\mathbf E_x[ C_a e^{-A_{\tau_a}} 1_{\{\tau_a< \tau_b \}} + C_b e^{- A_{\tau_b}} 1_{\{\tau_b< \tau_a \}} ], \] 
with $C_a,C_b \geq 0$ and $\tau_y=\inf\{ t \geq 0, X_t=y\}$. \\

\noindent

Then, $h$ is non-negative, $S$-convex and continuous on $(a,b)$, and the limits $h(a+)$ and $h(b-)$ exists and are given by 
\begin{equation}\label{eq_h(a+)}
h(a+)= \left\{ \begin{matrix} C_a & \text{if} & \int_a^{x} (S(y)-S(a))d\mu(y) < \infty \\ 0 & \text{otherwise} &  \end{matrix} \right. 
\end{equation}
\begin{equation}\label{eq_h(b-)}
 h(b-)= \left\{ \begin{matrix} C_b & \text{if} & \int_x^{b} (S(b)-S(y))d\mu(y) < \infty \\ 0 & \text{otherwise} &  \end{matrix} \right. .
\end{equation}
\end{lemma}

\bigskip
\noindent
{\bf Proof of Lemma \ref{tech-lemma}.}
$h$ is clearly non-negative. $S$-convexity follows directly from the Markov property applied for $h(\lambda x_1 +(1-\lambda)x_2)$ up to the stopping time $\tau_{x_1}\wedge \tau_{x_2}$. Precisely, we have
\[ h(\lambda x_1 +(1-\lambda)x_2)=\mathbf E_{\lambda x_1 +(1-\lambda)x_2}[ h(x_1) e^{-A_{\tau_{x_1}}} 1_{\{\tau_{x_1}< \tau_{x_2} \}} + h(x_2) e^{- A_{\tau_{x_2}}} 1_{\{\tau_{x_2}< \tau_{x_1} \}} ]. \]
Using that $e^{-A_t}\leq 1$, we obtain that $h$ is $S$-convex, i.e.
\[ h(\lambda x_1 +(1-\lambda)x_2) \leq  h(x_1) \frac{S(x_2)-S(\lambda x_1 +(1-\lambda)x_2)}{S(x_2)-S(x_1)}
+ h(x_2)\frac{S(\lambda x_1 +(1-\lambda)x_2)-S(x_1)}{S(x_2)-S(x_1)}  ]. \]
 
\noindent
$S$-convexity implies the continuity of $h$ on $(a,b)$. \\

\noindent
If $\int_a^{x_0} (S(y)-S(a))d\mu(y) < \infty$ for some $x_0 \in (a,b)$, then by property \eqref{equiv_finite_A}
$ h(x)=\mathbf E_x[C_b e^{- A_{\tau_b}} 1_{\{\tau_b< \tau_a \}} ]$ and thus
\[ 0 \leq h(x) \leq C_b \mathbf{P}_x(\tau_b < \tau_a) \underset{x \rightarrow a}{\longrightarrow} 0.\]
It follows that $h(a+)=0$.\\

\noindent
 If $\int_a^{x_0} (S(y)-S(a))d\mu(y) = \infty$ for some $x_0 \in (a,b)$, then by property \eqref{equiv_finite_A} $e^{-A_{\tau_a}}>0$ $\mathbf{P}_x$-a.s. 
If $a_n$ is a decreasing sequence with limit $a$ and $a_n<x$, applying the Markov property at $\tau_{a_n}$, we have
\[  h(x)=\mathbf E_x[ h(a_n) e^{-A_{\tau_{a_n}}} 1_{\{\tau_{a_n}< \tau_b \}} + C_b e^{- A_{\tau_b}} 1_{\{\tau_b< \tau_{a_n} \}} ]. \] 
Taking the limit along any subsequence such that $h(a_{n_k})$ converges to some $z$, we obtain
\[  h(x)=\mathbf E_x[ z e^{-A_{\tau_{a}}} 1_{\{\tau_{a}< \tau_b \}} + C_b e^{- A_{\tau_b}} 1_{\{\tau_b< \tau_{a} \}}], \]
and thus $z=C_a$ since $\mathbf E_x[e^{-A_{\tau_{a}}} 1_{\{\tau_{a}< \tau_b \}}]>0$. Therefore $h(a+)$ exists and equals $C_a$. \eqref{eq_h(b-)} follows by the same method. \hfill $\blacksquare$\\

\noindent
{\bf Strong Markov property:}
We  formalize in the next lemma a consequence of the strong Markov property that we will use several times throughout this Appendix.
\begin{lemma} \label{SMP}
If the  players use Markovian  randomized  stopping times with ccdf $\Gamma^1$ and $\Gamma^2$, then for every stopping time $T$ and every $x\in \CI$, their expected payoffs write as
\begin{equation}\label{MarkovpropertyJ}
J^i(x, \Gamma^1, \Gamma^2)=\mathbf E_x[ \int_{[0,T)}e^{-r s} R^i(X_s)\Lambda^j_{s-} d\Gamma^i_s+ \int_{[0,T)} e^{-r s}G^i(X_s)\Lambda^i_s d\Gamma^j_s + e^{-r T}\Lambda^j_{T-} \Lambda^i_{T -}  J^i(X_{T}, \Gamma^1, \Gamma^2)]
\end{equation}
\end{lemma}

\bigskip
\noindent
{\bf Proof of Lemma \ref{SMP}.}
It follows from Lemma \ref{ccdf} that
\begin{align*}
 J^i(x,\Gamma^1, \Gamma^2)& =\mathbf E_x \Big[ \int_{[0,T)}e^{-r s} R^i(X_s)\Lambda^j_{s-} d\Gamma^i_s+ \int_{[0,T)} e^{-r s}G^i(X_s)\Lambda^i_s d\Gamma^j_s  \\
& \qquad + e^{-r T} R^i(X_{T})\Lambda^j_{T-}  (\Gamma^i_{T}-\Gamma^i_{T-})+  e^{-r {T}}G^i(X_{T})\Lambda^i_{T} (\Gamma^j_{T}-\Gamma^j_{T-}) \\
& \qquad + \int_{(T,\infty)}e^{-r s} R^i(X_s)\Lambda^j_{s-} d\Gamma^i_s+ \int_{(T,\infty)} e^{-r s}G^i(X_s)\Lambda^i_s d\Gamma^j_s \Big] 
\end{align*}
Recall that the only jump of $\Lambda^i$ occurs at $\tau_{S^i}$ at which time the process jumps to zero and remains equal to zero. On the one hand, we have
\begin{align*}
e^{-r T} &R^i(X_{T})\Lambda^j_{T-}  (\Gamma^i_{T}-\Gamma^i_{T-})+  e^{-r {T}}G^i(X_{T})\Lambda^i_{T} (\Gamma^j_{T}-\Gamma^j_{T-}) \\
&=e^{-r T} R^i(X_{T})\Lambda^j_{T-}\Lambda^i_{T-}1_{\tau_{S^j} \geq T = \tau_{S^i}}+ e^{-r {T}}G^i(X_{T})\Lambda^i_{T-}\Lambda^j_{T-}1_{\tau_{S^i} > T = \tau_{S^j}} \\
&=e^{-r T}  J^i(X_{T})\Lambda^j_{T-}\Lambda^i_{T-}1_{\tau_{S^j} \geq T = \tau_{S^i}}+ e^{-r {T}} J^i(X_{T},\Gamma^1, \Gamma^2)\Lambda^i_{T-}\Lambda^j_{T-}1_{\tau_{S^i} > T = \tau_{S^j}}.
\end{align*}
On the other hand, we have
\begin{align*}
\int_{(T,\infty)}&e^{-r s} R^i(X_s)\Lambda^j_{s-} d\Gamma^i_s+ \int_{(T,\infty)} e^{-r s}G^i(X_s)\Lambda^i_s d\Gamma^j_s \\
&=e^{-r T}1_{T < \tau_{S^i} \wedge \tau_{S^j}}\left( \int_{(0,\infty)}e^{-r s} R^i(X_{T+s})\Lambda^j_{(T +s)-} d\Gamma^i_{T +s}+ \int_{(0,\infty)} e^{-r s}G^i(X_{T+s})\Lambda^i_{T+s} d\Gamma^j_{T +s} \right) \\
&=e^{-r T}1_{T < \tau_{S^i} \wedge \tau_{S^j}}\Lambda^j_{T} \Lambda^i_{T} \Big( \int_{(0,\infty)}e^{-r s} R^i(X_{s}\circ \theta_{T})(\Lambda^j_{s-} \circ \theta_{T}) d(\Gamma^i_{s}\circ \theta_T) \\
& \qquad\qquad\qquad + \int_{(0,\infty)} e^{-r s}G^i(X_{s}\circ \theta_{T})(\Lambda^i_{s}\circ \theta_{T}) d(\Gamma^2_{s}\circ \theta_{T}) \Big) \\
&=e^{-r T}1_{T < \tau_{S^i} \wedge \tau_{S^j}}\Lambda^j_{T-} \Lambda^i_{T-} \Big( \int_{[0,\infty)}e^{-r s} R^i(X_{s}\circ \theta_{T})(\Lambda^j_{s-} \circ \theta_{T}) d(\Gamma^i_{s}\circ \theta_T) \\
& \qquad\qquad\qquad + \int_{[0,\infty)} e^{-r s}G^i(X_{s}\circ \theta_{T})(\Lambda^i_{s}\circ \theta_{T}) d(\Gamma^j_{s}\circ \theta_{T}) \Big),
\end{align*}
where we used \eqref{mark} in the second equality.
Taking expectation and applying the Markov property at $T$, it leads to
\begin{align*}
\mathbf E_x \Big[ \int_{(T,\infty)}e^{-r s} R^i(X_s)\Lambda^j_{s-} d\Gamma^i_s+ \int_{(T,\infty)} e^{-r s}G^i(X_s)\Lambda^i_s d\Gamma^j_s \Big] = \mathbf E_x  \Big[ e^{-r T}1_{T < \tau_{S^i} \wedge \tau_{S^j}}\Lambda^j_{T-} \Lambda^i_{T-} J^i(X_T,\Gamma^1, \Gamma^2) \Big]
\end{align*}
The conclusion follows by putting the pieces together and noticing that $e^{-r T}\Lambda^j_{T-} \Lambda^i_{T -}  J^i(X_{T},\Gamma^1, \Gamma^2)$ vanishes on the event $\{ T> \tau_{S^i} \wedge \tau_{S^j}\}$. \hfill $\blacksquare$ \\

\noindent
{\bf Proof of Proposition \ref{geneprop}.}
Let us fix  Markov strategies $(\mu^i,S^i)$ for $i=1,2$. Let $A^i_t=\int_{\CI}L_{t}^a d\mu^i(a)$ denote the continuous additive functional associated to $\mu^i$, $\Lambda^i_t=e^{-A^i_t}1_{t<\tau_{S^i}}$ and $\Gamma^i_t=1-\Lambda^i_t$.\\

\noindent
The payoff of player $1$ for a pure stopping time $\tau^1$ is given by 
\begin{align*}
 J^1 &(x,\tau^1,\Gamma^2) = \\
 &\mathbf E_x[ e^{-r \tau^1}\Lambda^2_{\tau^1-}R^1(X_{\tau^1})+ \int_{[0,\tau^1 \wedge \tau_{S^2})} e^{-r s}G^1(X_s)d\Gamma^2_s +e^{-r \tau_{S^2}}\Lambda^2_{\tau_{S^2}-}G^1(X_{\tau_{S^2}})1_{\tau_{S^2}<\tau^1}].
\end{align*}
Note that $\bar  J^1$ is equal to $R^1$ on $S^1$ and to $G^1$ on $S^2\setminus S^1$.\\

\noindent
Using assumption A6, we first prove that $V_{R^1}(x) \leq \bar J^1(x) \leq G^1(x)$. 
These inequalities imply $\lim_{x \rightarrow \alpha+} \frac{\bar J^1(x)}{\phi(x)}=0$ and $\lim_{x \rightarrow \beta-} \frac{\bar J^1(x)}{\psi(x)}=0$.\\

\noindent 
For the first inequality, let $\tau^1_*=\tau_{(\alpha,x_{R^1}]}$, and $\hat\gamma^2(u)$ be defined by \eqref{inverse_stopping}. Using Lemma \ref{ccdf} and that $G^1\geq V_{R^1}$ we obtain
\begin{align*}
\bar J^1(x) \geq  J^1 &(x,\tau_*^1, \Gamma^2)= \int_0^1  J^1 (x,\tau_*^1, \hat\gamma^2(u))du  \geq \\
 & \int_0^1 \mathbf E_x[ e^{-r \tau_*^1}R^1(X_{\tau_*^1})1_{\tau_*^1\leq \hat\gamma^2(u)} + e^{-r \hat\gamma^2(u)}V_{R^1}(X_{\hat \gamma^2(u)})1_{\tau_*^1 > \hat\gamma^2(u)}] du.
\end{align*}
For all $u$, we have on $\{\tau_*^1 > \hat\gamma^2(u)\}$,  
\begin{equation}\label{eq-condition-VR1}
e^{-r \hat\gamma^2(u)}V_{R^1}(X_{\hat \gamma^2(u)})= \mathbf E_x[ e^{-r\tau_*^1}R^1(X_{\tau_*^1})|\CF_{\hat \gamma^2(u)}] \; \mathbf{P}_x-a.s.
\end{equation} 
Therefore, using the tower property of conditional expectation
\begin{align*}
\mathbf E_x[ e^{-r \tau_*^1}R^1(X_{\tau_*^1})1_{\tau_*^1\leq \hat\gamma^2(u)} + e^{-r \hat\gamma^2(u)}V_{R^1}(X_{\hat \gamma^2(u)})1_{\tau_*^1 > \hat\gamma^2(u)}]
=\mathbf E_x[ e^{-r \tau_*^1}R^1(X_{\tau_*^1})]=V_{R^1}(x),
\end{align*}
and we conclude that 
\[ \bar J^1(x) \geq  \int_0^1 V_{R^1}(x) du = V_{R^1}(x).\]
For the second inequality, we have by assumption  A6 $R^1\leq V_{R^1} \leq G^1$. 
Therefore for every $\tau^1\in \CT$
\begin{align*}
 J^1 (x,\tau^1, \Gamma^2) &= \int_0^1  J^1 (x,\tau^1, \hat\gamma^2(u))du  \\
 & \leq \int_0^1 \mathbf E_x[ e^{-r \tau^1}G^1(X_{\tau^1})1_{\tau^1\leq \hat\gamma^2(u)} + e^{-r \hat\gamma^2(u)}G^1(X_{\hat \gamma^2(u)})1_{\tau^1 > \hat\gamma^2(u)}] du \\
&=\int_0^1 \mathbf E_x[ e^{-r (\tau^1 \wedge \hat\gamma^2(u))}G^1(X_{\tau^1 \wedge \hat\gamma^2(u)})] du  \\
& \leq \int_0^1 G^1(x) du = G^1(x),
\end{align*}
where we used that $e^{-rt}G^1(X_t)$ is a supermartingale by assumption A7. \\

\noindent
Let us prove $(i)$: It is not optimal for player $1$ to stop immediately in $S^2$ if $R^1<G^1$. Therefore, for any PBR, $S^1\cap S^2 \cap(\alpha^1,\beta)=\emptyset$ which proves $(i)$. \\

\noindent
Let us prove $(ii)$: The inclusion $S^1 \subset C^1=\{ x | \bar J^1(x)=R^1(x)\}$ is obvious from the definition of a PBR.\\

\noindent
We now prove $(iii)$: $supp(\mu^1) \subset C^1 \cup S^2$.\\

\noindent
Assume that $x \in supp(\mu^1)\setminus S^2$. Recall that by Lemma \ref{ccdf} 
\[ \bar  J^1 (x) = \int_0^1 J^1(x,\hat\gamma^1(u), \Gamma^2)du,\]
where $\hat\gamma^1(u)=\inf \{ t \geq 0 \,|\, \Gamma^1_t>u\}$. Moreover, since for all $u\in [0,1]$, $J^1(x,\hat\gamma^1(u), \Gamma^2)\leq \bar  J^1 (x)$, the latter inequality is an equality for all $u\in M$ where $M$ has Lebesgue measure $1$.

From the definition of $\Gamma^1$, we have
  $\hat\gamma^1(u)= \inf\{ t \,|\, 1-e^{-A^1_t} > u\} \wedge \tau_{S^1}$ for $u\in [0,1)$. 
We have $\hat\gamma^1(u)>0$ $\mathbf{P}_x$-a.s. for all $u>0$ as $A^1$ is continuous. Using that $x \in supp(\mu^1)$, we also have $\lim_{u \rightarrow 0}\hat\gamma^1(u) =0$ $\mathbf{P}_x$-a.s.. Indeed, $\hat\gamma^1(u)$ is non-decreasing with respect to $u$ for all $\omega$ and converges to $\hat\gamma^1(0)=\inf\{ t \geq 0 | A^1_t>0\} \wedge \tau_{S^1}$. Then, arguing as in Lemma 2.16 p. 416 in Revuz and Yor (1999), we obtain that $\hat\gamma^1(0)=\tau_{supp(\mu^1)} \wedge \tau_{S^1}$ a.s., and therefore that $\hat\gamma^1(0)=0$ $\mathbf{P}_x$-a.s..\\

\noindent
Precisely, the local times of $X$ have a continuous version $(t,y)\rightarrow L^y_t$, and we know that for all $t>0$, $L^x_t>0$ $\mathbf{P}_x$-a.s.. We deduce that for every $\omega$ in a set of $\mathbf{P}_x$-probability one, $L^y_t>0$ for $y$ in a neighborhood of $x$, which implies 
\[ A^1_t = \int_{\CI} L^y_t d\mu^1(y) >0 \; \mathbf{P}_x-a.s.,\]
and thus $\hat\gamma^1(0)=0$ $\mathbf{P}_x$-a.s.. 

For all $u \in M$:
\[ \bar  J^1 (x)=J^1(x,\hat\gamma^1(u), \Gamma^2)=\mathbf E_x[ \int_{[0,\hat\gamma^1(u))} e^{-r s}G^1(X_s)d\Gamma^2_s + e^{-r \hat\gamma^1(u)}\Lambda^2_{\hat\gamma^1(u)-}R^1(X_{\hat\gamma^1(u)})]\]
Taking the limit as $u \in M$ goes to zero, we deduce that $\bar  J^1 (x)=R^1(x)$ by bounded convergence.\\

\noindent
Proof of $(iv)$. We  have 
\begin{equation}\label{point(iv)} 
\bar  J^1(x) \geq J^1(x,(0,(\alpha,x_{R^1}]),(\mu^2,S^2)) \geq J^1(x,(0,(\alpha,x_{R^1}]),(0,\emptyset)).
\end{equation}
The first inequality in \eqref{point(iv)} follows from the fact that $(\mu^1,S^1)$ is a PBR to $(\mu^2,S^2)$. For the second one, recall that by assumption A6 
\[ G^1(x) \geq V_{R^1}(x)=\sup_{\tau}\mathbf E_x[ e^{- r \tau}R^1(X_\tau)]=\mathbf E_x[ e^{- r \tau_*^1}R^1(X_{\tau_*^1})] ,\]
where $\tau_*^1=\tau_{(\alpha,x_{R^1}]}$. We have 
\begin{align*}
 J^1 (x,\tau_*^1, \Gamma^2) &= \int_0^1 J^1 (x,\tau_*^1, \hat\gamma^2(u))du \\
 &= \int_0^1 \mathbf E_x[ e^{-r \tau_*^1} R^1(X_{\tau_*^1})1_{\tau_*^1\leq \hat\gamma^2(u)} + e^{-r \hat\gamma^2(u)}G^1(X_{\hat\gamma^2(u)})1_{\tau_*^1 > \hat\gamma^2(u)}] du .
\end{align*}
On $\{\tau_*^1>\hat\gamma^2(u)\}$, using assumption A6, we have $\mathbf{P}_x$-a.s. 
\begin{equation}\label{eq-condition-f_1}
e^{-r \hat\gamma^2(u)}G^1(X_{\hat\gamma^2(u)}) \geq e^{-r \hat\gamma^2(u)}V_{R^1}(X_{\hat\gamma^2(u)})= \mathbf E_x[ e^{-r\tau_*^1}R^1(X_{\tau_*^1})|\CF_{\hat\gamma^2(u)}].
\end{equation} 
Therefore, using the tower property of conditional expectation
\begin{align*}
J^1 (x,\tau_*^1, \hat\gamma^2(u)) \geq \mathbf E_x[ e^{-r\tau_*^1}R^1(X_{\tau_*^1})],
\end{align*}
and the second inequality of \eqref{point(iv)} follows by integrating with respect to $u$.\\

\noindent
The conclusion follows by noticing  that $J^1(x,(0,(\alpha,x_{R^1}]),(0,\emptyset))>R^1(x)$ for $x>x_{R^1}$ and applying $(ii)$ and $(iii)$.\\

\noindent
Proof of $(v)$. We have from Lemma \ref{ccdf} 
\[ \bar J^1(x)= \int_0^1 J^1(x,\hat\gamma^1(u),\Gamma^2) du ,\]
and for all $u\in [0,1]$, $J^1(x,\hat\gamma^1(u),\Gamma^2)\leq \bar J^1(x)$, so that the latter inequality is an equality for all $u\in M$, where $M$ has Lebesgue measure $1$. 

We will use that $\hat\gamma^1(u) \rightarrow \tau_{S^1}$ as $u\rightarrow 1$. This property follows from the definition of $\Gamma^1$ since we have $\tau_{S^1}=\inf\{t \geq 0 \,|\, \Gamma^1=1\}$, and $\Gamma^1_t<1$ on $t< \tau_{S^1}$. Moreover $u \rightarrow \hat\gamma^1(u)$ is non-decreasing. 
We have for $u\in M$
\[ \bar J^1(x)= J^1(x,\hat\gamma^1(u),\Gamma^2)= \mathbf E_x[ \int_{[0,\hat\gamma^1(u))} e^{-r s}G^1(X_s)d\Gamma^2_s + e^{-r \hat\gamma^1(u)}\Lambda^2_{\hat\gamma^1(u)-}R^1(X_{\hat\gamma^1(u)})].\]
Taking the limit as $u\rightarrow 1$ with $u\in M$,we obtain by bounded convergence
\[ \bar J^1(x)= \mathbf E_x[ \int_{[0,\tau_{S^1})} e^{-r s}G^1(X_s)d\Gamma^2_s + e^{-r \tau_{S^1}}\Lambda^2_{\tau_{S^1}-}R^1(X_{\tau_{S^1}})]=J^1(x,\tau_{S^1},\Gamma^2),\]
which concludes the proof of the first assertion.\\

\noindent
For the second assertion, let $\tilde \Gamma^1$ be associated with $(\tilde \mu^1,S^1)$ and $\tilde \gamma^1(u) =\inf \{ t \geq 0 \,|\, \tilde \Gamma^1_t > u\}$. By assumption, we have that for all $u$
\begin{equation}\label{eq_J1R1} 
\bar J^1(X_{\tilde \gamma^1(u)})=R^1(X_{\tilde \gamma^1(u)}).
\end{equation}
Therefore
\[ J^1(x,\tilde \Gamma^1,\Gamma^2)=\int_0^1 \mathbf E_x[\int_{[0,\tilde\gamma^1(u))} e^{-r s}G^1(X_s)d\Gamma^2_s + e^{-r \tilde\gamma^1(u)}\Lambda^2_{\tilde\gamma^1(u)-}R^1(X_{\tilde\gamma^1(u)})] du .\]
Applying the Markov property at $\tilde\gamma^1(u)$ and using that $\bar  J^1(.) = J^1(.,(0,S^1),\Gamma^2)$, we have
\begin{align*}
\bar  J^1(x) &= \mathbf E_x[\int_{[0,\tilde\gamma^1(u))} e^{-r s}G^1(X_s)d\Gamma^2_s + e^{-r \tilde\gamma^1(u)}\Lambda^2_{\tilde\gamma^1(u)-}\bar  J^1(X_{\tilde\gamma^1(u)})]\\
&= \mathbf E_x[\int_{[0,\tilde\gamma^1(u))} e^{-r s}G^1(X_s)d\Gamma^2_s + e^{-r \tilde\gamma^1(u)}\Lambda^2_{\tilde\gamma^1(u)-} R^1(X_{\tilde\gamma^1(u)})],
\end{align*}
where the second equality follows from \eqref{eq_J1R1}.
We conclude that
\[ J^1(x,\tilde \Gamma^1,\Gamma^2)=\int_0^1 \bar  J^1(x) du= \bar  J^1(x),\] 
which ends the proof.  \hfill $\blacksquare$

\bigskip

\noindent
{\bf Proof of Proposition \ref{lem:continuity}}

We first prove the following Lemma on the continuity of brvf.
\begin{lemma} \label{iii}
If $(\mu^1,S^1)$ is a pbr to $(\mu^2,S^2)$ with associated brvf $\bar J^2$,  then $\bar J^1$ is continuous on any interval $[a,b]$ such that $(a,b) \subset \CI\setminus (S^1\cup S^2)$.
\end{lemma}

\noindent
{\bf Proof of Lemma \ref{iii}}
Given $x \notin S^1 \cup S^2$, define $T_\varepsilon= \tau_{x+\varepsilon}\wedge \tau_{x-\eta} $ such that $\varepsilon,\eta>0$ and $[x-\eta,x+\varepsilon] \subset (S^1)^c \cap (S^2)^c$. Lemma \ref{SMP}  applied at $T_\varepsilon$ implies:
\begin{align*}
 \bar J^1(x)& = \mathbf E_x[ \int_{[0,T_\varepsilon)}e^{-r s} R^1(X_s)\Lambda^2_{s-} d\Gamma^1_s+ \int_{[0,T_\varepsilon)} e^{-r s}G^1(X_s)\Lambda^1_s d\Gamma^2_s + e^{-r T_\varepsilon}\Lambda^2_{T_\varepsilon -} \Lambda^1_{T_{\varepsilon} -} \bar{J}^1(X_{T_\varepsilon})] \\
& =\mathbf E_x[ \int_{[0,T_\varepsilon)}e^{-r s} R^1(X_s)\Lambda^2_s d\Gamma^1_s+ \int_{[0,T_\varepsilon)} e^{-r s}G^1(X_s)\Lambda^1_s d\Gamma^2_s + e^{-r T_\varepsilon}\Lambda^2_{T_\varepsilon} \Lambda^1_{T_{\varepsilon}} \bar J^1(X_{T_\varepsilon})].
\end{align*}
The last equality follows from the fact that on $\{t \leq T_\varepsilon\}$, we have $\Lambda^i_{t-}=\Lambda^i_t$.

$T_\varepsilon$ goes to zero in probability as $\varepsilon$ goes to zero ($\eta$ being fixed) (even a.s. since it is decreasing with $\varepsilon$) and $\mathbf{P}_x(T_\varepsilon= \tau_{x+\varepsilon})$ goes to $1$. Therefore bounded convergence implies that $\bar  J^1(x)=\bar  J^1(x+)$. Similarly, we can prove $\bar  J^1(x)=\bar  J^1(x-)$. \\

\noindent
Now, let us consider an interval $(a,b) \subset (S^1)^c \cap (S^2)^c$.  
$\bar  J^1$ is continuous on $(a,b)$ thanks to the preceding argument. Le us assume that $a>\alpha$ and belongs to $S^1$ so that $\bar  J^1(a)=R^1(a)$, we will prove that $\bar  J^1$ is right-continuous at $a$.\\

\noindent
Using Proposition \ref{geneprop}-$(v)$, $\bar J^1(x)=J^1((0,S^1),(\mu^2,S^2))$.  Applying Lemma \ref{SMP} with $\tau=\tau_a \wedge \tau_b$, we have for $x\in(a,b)$:
\begin{equation}\label{eq_J13}
\bar  J^1(x) =\mathbf E_x[  \int_{[0,\tau)} e^{-r s}G^1(X_s) d\Gamma^2_s + e^{-r\tau}e^{- A^2_\tau} \bar  J^1(X_{\tau})] .
\end{equation}
Since $a\notin S^2$, the measure $\mu^2$ is locally finite at $a$. Therefore,
\begin{equation}\label{eq_upG1} 
0\leq \mathbf E_x[\int_{[0,\tau)} e^{-r s}G^1(X_s)d\Gamma^2_s ] \leq  C \mathbf E_x[1-e^{-A^2_\tau}],
\end{equation}
where $C$ is an upper bound on $G^1$ on $[a,b]$. Applying Lemma \ref{tech-lemma} (with $C_a=C_b=1$ and $\mu=\mu^2$), $\mathbf E_x[1-e^{-A^2_\tau}]$ goes to $0$ as $x$ converges to $a$.\\

\noindent
Lemma \ref{tech-lemma} also implies that $\mathbf E_x[e^{-r\tau}e^{- A^2_\tau}\bar  J^1(X_{\tau})]$ goes to $\bar  J^1(a)=R^1(a)$ as $x$ converges to $a$ (by taking $\mu=\mu^2+r Leb$), and thus $\bar J^1$ is right-continuous at $a$.\\

\noindent
Let us now consider the case $a \in S^2$ so that $\bar  J^1(a)=G^1(a)$, we will prove that $\bar  J^1$ is right-continuous at $a$. Let $0<\varepsilon < b-a$. As before with $\tau=\tau_a \wedge \tau_{a+\varepsilon}$, we have for $x\in(a,a+\varepsilon)$:
\begin{equation}\label{eq_J14}
\bar  J^1(x) =\mathbf E_x[  \int_{[0,\tau)} e^{-r s}G^1(X_s) d\Gamma^2_s + e^{-r\tau}e^{- A^2_\tau} \bar  J^1(X_{\tau})] .
\end{equation}
If $\int_{a}^x (S(y)-S(a))d\mu^2(y)<\infty$, the proof is completely similar as for the case $a\in S^1$. Let us assume that $\int_{a}^x (S(y)-S(a))d\mu^2(y)=\infty$.
Lemma \ref{tech-lemma} implies that $\mathbf E_x[e^{-r\tau}e^{- A^2_\tau}\bar  J^1(X_{\tau})]$ goes to zero  as $x$ converges to $a$ (by taking $\mu=\mu^2+r Leb$). On the other hand
\begin{align*}
\mathbf E_x[  \int_{[0,\tau)} e^{-r s}G^1(X_s) d\Gamma^2_s ] \geq (\min_{[a,a+\varepsilon]} G^1) \mathbf E_x [e^{-r \tau} (1-e^{- A^2_\tau})] = (\min_{[a,a+\varepsilon]} G^1) \mathbf E_x [e^{-r \tau}  - e^{-r \tau}e^{- A^2_\tau})].
\end{align*}
By Lemma \ref{tech-lemma}, the last expectation goes to $1$ as $x$ converges to $a$.
We deduce that $\liminf_{x \rightarrow a+} \bar J^1(x) \geq (\min_{[a,a+\varepsilon]} G^1)$ and the conclusion follows by sending $\varepsilon$ to zero.
The cases  $b\in S^1, b\in S^2$ can be proven in a similar way. Lemma \ref{iii} is proven \hfill $\blacksquare$.\\

We now prove a second lemma related to the pbr on the intervals $(\alpha,\alpha^i]$.
\begin{lemma}\label{lem:pbralphai}
The following properties hold: 
\begin{itemize}
\item If $(\mu^i,S^i)$ is a pbr to $(\mu^j,S^j)$, then $(\alpha,\alpha^i]\subset S^1\cup S^2$. 
\item If $((\mu^1,S^1),(\mu^2,S^2))$ is a MPE, then either $[\alpha^1 \wedge \alpha^2, \alpha^1\vee \alpha^2] \subset S^1$ or $[\alpha^1 \wedge \alpha^2, \alpha^1\vee \alpha^2] \subset S^2$.  
\end{itemize}
\end{lemma}

\noindent
{\bf Proof of Lemma \ref{lem:pbralphai}}
We prove the first point with $i=1$. Note that $\bar J^1=R^1=V_{R^1}=G^1$ on $(\alpha,\alpha^1]$. 
Assume by contradiction that $x \in (\alpha,\alpha^1)$ and $x\notin S^1\cup S^2$. Let $(a,b) \subset \CI\setminus (S^1\cup S^2)$ with $b<\alpha^1$ which contains  $x$. Since $(0,S^1)$ is a PBR to $(\mu^2,S^2)$, and using Lemma \ref{SMP} at $\tau$ the exit time of $(a,b)$ we have
\begin{align*} 
\bar J^1(x)&=\mathbf E_x[\int_{[0,\tau)} e^{-r s}G^1(X_s)d\Gamma^2_s + e^{-r \tau}\Lambda^2_{\tau-}\bar  J^1(X_{\tau})] \\
&= \mathbf E_x[\int_{[0,\tau)} e^{-r s}R^1(X_s)d\Gamma^2_s + e^{-r \tau}\Lambda^2_{\tau-}R^1(X_{\tau})] \\
&= \int_0^1\mathbf E_x[e^{-r \hat \gamma^2(u,.)}R^1(X_{\hat \gamma^2(u,.)})1_{\hat \gamma^2(u,.)<\tau} + e^{-r \tau}1_{\tau \leq \hat \gamma^2(u,.)}R^1(X_{\tau})] du \\
&= \int_0^1 \mathbf E_x[e^{-r (\hat \gamma^2(u,.) \wedge \tau)}R^1(X_{\hat \gamma^2(u,.)\wedge \tau})] du \\
&< R^1(x),
\end{align*}
where the third equality is obtained as in Lemma \ref{ccdf}, and the last inequality follows from assumption A3 together with the fact that for for all $u>0$, $\tau \wedge \hat \gamma^2(u,.) >0$ $\mathbf{P}_x$-almost surely since $\Gamma^2$ is continuous on $[0,\tau_{S^2})$ and $\tau_{S^2}>0$ $\mathbf{P}_x$-almost surely. The above inequality contradicts $\bar J^1 \geq R^1$ and this concludes the proof of the first point.\\

\noindent
Let us prove the second point. Assume without loss of generality that $\alpha^2 < \alpha^1$. By Proposition \ref{geneprop}, we have $S^1\cap S^2 \cap (\alpha^2,\alpha^1]=\emptyset$. However, the preceding point implies that $(\alpha^2,\alpha^1] \subset S^1\cup S^2$. Using a connectedness argument, it must be that either $S^1\cap(\alpha^2,\alpha^1]$ or $S^2\cap(\alpha^2,\alpha^1]$ is empty. This concludes the proof of the second point.  \hfill $\blacksquare$.\\

\noindent
Let us now prove the Proposition. Assume without loss of generality that $\alpha^2\leq \alpha^1$.\\

\noindent
For $i=1,2$, the payoff at equilibrium $\bar {J}^i$ is  continuous on $(\alpha,\alpha^1]$ by the preceding points since it is equal to $R^i$ or $G^i$ on this interval. Note that $\bar J^1$ is right-continuous at $\alpha^1$ since $R^1 \leq \bar J^1\leq G^1$, and that the same is true for $\bar J^2$ if $\alpha^2=\alpha^1$. \\

By Lemma \ref{iii},  $\bar {J}^i$ is continuous on any interval $[a,b]$ such that $(a,b) \subset \CI\setminus (S^1\cup S^2)$. Moreover, $\bar {J}^i$ is also continuous on $S^1$ and on $S^2$.  Therefore, if  $\bar J^1$ is not right-continuous at $x$, it must be that $x\geq \alpha^1$, that $x\in S^1\cup S^2$ and that for every $\varepsilon>0$, $[x,x+\varepsilon)$ intersects both $S^1\cup S^2$ and $\CI \setminus (S^1 \cup S^2)$. 
Let us  consider the case that $x\in S^1$  and  $x> \alpha^1$ and assume by contradiction that $\bar {J}^i$ is not right-continuous at $x$. Since $(\alpha^1,\beta)\cap S^1\cap S^2=\emptyset$, we may choose $\varepsilon$ sufficiently small so that $[x,x+\varepsilon) \cap S^2 =\emptyset$. 
If $(a,b)$ is a connected component of the open set $[x,x+\varepsilon) \setminus S^1$, so that $a,b\in S^1$, then it must be that $\mu^2((a,b))>0$. By contradiction, if $\mu^2((a,b))=0$, then for $y \in(a,b)$: 
\[ \bar J^1(y)={J}^1(y,(0,S^1),(\mu^2,S^2))=\mathbf E_y[e^{-r\tau_{S^1}}R^1(X_{\tau_{S^1}})]< R^1(y),\]
where we used that $b\leq x_{R^1}$ and assumption A3. This contradicts the fact that $\bar J^1 \geq R^1$.  
Therefore, there exists some $y\in (a,b)$ such that $\bar J^2(y)=R^2(y)$.
 As this is true for every connected component, it must be that there exists a decreasing sequence $y_n$ with limit $x$ such that $\bar J^2(y_n)=R^2(y_n)$ and a sequence of connected components $(a_n,b_n)$ such that $y_n\in (a_n,b_n)$ whose length goes to zero.
By Proposition \ref{geneprop}, we have $\bar J^2(a_n)= G^2(a_n)$ and $\bar J^2(b_n)= G^2(b_n)$. Recall also that $G^2(x) - R^2(x) >0$. On the interval $(a_n,b_n)$, since $(0,S^2)$ is a best reply to $(\mu^1,S^1)$, we have by Lemma \ref{SMP} 
\begin{align*}
 \bar J^2 (y_n) &=\mathbf E_{y_n}[ \int_{[0,{\tau_n})} e^{-rs} G^2(X_s) d\Gamma^1_s +e^{-r \tau_n}\Lambda^1_{\tau_n-} \bar J^2(X_{\tau_n})] \\
&=\mathbf E_{y_n}[ \int_{[0,{\tau_n})} e^{-rs} G^2(X_s) d\Gamma^1_s +e^{-r \tau_n}\Lambda^1_{\tau_n-}  G^2(X_{\tau_n})] 
\end{align*}  
where $\tau_n:=\tau_{a_n}\wedge \tau_{b_n}$.
$G^2$ and $R^2$ being locally Lipschitz, there exists $\varepsilon>0$ such that  for all $n$ sufficiently large
\[ \forall y \in (a_n,b_n), \; G^2(y) > R^2(y_n) +\varepsilon.\]
We deduce that 
\begin{align*}
 \bar J^2 (y_n) & \geq  (R^2(y_n) + \varepsilon) \mathbf E_{y_n}[ \int_{[0,{\tau_n})} e^{-rs} d\Gamma^1_s +e^{-r \tau_n}\Lambda^1_{\tau_n-}  ]  \\
 &\geq (R^2(y_n) + \varepsilon) \mathbf E_{y_n}[ e^{-r \tau_n} ].
\end{align*}  
We have $\mathbf E_{y_n}[ e^{-r \tau_n} ]=A_n \psi(y_n)+ B_n \phi(y_n)$, where the coefficients $A_n,B_n$ are such that 
\[ A_n \psi(a_n)+ B_n \phi(a_n) = A_n \psi(b_n)+ B_n \phi(b_n) =1.\]
If follows easily that these coefficients are bounded, and therefore that $\mathbf E_{y_n}[ e^{-r \tau_n} ]\rightarrow 1$ by using that $\psi,\phi$ are locally Lipschitz. For $n$ sufficiently large, it contradicts the fact that $ \bar J^2 (y_n)=R^2(y_n)$. We conclude that the functions $\bar {J}^i$ for $i=1,2$ are right-continuous at $x$.\\
\noindent
The cases  ($x>\alpha^1$ and $x\in S^2$) and ($x=\alpha^1$ and $\alpha^1>\alpha^2$ and $x\in S^1$) are similar.\\

\noindent
To conclude that the payoff functions are right-continuous on $\CI$, it remains only to prove that $\bar J^2$ is right-continuous at $x=\alpha^1$ if $x\in S^2$ and $\alpha^2<\alpha^1$. Assume by contradiction that $\bar J^2$ is not right-continuous at $x$. Since $(\alpha^2,\beta)\cap S^1\cap S^2=\emptyset$, we may choose $\varepsilon$ sufficiently small so that $[x,x+\varepsilon) \cap S^1 =\emptyset$ as well as 
\[ G^2-R^2 > \varepsilon \text{ and } \CL R^2 - r R^2 < - \varepsilon \; \text{ on $[x,x+\varepsilon)$} .\]
Note that $\mu^1([x,x+\varepsilon))<\infty$ since $\mu^1$ is locally finite on $\CI \setminus S^1$. 
If $(a,b)$ is a connected component of the open set $(x,x+\varepsilon) \setminus S^2$, so that $a,b\in S^2$, then for $y \in(a,b)$ and $\tau=\tau_a \wedge \tau_b$: 
\[ \bar J^2(y)={J}^2(y,(\mu^1,S^1),(0,S^2))=\mathbf E_y[\int_{[0,\tau)} e^{-rs} G^2(X_s) d\Gamma^1_s + e^{-r\tau}\Lambda^1_{\tau-}R^2(X_{\tau})] \geq  R^2(y).\]
We deduce that 
\begin{align*}
\varepsilon \mathbf E_y[ \Gamma^1_{\tau -}] & \geq \mathbf E_y[\int_{[0,\tau)} e^{-rs} (G^2(X_s)-R^2(X_s)) d\Gamma^1_s] \\
& \geq  \mathbf E_y [ R^2(y) -e^{-r\tau}\Lambda^1_{\tau-}R^2(X_{\tau}) -\int_{[0,\tau)} e^{-rs} R^2(X_s) d\Gamma^1_s] \\
&= \int_0^1 \mathbf E_y [  R^2(y) - e^{- r (\tau \wedge \hat \gamma^1(u,.))} R^2(X_{\tau \wedge \hat \gamma^1(u,.))}) ] du \\
&= -\int_0^1 \mathbf E_y [  \int_0^{\tau \wedge \hat \gamma^1(u,.)} (\CL R^2- r R^2)(X_s) ds ] du \\
& \geq \varepsilon \mathbf E_y [  \int_0^{\tau} \Lambda^1_s ds ]\\
& \geq \varepsilon \mathbf E_y[ \Lambda^1_\tau ].
\end{align*}
Using that $\Gamma^1_{\tau -}=1- \Lambda^1_\tau$ and dividing by $\varepsilon$, we obtain
\[ 1 \geq 2 \mathbf E_y[ \Lambda^1_\tau ] =2 \mathbf E_y[e^{-A^1_{\tau}}] \geq 2 (1- \mathbf E_y[ A^1_\tau ]).\]
Using the formula \eqref{eq:formula_E_A}, we have 
\[ \mathbf E_y[ A^2_\tau ]= \int_a^b 2 (S'(z))^{-1} \Phi_{a,b}(y,z) d \mu^1 (z) \leq C \mu^1((a,b)),\]
for some positive constant $C$ using that $2 (S'(x))^{-1} \Phi_{a,b}(y,x)$ is uniformly bounded on $[x,x+\varepsilon)$. We obtain finally
\[ 1 \geq 2 (1- C \mu^1((a,b))).\]
There exists a sequence of connected components $(a_n,b_n)$ of  $(x,x+\varepsilon) \setminus S^2$ whose length goes to zero. Since $\mu^1$ is locally bounded at $x$, it must be that $\mu^1((a_n,b_n))$ goes to zero, and the previous inequality implies $1\geq 2$, a contradiction. We conclude  that $\bar J^2$ is right-continuous at $x$.\\

\noindent
That the payoff functions are left-continuous follows from similar arguments.\\

\noindent
We conclude that at equilibrium, the payoff functions of both players are continuous. \hfill $\blacksquare$.\\

\renewcommand{\thesection}{Appendix B: Proofs for Section 4}

\renewcommand{\theequation}{B.\arabic{equation}}
\setcounter{equation}{0}

\section{}\label{appendix_on_section4}

The next Lemma solves equation (\ref{barb1})
\begin{lemma}\label{crucialthreshold}
The equation 
\begin{equation} \label{barb11}
 R^i(x_{R^i}) = \frac{ \phi (x_{R^i})}{\phi(x)}G^i(x)
\end{equation}
has a unique  solution $\underline x^j \in (\alpha^i, x_{R^i})$ and,
 $ R_i(x_{R^i}) <\frac{ \phi (x_{R^i})}{\phi(x)} G^i(x)$ on $(\underline x^j, \beta)$.
\end{lemma}

The proof of Lemma \ref{crucialthreshold} is a direct consequence of a  useful and elegant change of variable introduced by Dayanik and Karatzas (2003) which we will use in several proofs.
For each $x \in \mathcal I$, let us define $\zeta(x) := \frac{\phi(x)}{\psi(x)}$ which is strictly decreasing in $x$ and  maps $\mathcal I$ onto $(0,\infty)$. For any  function $g : \mathcal I \to \RR$, define the function $\hat{g}$ by
\begin{eqnarray} \label{changeofvar}
\hat g (y) :=\frac{g}{\psi} \circ \zeta^{-1} (y), \quad y \in(0,\infty).
\end{eqnarray}
 Observe  that $\hat \phi (y) = y$ and $\hat  \psi(y)=1$  for every $y \in (0, \infty)$. A direct computation shows that, if $g \in {\cal C}^2({\cal I})$, then for any $x \in {\cal I}$
\begin{equation} \label{cov}
\hat g''(\zeta(x))= \frac{2 \phi(x)^3}{(\gamma S'(x))^2 \sigma^2(x)}(\mathcal{L}g - rg)(x).
\end{equation}
Thus, we deduce  from assumption A7  that $(\hat G^i)^{''} \leq  0$ everywhere $(\hat G^i)^{''}$ is defined. From assumption A3, $(\hat R^i)^{''}(\zeta(x)) <0$ for every $x\in (\alpha, x_0^i)$  or, equivalently,  $(\hat R^i)^{''}(y) <0$ for every $y\in ( \zeta (x_0^i), \infty)$. This implies that 
$(\hat R^i)^{''}(y) <0$ for every $y\in ( \zeta (x_{R^i}), \infty)$ since $x_{R^i} < x_0^i$. 
We will use also the following remarks. From Lemma \ref{lem_sign} and assumption A6,  we have $G^i > 0$ over $\mathcal I$. Thus, $\hat G^i >0$ over $(0, \infty)$, and (\ref{gp'}) implies
\begin{eqnarray} 
\lim_{y \to 0 }  \hat G^i (y) = 0, \label{lim0}\\
\lim_{y \to \infty} \frac{\hat{G^i}(y)}{y}=0   \label{gp3}.
\end{eqnarray}

\noindent
{\it Proof of Lemma \ref{crucialthreshold}} Let  define $f(x) := \frac{ R^i(x_{R^i}) }{ \phi (x_{R^i})}\phi(x)$. 
Notice that  $f=V_{R^i}\geq R^i$ on $[x_{R^i},\beta)$ and that thanks to the smooth-fit property, $f$ is tangent to $R^i$ at $x_{R^i}$. Applying to $f$ the change of variable formula \eqref{changeofvar}, 
a direct computation shows that, $\underline  x^i$ is a solution to (\ref{barb11}) iff $\zeta (\underline x^i)$ is a solution to 
\begin{equation} \label{barbb1'}
\hat f(y)= \hat G^i(y),
\end{equation}
which is equivalent to 
\[\frac{R^i(x_{R^i})}{\phi (x_{R^i})} y = \hat G^i(y).\] 
Using that $f=V_{R^i}$ on $[x_{R^i},\beta)$, we deduce from assumption A6 that $\hat f(y) < \hat G^i(y)$ on $(0, \zeta(x_{R^i})]$. 
It follows that  Equation (\ref{barb11}) admits a unique  solution $\underline x^j < x_{R^i}$ because $\hat G^i$ is positive concave and satisfies  (\ref{gp3}). Moreover, we have $\frac{\phi(x)}{\phi (x_{R^i})} R^i(x_{R^i}) < G^i(x)$ on $(\underline x^j, \beta)$ and $\frac{\phi(x)}{\phi (x_{R^i})} R^i(x_{R^i}) > G^i(x)$ on $(\alpha,\underline x^j)$. Finally, note that $\alpha^i\leq x_{R^i}\leq x^i_0$ and that $\hat R^i$ is strictly concave on $(\zeta(x^i_0),+\infty)$. Therefore, $\hat f > \hat R^i$ on $(\zeta(x_{R^i}),+\infty)$ since $\hat f$ is linear and tangent to $\hat R^i$ at $\zeta(x_{R^i})$. If $\alpha^i>\alpha$, because $G^i = R^i$ over $(\alpha, \alpha^i]$, it must be that $\alpha^i < \underline x^j$. \hfill $\blacksquare$ \\

\noindent
{\bf Proof of Proposition \ref{armchair}.} We show in the next lemma that the variational systems (\ref{Gonsolin})-(\ref{var15e}) and (\ref{Ermet})-(\ref{var5'e}) admit  a solution. Proposition \ref{armchair} is then a direct consequence  of Theorem \ref{CSmixed} as we explained in the main text.


\begin{lemma} \label{forarmchair}
In the running example, if the firms' liquidation values $l^1 \leq l^2$ are close to each other, $m$ is sufficiently large and $b>0$ then, there exists a constant $a>0$, and two functions $w^1 \in {\cal C}^0((\alpha, \beta))\cap {\cal C}^2 ((\alpha, \beta) \setminus \{  \underline x^1 \})$ and $w^2 \in {\cal C}^0((\alpha, \beta)) \cap {\cal C}^2 ((\alpha,\beta) \setminus \{\underline x^1,  x_{R^1} \})$ solution to the variational systems (\ref{Gonsolin})-(\ref{var15e}) and (\ref{Ermet})-(\ref{var5'e}).
\end{lemma}

\noindent {\bf Proof of Proposition \ref{armchair}.}
Using notations of section \ref{running}, we have $\phi (x) = x^{\rho^-}$, $\psi(x) = x^{\rho^+}$ and
\begin{eqnarray}
V_{R^i}(x) \equiv \sup_{\tau \in {\cal T}_X} \mathbf E_x[e^{-r \tau} R^i(X_\tau)] = \left\{ \begin {matrix} \frac{\phi(x)}{\phi(x_{R^i})} \,(-\frac{ 1}{r - \mu}x_{R^i} +l^i)& \text{if} & x > x_{R^i}, \\  -\frac{1}{r - \mu} x + l^i & \text{if} & x \leq x_{R^i}. \end{matrix} \right. \label{VRex1}
\end{eqnarray}
\noindent
Also,
\begin{eqnarray}
V^i_m (x) & = & \sup_\tau \mathbf E_x [\int_0^\tau e^{-rs} mX_s \,  ds + l^i] \\
& =& \left\{ \begin {matrix} \frac{m}{r- \mu} x + \frac{\phi(x) }{\phi(\alpha^i)}(-\frac{m}{r- \mu} \alpha^i +l^i  ) & \text{if} & x \geq \alpha^i, \\   l^i&  \text{if} & x \leq \alpha^i, \end{matrix} \right.
\end{eqnarray}
with $\alpha^i = \frac{1}{m} x_{R^i}  < x_{R^i}$ and $x_{R^i}=\frac{\rho^- (r - \mu) }{\rho^- - 1} l_i$. Thus,
\begin{eqnarray}
G^i (x) & = & (V^i_m - E)(x) \\
& =& \left\{ \begin {matrix} \frac{m -1}{r- \mu}x   + \frac{\phi(x) }{\phi(\alpha^i)}(-\frac{m}{r- \mu} \alpha^i + l^i ) & \text{if} & x \geq \alpha^i, \\ -\frac{1}{r -\mu}x + l^i & \text{if} & x \leq \alpha^i. \end{matrix} \right.\\
& =& \left\{ \begin {matrix} \frac{m -1}{r- \mu}x   + l^i \frac{\phi(x) }{\phi(\alpha^i)(1-\rho^-)} & \text{if} & x \geq \alpha^i, \\ -\frac{1}{r -\mu}x + l^i & \text{if} & x \leq \alpha^i. \end{matrix} \right.
\end{eqnarray}
\noindent
Equation (\ref{barb1}) that we write under the form
\begin{equation}  \label{exitthreshold}
\frac{\phi(x) }{\phi(x_{R^1})} R^1 (x_{R^1})=G^1(x),
\end{equation}
has a unique solution $\underline x^2 \in (\alpha^1, x_{R^1})$  which is explicit.
Precisely,   (\ref{exitthreshold})  is equivalent to
\[\frac{x^{\rho^-}}{\left( l^1 (r-\mu)\frac{\rho^-}{\rho^--1}  \right)^{\rho^-}} \left[ l^1 - \frac{1}{(r-\mu)}l^1 (r-\mu)\frac{\rho^-}{\rho^--1} \right]=\frac{m -1}{r- \mu}x   + l^1 \frac{x^{\rho^-} m^{\rho^-}}{\left( l^1 (r-\mu)\frac{\rho^-}{\rho^--1}  \right)^{\rho^-}(1-\rho^-)}.\]
\noindent
Multiplying by $\left( l^1 (r-\mu)\frac{\rho^-}{\rho^-
-1}  \right)^{\rho^-}$ and dividing by $xl^1$, we obtain
\[x^{\rho^--1} \frac{1}{1-\rho^-} =\frac{m -1}{r- \mu} (l^1)^{\rho^- - 1} \left( (r-\mu)\frac{\rho^-}{\rho^--1}  \right)^{\rho^-}   +  \frac{x^{\rho^--1} m^{\rho^-}}{(1-\rho^-)}.\]
\[x^{\rho^--1} \frac{1-m^{\rho^-}}{1-\rho^-} =(m -1)(r- \mu)^{\rho^--1} (l^1)^{\rho^- - 1} \left(\frac{\rho^-}{\rho^--1}  \right)^{\rho^-}  .\]
\[x = l^1 (r-\mu) \frac{\rho^-}{\rho^--1}   \left[ \frac{(1-m^{\rho^-})}{(-\rho^-)(m -1)} \right]^{\frac{1}{1-\rho^-}}  .\]
We conclude that 
\[\underline{x}^2 = l^1 (r-\mu) \frac{\rho^-}{\rho^--1} \theta = \theta x_{R^1} .\]
with $\theta:=\theta(m,\rho^-)= \left[ \frac{(1-m^{\rho^-})}{(-\rho^-)(m -1)} \right]^{\frac{1}{1-\rho^-}} \in (1/m,1)$. \\
\noindent
It follows that the function $w^1$ defined by
\begin{eqnarray}
w^1(x) = \left\{ \begin {matrix} \frac{\phi(x)}{\phi(x_{R^1})} \,(-\frac{ 1}{r - \mu}x_{R^1} +l^i)& \text{if} & x > \underline x^2, \\  G^1(x) & \text{if} & x \leq \underline x^2. \end{matrix} \right. \label{VRex1}
\end{eqnarray}
is, by construction, solution to the variariational system (\ref{Gonsolin})-(\ref{var15e}).\\

\noindent
If a solution $w^2$ to (\ref{Ermet})-(\ref{var5'e}) exists, then letting $T^2_x$ denote the unique solution to $\CL u - ru =0$ which is tangent to $R^2$ at $x$, it must be that $w^2=T^2_{\underline x ^2} $ on $(\underline x^2, x_{R^1})$ and that $w^2 =A \phi $ on $(x_{R^1}, \infty)$ for some constant $A$.
Precisely, we have $T^2_{\underline x ^2}= B \psi + C \phi$ with $B,C>0$ given by
\begin{eqnarray*}
B= \frac{\phi'(\underline x^2)(\frac{1}{r-\mu}\underline x^2 - l^2) - \frac{1}{r -\mu} \phi(\underline x^2)}{\psi'(\underline x^2) \phi(\underline x^2)-\psi(\underline x^2) \phi'(\underline x^2)} >0\\
C = \frac{\psi(\underline x^2) \frac{1}{r - \mu} + \psi'(\underline x^2) (- \frac{1}{r - \mu}\underline x^2 + l^2)}{\psi'(\underline x^2) \phi(\underline x^2)-\psi(\underline x^2) \phi'(\underline x^2)}>0
\end{eqnarray*}
Continuity at $x_{R^1}$ yields
$$A = B \frac{\psi(x_{R^1})}{\phi(x_{R^1})} + C,$$
which leads to$$\Delta (w^2)'(x_{R^1}) = B(  \frac{\psi(x_{R^1})}{\phi(x_{R^1})}\phi'(x_{R^1}) -\psi'(x_{R^1})) < 0.$$
\noindent
We deduce that, if  
\begin{equation} \label{ineqforrunning}
 G^2(x_{R^1}) >  T^2_{\underline x ^2}(x_{R^1}) > T^2_{x_{R^2}}(x_{R^1})
\end{equation}
then, $w^2= R^2 1_{(0, \underline x^2]} +T^2_{\underline x ^2} 1_{[\underline x^2, x_{R^1}] }+ 
A \phi 1_{(x_{R^1}, \infty)}$ and, $a = - \frac{1}{2} \,\Delta {w^{2\prime}}(x_{R^1})/(G^2(x_{R^1}) - w^2(x_{R^1}))$ is a solution to the variational system (\ref{Ermet})-(\ref{var5'e}).
The first inequality in (\ref{ineqforrunning}) ensures that $a >0$, the second inequality ensures that $\underline x^2 < x_{R^2}$, and it also implies that $w^2 \geq R^2$ on $[x_{R^1}, \infty)$. The inequality $w^2 \geq R^2$ on $(\underline x^2, x_{R^1})$ follows from the convexity of $T^2_{\underline x ^2}$.
We show below that, if $l^1, l^2$ are close to each other$,$  $m$ is sufficiently large,  and $b>0$ then, (\ref{ineqforrunning}) holds true.\\

\noindent
Letting $\gamma = \rho^+ - \rho^-$, direct computations lead to
\begin{align*}
B & =  \frac{1}{\gamma} \left[- \frac{(1-\rho^-)}{r - \mu} (\underline x^2)^{1-\rho^+} - l^2 \rho^- (\underline x^2)^{-\rho^+})\right] \\
&=  \frac{\rho^-(\underline x^2)^{-\rho^+}}{\gamma} \left[\theta l^1  - l^2 \right]
\end{align*}
\begin{align*}
C & =  \frac{1}{\gamma} \left[\frac{1-\rho^+}{r - \mu} (\underline x^2)^{1-\rho^-} + l^2 \rho^+ (\underline x^2)^{-\rho^-}\right] \\
&=\frac{\rho^+ (\underline x^2)^{-\rho^-}}{\gamma} \left[ l^2 - \theta \frac{\rho^+-1}{\rho^+}\frac{\rho^-}{\rho^--1} l^1\right] 
\end{align*}
\noindent
Using that $\underline{x}^2=\theta x_{R^1}$, we deduce that 
\begin{align*}
T^2_{\underline x ^2}(x_{R^1})& =B (x_{R^1})^{\rho^+}+C (x_{R^1})^{\rho^-}\\
&= \frac{\rho^- \theta^{-\rho^+}}{\gamma} \left[\theta l^1  - l^2 \right]+\frac{\rho^+ \theta^{-\rho^-}}{\gamma} \left[ l^2 - \theta \frac{\rho^+-1}{\rho^+}\frac{\rho^-}{\rho^--1} l^1\right] 
\end{align*}
\noindent
We have $T^2_{x_{R^2}}=\frac{l^2}{(1-\rho^-)\phi(x_{R^2})} \phi$  so that
\begin{align*}
T^2_{x_{R^2}}(x_{R^1})= \frac{l^2}{(1-\rho^-)\phi(x_{R^2})} \phi(x_{R^1})=\frac{l^2}{(1-\rho^-)} \left(\frac{l^1}{l^2}\right)^{\rho^-}.
 \end{align*}
We also have (with equality  if $x_{R^1} \geq \alpha^2$)
\begin{align*}
G^2(x_{R^1}) &\leq  \frac{m -1}{r- \mu}x_{R^1}   + l^2 \frac{\phi(x_{R^1}) }{\phi(\alpha^2)(1-\rho^-)}\\
&= (m-1)\frac{\rho^-}{\rho^--1} l^1 +  l^2 \left( \frac{l^1}{l^2}\right)^{\rho^-}\frac{m^{\rho^-}}{1-\rho^-}.
\end{align*}
\noindent
Therefore, a sufficient condition for (\ref{ineqforrunning}) is
\begin{align*}
(m-1)\frac{\rho^-}{\rho^--1} l^1 &+  l^2 \left( \frac{l^1}{l^2}\right)^{\rho^-}\frac{m^{\rho^-}}{1-\rho^-}  \\
&>  \frac{\rho^- \theta^{-\rho^+}}{\gamma} \left[\theta l^1  - l^2 \right]+\frac{\rho^+ \theta^{-\rho^-}}{\gamma} \left[ l^2 - \theta \frac{\rho^+-1}{\rho^+}\frac{\rho^-}{\rho^--1} l^1\right] \\ 
&> \frac{l^2}{(1-\rho^-)} \left(\frac{l^1}{l^2}\right)^{\rho^-}. 
\end{align*}
\noindent
Let us rewrite these inequalities with $l^1 = l^2$:
\begin{align*}
(m-1)\frac{\rho^-}{\rho^--1}  &+ \frac{m^{\rho^-}}{1-\rho^-}  \\
&>  \frac{\rho^- \theta^{-\rho^+}}{\gamma} \left[\theta  - 1\right]+\frac{\rho^+ \theta^{-\rho^-}}{\gamma} \left[ 1 - \theta \frac{\rho^+-1}{\rho^+}\frac{\rho^-}{\rho^--1} \right] \\ 
&> \frac{1}{(1-\rho^-)}.
\end{align*}
\noindent
We observe that, for $m$ large enough
$$(m-1)\frac{\rho^-}{\rho^--1}  + \frac{m^{\rho^-}}{1-\rho^-}  >
\frac{1}{(1-\rho^-)}.$$
Noticing that $1 - \theta \frac{\rho^+-1}{\rho^+}\frac{\rho^-}{\rho^--1}>1 - \frac{\rho^+-1}{\rho^+}\frac{\rho^-}{\rho^--1}>0$, we also have that
$$
\frac{\rho^- \theta^{-\rho^+}}{\gamma} \left[\theta  - 1\right]+\frac{\rho^+ \theta^{-\rho^-}}{\gamma} \left[ 1 - \theta \frac{\rho^+-1}{\rho^+}\frac{\rho^-}{\rho^--1} \right] \\ 
> \frac{1}{(1-\rho^-)}.$$
\noindent
We finally show that, for $m$ large enough and $b>0$,
\begin{equation} \label{lastineq}
(m-1)\frac{\rho^-}{\rho^--1}  + \frac{m^{\rho^-}}{1-\rho^-}  
>  \frac{\rho^- \theta^{-\rho^+}}{\gamma} \left[\theta  - 1\right]+\frac{\rho^+ \theta^{-\rho^-}}{\gamma} \left[ 1 - \theta \frac{\rho^+-1}{\rho^+}\frac{\rho^-}{\rho^--1} \right].
\end{equation}
Recall that $\theta= \left[ \frac{(1-m^{\rho^-})}{(-\rho^-)(m -1)} \right]^{\frac{1}{1-\rho^-}}$. We deduce that a sufficient condition for (\ref{lastineq}) is that, for $m$ large enough
$$m > \frac{\rho^- m^{\frac{\rho^+}{1 - \rho^-}}}{\gamma} (m^{\frac{1}{\rho^- -1}} - 1)  +\frac{\rho^+ m^{\frac{\rho^-}{1 - \rho^-}}}{\gamma} ( 1 - m^{\frac{1}{\rho^- -1}} \frac{\rho^+-1}{\rho^+}\frac{\rho^-}{\rho^--1}).$$
\noindent
The latter inequality holds true for $m $ large when $\rho^+ + \rho^- < 1 $ which is equivalent to $b >0$. A simple continuity argument ends the proof of the Lemma. \hfill $\blacksquare$

\renewcommand{\thesection}{Appendix C: Proofs for Section 5}

\renewcommand{\theequation}{C.\arabic{equation}}
\setcounter{equation}{0}
\section{}\label{appendix_on_section5}

The proof of  Theorem \ref{CNmixed} is mainly based on Proposition \ref{geneprop} and on assertions (i), (ii) and (iii) of the following simple lemma. We will use later assertion (iv) in the proof of Theorem \ref{CSmixed}.

\begin{lemma}\label{one_point}
Let $u$ be a $C^2$ function defined on an open interval $(a,b) \subset \CI$ which satisfies 
\[ \forall x\in (a,b),\; {\cal L} u (x) - ru (x) =0.\]
Then, we have 
\begin{itemize}
\item[(i)] If $b=\beta$, $u(\beta-)=0$, $u(a+)=R^i(a)$ and $u(x)\geq V_{R^i}(x)$ for all $x\in(a,\beta)$, then $a= x_{R^i}$.  
\item[(ii)] If $u(x)\geq V_{R^i}(x)$ for all $x\in(a,b)$ then $\{x \in (a,b) \,|\, u(x)=R^i(x)\}$ contains at most one point.
\item[(iii)] If there exist two points $a',b'\in(a,b)$ with $b'\leq x_{R^i}$ and $u(a')=R^i(a')$, and $u(b')=R^i(b')$, then $u(x)<R^i(x)$ for all $x\in (a',b')$. 
\item[(iv)] Assume that $\alpha< a \leq x_{R^i}$, $u \geq R^i$ on $(a,b)$, $u(a)=R^i(a)$ and $u'(a+)>(R^i)'(a)$. Then, for every $\varepsilon>0$ sufficiently small, the function $f$ solution of $\CL f - rf =0$ on $(a-\varepsilon,a+\varepsilon)$ with $f(a-\varepsilon)=R^i(a)$ and $f(a+\varepsilon)=u(a+\varepsilon)$ satisfies $f(a)> u(a)$.   
\end{itemize}
\end{lemma}

\noindent
{\bf Proof of Lemma \ref{one_point}.}
Recall  that a solution $u\in {\cal C}^2((a,b))$ of the equation ${\cal L} u(x) - ru(x) = 0$  writes under the form $u(x) = A \psi(x)  + B \phi(x)$ where  $A$ and $B$ are real numbers. Using the change of variables  (\ref{changeofvar}), the inequation  $u(x) \geq V_{R^i}(x)$ on $(a,b)$ is equivalent to  $\hat u(z)=A z + B \geq   \hat V_{R^i}(z)$ for all  $z\in (\zeta(b), \zeta(a))$.

Recall that $\hat V_{R^i}$ is $\CC^1$, that $\hat V_{R^i}(z)=C^i z > \hat R^i(z)$ for all  $z\in (0, \zeta(x_{R^i}))$ for some $C^i>0$, and that $\hat V_{R^i}=\hat R^i$ and is $\CC^2$ and strictly concave on $[\zeta (x_{R^i}), \infty)$.

Proof of $(i)$. The assumption $u(\beta-)=0$ implies $B=0$, and thus $\hat u(0+)=0$. The assumption $\hat u \geq \hat V_{R^i}$ on $(0, \zeta(a))$, implies that $A \geq C^i$. If this inequality is strict, then $A z > C^i z \geq \hat R^i (z)$ for all $z>0$, since $\hat V_{R^i}$ is concave, which would contradict the assumption $\hat u(\zeta(a)-)=A \zeta(a)=\hat R^i (\zeta(a))$. Therefore $A=C^i$ and from the properties of $\hat V_{R^i}$, the unique solution of $A z= \hat R^i(z)$ is $\zeta(x_{R^i})$.
   
Proof of $(ii)$.
Note that $V_{R^i}>R^i$ on $(x_{R^i},\beta)$, so that $\hat V_{R^i} > \hat R^i$ on $(0,\zeta(x_{R^i}))$.
Therefore, if there exists $\bar z\in (\zeta(b),\zeta(a))$ such that $\hat{u}(\bar z)=\hat R^i(\bar z)$, it must be that $\bar z\geq  \zeta(x_{R^i})$. 
In such a case, $\hat{u}$ is a tangent to the concave $\CC^1$ map $\hat{V}_{R^i}$ at $\bar z$.
On the interval $[\zeta(x_{R^i}),+\infty)$, $\hat V_{R^i}=\hat R^i$ is strictly concave and thus $\hat{u}(z)>\hat R^i(z)$ for all $z \neq \bar z$ in $[\zeta(x_{R^i}),+\infty) \cap (\zeta(b),\zeta(a))$, which concludes the proof.

Proof of $(iii)$. $\hat u$ is an affine map on the interval $[\zeta(b'),\zeta(a')]$ which is equal to $\hat R^i$ at both boundaries. The inequality follows directly from the fact that $\hat R^i$ is strictly concave on $[\zeta (x_{R^i}), \infty)$ given that $b' \leq x_{R^i}$.

Proof of $(iv)$. We have $\hat u (z)= A z+ B$ for some constants $A,B$ on $[\zeta(b),\zeta(a)]$. Since $a\leq x_{R^i}$, the map $\hat R^i$ is concave on $[\zeta(a),+\infty)$. 
A direct computation shows that 
\[ A=\hat u'(\zeta(a)-)= \frac{\psi(a)u'(a+)- \psi'(a)u(a)}{\psi(a)^2 \zeta'(a)} , \; (\hat R^i)'(\zeta(a))=\frac{\psi(a)(R^i)'(a) - \psi'(a)R^i(a)}{\psi(a)^2\zeta'(a)}, \]
so that 
\[(\hat R^i)'(\zeta(a))-\hat u'(\zeta(a)-)=  \frac{(R^i)'(a)- u'(a+)}{\psi(a) \zeta'(a)} >0,\]
using that $u(a)=R^i(a)$, $\zeta'(a)<0$ and $u'(a+)>(R^i)'(a)$. 
We deduce that for $\varepsilon$ small enough, $\hat R^i (\zeta(a-\varepsilon)) > A \zeta(a-\varepsilon) +B$.
The function $f$ satisfying $\CL f - r f=0$ on $(a-\varepsilon,a+\varepsilon)$, 
$f(a-\varepsilon)=R^i(a-\varepsilon)$ and $f(a+\varepsilon)=u(a+\varepsilon)$ is such that $\hat f (z)=A'z+B'$ for some constants $A',B'$ with $\hat f (\zeta(a-\varepsilon))=\hat R^i(\zeta(a-\varepsilon))$ and $\hat f (\zeta(a+\varepsilon))=\hat u (\zeta(a+\varepsilon))$. We deduce that 
\[ A \zeta(a+\varepsilon)+B= A'\zeta(a+\varepsilon) +B \text{ and } A'\zeta(a-\varepsilon) +B'> A \zeta(a-\varepsilon) +B, \]
which imply $A'\zeta(a) +B'> A \zeta(a) +B$ or equivalently $f(a)>u(a)$.  \hfill $\blacksquare$\\

\noindent \textbf{Proof of Theorem \ref{CNmixed}.}
For $i \in \{1, 2 \}$, let us define $q^{ i}_1 := \max ( \mbox{supp}(\mu^i))$ with the convention $\max\emptyset = \alpha$. We assume without loss of generality that $q^{ 1}_1 \geq q^{
 2}_1$.
From  assertion $(iv)$ of Proposition \ref{geneprop}, we have that $x_{R^i} \geq q^{
 i}_1 \vee s^i$.
The proof consists of four steps and uses repeatedly Lemma \ref{one_point}. 
In the following,  $T^i_q$ denotes the curve solution to $\CL u - ru =0$ and tangent to $R^i$ at $q$ where   $q\leq x_{R^i}$.\\

\noindent
{\bf Step 1.} We show that
\begin{itemize}
\item[(i)]
If $q^{ 1}_1 > s^1 \vee s^2$ then $q^{ 1}_1= x_{R^1} >q^{ 2}_1 $, and $x_{R^1}$  is an isolated point of the support of $\mu^1$.
\item[(ii)]
If   $q^{ 2}_1 > s^1 \vee s^2 $ then,
$q^{ 2}_1$ is an isolated point of the support of $\mu^2$. 
\end{itemize}

Note that the assumption $q_1^1>s^1\vee s^2$ implies $q_1^1>\alpha$. 
From Proposition \ref{geneprop}-$(v)$, the strategy $(0, S^1) $ is a PBR to strategy  $(\mu^2, S^2)$. Therefore, $\tau_{S^1}$ is an optimal solution to
 $ \bar J^1 ( x, (\mu^2, S^2))= \sup_{\tau^1 \in {\cal T}} J^1 (x, \tau^1, (\mu^2, S^2))$.
 Letting $y=q^{ 2}_1\vee s^1\vee s^2$, we have $\tau_{S^1} \geq \tau_{y}$ $\mathbf{P}_x$-almost surely for all $x \in (y, \beta)$. 
We deduce that the best reply value function $\bar J^1$ satisfies,  for any $x \in (y,\beta)$,
\begin{align*}
\bar J^1(x):= \bar J^1 (x,(\mu^2, S^2)) &= J^1(x, (\mu^1, S^1), (\mu^2, S^2)) = J^1(x, (0, S^1), (\mu^2, S^2)) \\
& = \mathbf E_x \left [ e^{- r \tau_{y}} \bar J^1(X_{\tau_{y}}) \right ]=  \mathbf E_x \left [ e^{- r \tau_{y}} \bar J^1(y) \right ],
\end{align*}
where the strong Markov Property (\ref{MarkovpropertyJ}) yields the penultimate equality.
It follows  from (\ref{laplace}) that 
\begin{equation} \label{sym1} \lim_{x\rightarrow \beta}\bar J^1(x) = 0 \;\;\mbox{and} \;\; \CL \bar J^1(x) - r \bar J^1(x) = 0 \;\; \mbox{on} \;\; (y,\beta).
\end{equation}
Furthermore, we have  that
\begin{equation} \label{cons2}
\bar J^1(x) \geq V_{R^1}(x) \; \mbox{ on } \; (y,\beta), \; \mbox{and}, \;
\bar J^1 (q^{ 1}_1)= R^1( q^{ 1}_1),
\end{equation}
where the last equality follows from assertions (ii) and (iii) of Proposition \ref{geneprop}.
Finally, Lemma \ref{one_point}-(i) together with \eqref{sym1} and \eqref{cons2}  imply that $q^{ 1}_1 = x_{R^1}$. 

By way of contradiction, assume  that $q^{ 1}_1 = q^{ 2}_1$. Applying the same arguments than above and exchanging the roles of the two players, we obtain $q^{ 2}_1 = x_{R^2}\neq x_{R^1}$, a contradiction.
Therefore, $q^{ 1}_1 > q^{ 2}_1$ and thus $q^{ 1}_1>y$. 
Lemma \ref{one_point}-(ii) together with \eqref{sym1} and \eqref{cons2}  imply
\begin{equation} \label{cons3}
q^{ 1}_1 = x_{R^1}, \; \mbox{and,  } \; (y,\beta) \cap \{ x \in (\alpha, \beta) \, | \, \bar J^1(x) = R^1(x) \} = \{x_{R^1} \}.
\end{equation} 
Thus,  $q^{ 1}_1 =x_{R^1}$ is an isolated point of the support of $\mu^1$.  Observe  that $\bar J^1 = T^1_{q_1^{ 1}}$ on $(y,\beta)$.

To prove (ii), we first define for $i=1,2$
\[ q^{ 1}_2= \max (\mbox{supp}(\mu^1)\setminus \{q^{ 1}_1\} )< q^{ 1}_1,\]
where the strict inequality follows from (i). Note that the assumption   $q^{ 2}_1 > s^1 \vee s^2 $ implies $y=q^{ 2}_1 > s^1\vee s^2$. It follows  that  $q^{ 2}_1 > q^{ 1}_2$. Indeed, if $q^{ 1}_1 > q^{ 1}_2 \geq q^{2}_1$, then the strong Markov property implies that, for any $x \in (q^{1}_2, q^{ 1}_1)$, 
 $\bar J^1(x)  = \mathbf E_x[e^{- r \tau_{q^{ 1}_1} \wedge  \tau_{q^{ 1}_2}} \bar J^1 (X_{ \tau_{q^{ 1}_1} \wedge  \tau_{q^{ 1}_2}})]$. We deduce from  a standard computation  that $\bar J^1$ satisfies ${\cal L}\bar J^1 - r\bar J^1 = 0$ on $(q^{ 1}_2, q^{ 1}_1)$. We also  have
 $\bar J^1 (q^{1}_1) = R^1(q^{ 1}_1) $, and $ \bar J^1 (q^{1}_2)= R^1( q^{ 1}_2)$. It follows from  Lemma \ref{one_point} - (iii) that   $\bar J^1(x)  < R^1(x)$ on $(q^{1}_2, q^{ 1}_2)$, thus a contradiction.

Thus, we have that  $q^{1}_1 > q^{2}_1 > q^{1}_2$. Then, given $x \in (q^{1}_2, q^{1}_1)$,  the strong Markov property leads to  $\bar J^2(x) = \mathbf E_x[e^{- r \tau_{q^{ 1}_1} \wedge  \tau_{q^{ 1}_2}} \bar J^2 (X_{ \tau_{q^{1}_1} \wedge  \tau_{q^{1}_2}})]$ with
\begin{equation}
 \bar J^2 (x) \geq  V_{R^2}(x) \; \mbox{on} \;  (q^{1}_2, q^{ 1}_1), \; \mbox{and} \;  
\bar J^2 (q^{2}_1) = R^2( q^{2}_1).
\end{equation}
We  apply Lemma \ref{one_point}-(ii) to $\bar J^2$ and we
 get that
\begin{equation} 
\bar J^2 (x)> R^2(x) \; \mbox{on }  (q^{1}_2, q^{1}_1) \;\setminus \{q^2_1 \}.
\end{equation}
Thus,  $q^{ 2}_1$ is an isolated point of the support of $\mu^2$.
Observe that $\bar J^2 = T^2_{q^{2}_1}$ on $(q^{ 1}_2, q^{1}_1)$.

Define then 
\[ q^{ 2}_2= \max (\mbox{supp}(\mu^2)\setminus \{q^{ 2}_1\} )< q^{2}_1.\]

\noindent
{\bf  Step 2.} 
 We show that
 \begin{itemize}
\item[(i)]
If $ q^{ 1}_2 > s^1 \vee s^2$ then, $q^{ 1}_2$ is an isolated point of the support of $\mu^1$. 
\item[(ii)]
If $ q^{ 2}_2 > s^1 \vee s^2$ then, 
$q^{ 2}_2 $ is an isolated point  of the support of $\mu^2$.
 \end{itemize}

Let us consider that $q^{ 1}_2 > s^1 \vee s^2$ and let us show that $q^{ 1}_2$ is an isolated point of the support of $\mu^1$. The same reasoning than in step 1-$(ii)$ applies, and the inequality  
 $q^{ 1}_2 >  q^{ 2}_2$ is proven by contradiction using Lemma \ref{one_point}-(iii).  
Then, we have
\begin{equation} \label{cons4}
 \bar J^1 (x) = \mathbf E_x[e^{- r \tau_{q^{ 2}_1} \wedge  \tau_{q^{ 2}_2}} \bar J^1 (X_{ \tau_{q^{ 2}_1} \wedge  \tau_{q^{ 2}_2}})] \geq V_{R^1}(x) \; \mbox{ on } (q^{ 2}_2, q^{ 2}_1) \; \mbox{and} \;  \bar J^1 (q^{ 1}_2) = R^1(q^{ 1}_2 ),  
 \end{equation}
and Lemma \ref{one_point}-(ii) implies
\begin{equation} \label{cons5}
\bar J^1 (x)> R^1(x) \; \mbox{on } \;(q^{ 2}_2, q^{ 2}_1) \setminus \{q^{ 1}_2 \}.
\end{equation}
Thus,  $q^{ 1}_2$ is an isolated point of the support of $\mu^1$.

Observe that $\bar J^1 = T^1_{q^{ 1}_2}$ on $( q^{ 2}_2, q^{ 2}_1)$. Because $\bar J^1$ is continuous from Proposition \ref{lem:continuity}, 
$T^1_{q^{ 1}_1}$ and  $T^1_{q^{ 1}_2}$ intersects at $q^{ 2}_1$ and we have that $\bar J^1 (q^{ 2}_1) = T^1_{q^{ 1}_1} (q^{ 2}_1) =T^1_{q^{ 1}_2} (q^{ 2}_1)$.

Let us assume that $ q^{ 2}_2 > s^1 \vee s^2$. An analogous reasoning yields that  $ q^{ 2}_2 $ is an isolated point  of the support of $\mu^2$. The curve  $T^2_{q^{ 2}_1}$ and $T^2_{q^{ 2}_2}$ intersects at $q^{ 1}_2$ and we have 
$\bar J^2 (q^{ 1}_2) = T^2_{q^{ 2}_1} (q^{ 1}_2) =T^2_{q^{ 2}_2} (q^{ 1}_2)$.\\

\noindent
{\bf Step 3.} Proceeding by induction, let us define for all $k\geq 1$ and $i=1,2$
\[ q^{ i}_{k+1} = \max ( \mbox{supp}(\mu_i) \setminus \{q^{ i}_1,q^{ i}_2,...,q^{ i}_k\}) .\] 
Assume that $(q^{ 1}_l)_{1\leq l \leq  k}$ and $(q^{ 2}_l)_{1\leq l \leq  {k}}$ are such that  $q^{ 1}_1 >q^{ 2}_1>q^{ 1}_2>q^{ 2}_2>...>q^{ 1}_{k-1}>q^{ 2}_{k-1}>q^{ 1}_k >q^{ 2}_k$ and that:
\begin{itemize} 
\item for any  $2 \leq l \leq k$,  $\bar J^1 = T^1_{q^{ 1}_l}$ on $(q^{ 2}_{l}, q^{ 2}_{l-1})$ and, $T^1_{q^{ 1}_{l - 1}}$ and $T^1_{q^{ 1}_{l}}$  intersects at $q^{ 2}_{l-1}$.
\item for any  $2 \leq l \leq k$,   $\bar J^2 = T^2_{q^{ 2}_l}$ on $(q^{ 1}_{l+1}, q^{ 1}_{l})$ and, 
$T^2_{q^{ 2}_{l - 1}}$ and $T^2_{q^{ 2}_{l}}$  intersects at $q^{ 1}_{l}$.
\end{itemize}

The same arguments as in step 2 show that
\begin{itemize}
\item[(i)]
 If $q^{ 1}_{k+1} > s^1 \vee s^2 $ then,
$q^{ 1}_{k+1}$ is an isolated point of the support of $\mu^1$ and, $T^1_{q^{ 1}_{k}}$ and $T^1_{q^{ 1}_{k+1}}$ intersects at $q^{ 2}_{k}$.
\item[(ii)]
 If $q^{ 2}_{k+1} > s^1 \vee s^2$ then, $q^{ 2}_{k+1}$ is an isolated point of the support of $\mu^2$ and, $T^2_{q^{ 2}_{k}} $ and $T^2_{q^{ 2}_{k+1}}$ intersects at $q^{ 1}_{k+1}$.
\end{itemize}

 \noindent
{\bf Step 4.} If for some $k \geq 1$ and $i=1,2$, we have $q^{ i}_{k} \leq s^1\vee s^2$, then applying the induction step a finite number of times, we conclude that the restrictions of $\mu^1$ and $\mu^2$ to $(s^1\vee s^2, \beta)$ have finite support, which concludes the proof.

Otherwise, we may apply the induction step infinitely many times, and we obtain two decreasing sequences $(q^{ 1}_n)_{n\geq 1}$,  $(q^{ 2}_n)_{n\geq 1}$ with 
$q^{ 1}_1 > q^{ 2}_1 >q^{ 1}_2 > q^{ 2}_2>...>q^{ 1}_n > q^{ 2}_n > ...$
where 
$q^{ 1}_1 = x_{R^1}$, $q^{ 2}_1 <x_{R^1}$ is given and where for any $n \geq 2$, 
\begin{itemize}
\item
 $T^1_{q^{ 1}_n}$ and $T^1_{q^{ 1}_{n-1}}$ intersects at $q^{ 2}_{n- 1}$.
\item
$T^2_{q^{ 2}_{n-1}}$ and $T^2_{q^{ 2}_n}$ intersects at $q^{ 1}_{n}$.
\end{itemize}
Lemma \ref{conv} below  shows that any decreasing sequences $(q^{ 1}_n)_{n\geq 1}$,  $(q^{ 2}_n)_{n\geq 1}$ satisfying these properties converge to $\alpha$. Thus, if for any $n \geq 1$, $q_n^i > s^1\vee s^2$ then,  $S^1=S^2=\emptyset$ and, $\mbox{supp} \mu^i =\{q^{ i}_{k},  \; k \geq 1\}$, for $i=1,2$.  
This concludes  the proof of Theorem \ref{CNmixed}. \hfill $\blacksquare$

\begin{lemma}\label{conv}
The two decreasing sequences  $(q^{ 1}_n)_{n \in \NN}$ and $(q^{ 2}_n)_{n \in \NN}$ converge to $\alpha$. 
\end{lemma}

\bigskip

\noindent
{\bf Proof of Lemma \ref{conv}.}

The proof makes extensive use of change of variable  \eqref{changeofvar}.
As a preliminary, observe that, using the change of variable \eqref{changeofvar}, assumption {\bf A8}
implies that $(\hat{R}^i)''$ is locally Lipschitz.\\

\noindent 
Let us assume by contradiction that $\lim q_{ i}^n = \bar q > \alpha$.\\

\noindent 
Let $T^i(q)$ denote the curve solution to $\CL u - ru =0$ and tangent to $R^i$ at $q$ for $q\geq x_{R^i}$.
This curve is above $R^i$ for $x> x_{R^i}$ and is equal to $A(q) \psi+B(q) \phi$ for some positive coefficients $A(q),B(q)$. \\

\noindent 
 Define $z^i= \zeta(x_{R^i})$ and note that for any $z \geq z^i$, then $\hat{T}^i_z= \widehat{ T^i_{\zeta^{-1}(z)}}$ is the affine map  tangent line to the function $\hat{R}^i$ at $z$ given by: 
\[  \hat{T}^i_z(y)= A (\zeta^{-1}(z)) y+  B (\zeta^{-1}(z))= \hat R^i(z)+ \hat R'^i(z)(y-z).\]\\
\noindent
We define $y_{2n}=\zeta(q^2_n)$ and $y_{2n-1}=\zeta(q^1_n)$ for $n\geq 1$.\\
\noindent 
As a consequence of the proof of Theorem \ref{CNmixed}, the following recursive relations are satisfied: 
\[ \hat R^1( y_{2n-1} )+ (\hat R^1)'(y_{2n-1})(y_{2n} - y_{2n-1})=\hat R^1(y_{2n+1})+ (\hat R^1)'(y_{2n+1})(y_{2n} - y_{2n+1}).\] 
\[ \hat R^2( y_{2n} )+ (\hat R^2)'(y_{2n})(y_{2n+1} - y_{2n})=\hat R^2(y_{2n+2})+ (\hat R^2)'(y_{2n+2})(y_{2n+1} - y_{2n+2}).\] 
In words,
\begin{itemize}
\item
 The tangent line $\hat T^1_{y_{2n-1}}$ to the function $\hat R^1$ at $y_{2n -1}$ and the tangent line $\hat T^1_{y_{2n+1}}$ to the function $\hat R^1$ at $y_{2n +1}$ intersect at $y_{2n}$.
 \item
  The tangent line $\hat T^2_{y_{2n}}$ to the function $\hat R^2$ at $y_{2n }$ and the tangent line $\hat T^2_{y_{2n+2}}$ to the function $\hat R^2$ at $y_{2n +2}$ intersect at $y_{2n+1}$.
\end{itemize}
\noindent 
The above equations for $i=1,2$ can be written as (with $x<y<z$ three consecutive terms of the sequence $(y_n)$)
\begin{equation}\label{Phi_equiv}
 \hat R^i(x)+ (\hat R^i)'(x)(y-x) - \hat R^i(y)= \hat R^i(z) + (\hat R^i)'(z)(y-z) - \hat R^i(y)
\end{equation}
Using Taylor's theorem with integral remainder, \eqref{Phi_equiv} is equivalent to 
\begin{equation}\label{eq_integral_R}
- \int_x^y (y-s) (\hat R^i)''(s) ds= -\int_y^z (s-y) (\hat R^i)''(s)ds. 
\end{equation}
Note that $(\hat R^i)''<0$ on $[y_1,\infty)$, so that the right-hand side of \eqref{eq_integral_R} is increasing in $z$. Therefore, given $y_1\leq x <y$, if a solution $z>y$ exists, it is unique.\\

\noindent
By assumption, we have $y_n \rightarrow \bar y= \zeta(\bar q)<\infty$.\\

\noindent 
Since $(\hat R^i)''$ is locally Lipschitz, there exists $K$ such that for all $n$ and $i=1,2$,  
\[ \forall s,y \in [y_1,\bar y],  \;  |(\hat R^i)''(s) - (\hat R^i)''(y)|\leq  K |s-y| .\] 
We deduce that
\begin{align*}
- \int_x^y (y-s) (\hat R^i)''(s) ds \geq - (\hat R^i)''(y) \frac{(y-x)^2}{2} -K \frac{(y-x)^3}{3}
\end{align*}
\begin{align*}
- \int_y^z (s-y) (\hat R^i)''(s) ds \leq - (\hat R^i)''(y) \frac{(z-y)^2}{2} +K \frac{(z-y)^3}{3}
\end{align*}
From equation \eqref{eq_integral_R}, we obtain
\begin{align*}
(z-y)^2 + 2K \frac{(z-y)^3}{3|(\hat R^i)''(y)|} \geq (y-x)^2 - 2K \frac{(y-x)^3}{3|(\hat R^i)''(y)|}
\end{align*}
Let $C$ such that for all $y \in [y_1,\bar y]$ and $i=1,2$
\[ \frac{2K}{3|(\hat R^i)''(y)|} \leq C.\] 
Then, we have
\begin{align*}
(z-y)^2 + C(z-y)^3 \geq (y-x)^2 - C(y-x)^3.
\end{align*}
Define $u_n=y_{n+1}-y_n$ and 
\[ g(u)=u^2-Cu^3, \;  h(u)=u^2+Cu^3\]
so that the previous analysis leads to $h(u_{n+1})\geq g(u_n)$.\\

\noindent 
By assumption $y_1 +\sum_{n\geq 1} u_n =\bar y <\infty$, which implies that $u_n \rightarrow 0$. 
Therefore, for $n$ sufficiently large $g(u_n)>0$ and we have  $u_{n+1} \geq h^{-1}(\phi(u_n))$,
where $h^{-1}$ denotes the inverse of $h$ restricted to $[0,\infty)$.\\

\noindent 
We have $h^{-1}(z)=\sqrt{z} - \frac{C}{2}z + o(z)$, so that $h^{-1}(g(u))= u -Cu^2 + o(u^2)$. 
We deduce that 
\[ u_{n+1} \geq  u_n -C u_n^2 +o(u_n^2).\]
It follows that 
\[ \frac{1}{u_{n+1}} - \frac{1}{u_n} \leq \frac{1}{u_n}\left[\frac{1}{1- Cu_n +o(u_n)} - 1  \right] =   C + o(1).\]
We obtain 
\[ \frac{1}{u_{n}} = \frac{1}{u_1} + \sum_{k=1}^{n-1} \left(\frac{1}{u_{k+1}} - \frac{1}{u_k}  \right)  \leq  n C +o(n),\]
and finally 
\[ u_n \geq  \frac{1}{n C} + o(\frac{1}{n}),\]
which is a contradiction as $y_1 +\sum_{n\geq 1} u_n =\bar y <\infty$. \hfill $\blacksquare$

\bigskip

\noindent
{\bf Proof of Theorem \ref{CSmixed}.} Let us assume that $$((\mu^1, S^1), (\mu^2, S^2)) = \left (\left (\sum_{i =1}^n a_i \delta_{q^{ 1}_i}, \emptyset \right), \left (\sum_{i = 1}^{n- 1} b_i \delta_{q^{ 2}_i}, (\alpha,s^{ 2}]\right )\right )$$ is a mixed strategy MPE. Let $\Lambda^i$ for $i=1,2$ denote the {\it ccdf } associated to the strategy $(\mu^i,S^i)$, so that 
\[ 
\Lambda^1_t= e^{- \sum_{i =1}^{n} a_i L_t^{q^{ 1}_i}}, \; 
\Lambda^2_t = e^{- \sum_{i =1}^{n-1} b_i L_t^{q^{ 2}_i}}1_{t < \tau^2_*}, 
\]
where for the sake of lighter notation $\tau^2_*=\tau_{(\alpha,s^2]}$.\\
\noindent
Let us consider $\bar J^2$, the best reply value function to strategy $(\mu^1, S^1)$. We will show that $\bar J^2$ is a solution to the system $\CV\CS^2$.  
From Proposition \ref{geneprop} we have that, $V_{R^2} \leq \bar J^2 \leq G^2$, thus  $\bar J^2 (x)  \geq R^2(x)$ on ${\cal I}$.
From Proposition \ref{lem:continuity}, $\bar J^2$ is continuous on ${\cal I}$.
From the proof of Proposition \ref{CNmixed}, we have that $\lim_{x \rightarrow \beta} \bar J^2(x)=0$, that $\CL \bar J^2 - r\bar J^2 =0$ on $(q^1_1,\beta)$ and on $(s^2, q^1_n)$, and that $\bar J^2 = T^2_{q^{ }_l}$  on $(q^{ 1}_{l+1}, q^{ 1}_{l})$ for $1 \leq l \leq n-1$  where $T^2_q$ denotes the curve solution to ${\cal L } u - r u = 0$ and tangent to $R^2$ at $q$. Therefore, if
$((\mu^1,S^1), (\mu^2, S^2))$ is a MPE then, $\bar J^2$ satifies (\ref{var11g}), (\ref{var13g}), (\ref{var15g}),\eqref{var16g}.
Also, $\bar J^2(x) = R^2(x)$ for $x \leq s^{ 2}$ because the payoff of player 2 is $R^2$ on $S^2$. Therefore, (\ref{var11g'}) is satisfied. 
Therefore, 
 $\bar J^2$ is in ${\cal C}^0({\cal I})  \cap {\cal C}^2({\cal I} \setminus (\{(q^{ 1}_i)_{1 \leq i \leq n}\} \cup  \{s^{ 2} \} ))$, and also $|(\bar {J^2})' (x+)| < \infty $ and  $|(\bar{ J^2})' (x-)| < \infty $ for $x\in \{(q^{ 1}_i)_{1 \leq i \leq n}\} \cup  \{s^{ 2} \}$.
  It remains to prove the conditions \eqref{var14g} and \eqref{var14g'} to show that $\bar J^2$ is a solution of the system $\CV\CS^2$. \\
  
  Let us  prove that $\bar J^2$ satisfies  (\ref{var14g}). 
Since $\bar J^2 \geq R^2$ with equality at $s^2$, it must be that $(\bar J^2)'(s^2) \geq (R^2)'(s^2)$. Assume by contradiction that this inequality is strict. Let us consider the stopping time $T_{\varepsilon} := \inf \{t \geq 0 \, | \, X_t \notin (s^{ 2} - \varepsilon, s^{ 2} +\varepsilon) \}$ where $\varepsilon > 0$ is such that $\alpha < s^2 - \varepsilon < s^2+\varepsilon < q^1_n$. Define $f(x)=\mathbf E_x[ e^{-r T_\varepsilon} \bar J^2(X_{T_\varepsilon})]$ for $x\in [s^2 - \varepsilon, s^2 +\varepsilon]$. Noting that $\tau^2_*$ is a best reply to $(\mu^1,S^1)$ and applying the strong Markov property, we obtain that $f$ is the payoff of player $2$ against $(\mu^1,S^2)$ when using the (non Markovian) stopping time $T_\varepsilon+ \tau_{(\alpha,s^2]} \circ \theta_{T_\varepsilon}$ (i.e. waiting up to $T_\varepsilon$ and then stopping the first time $X$ is below $s^2$  in the continuation game). Applying Lemma \ref{one_point}-$(iv)$ with $i=2$, $a=s^2$, $b=q^1_n$ and  $u=\bar J^2$, we deduce that for sufficiently small $\varepsilon$, $f(s^2)> \bar J^2(s^2)$, which contradicts the fact that $(\mu^2,S^2)$ is a PBR to $(\mu^1,S^1)$. 
Thus  $\bar J^2$ satisfies  (\ref{var14g}).\\
  
That $\bar J^2$ satisfies (\ref{var14g'}) follows from
  Lemma  \ref{CN} below.
 \begin{lemma} \label{CN}
If $((\mu^1, S^1), (\mu^2, S^2)) = \left (\left (\sum_{i =1}^n a_i \delta_{q^{ 1}_i}, \emptyset \right), \left (\sum_{i = 1}^{n- 1} b_i \delta_{q^{ 2}_i}, (\alpha,s^{ 2}])\right )\right )$ is a mixed strategy MPE then, for every stopping time $\tau \in {\cal T}$ and every $x \in {\cal I}$,
\begin{eqnarray}
\bar J^2(x)  & = 
 \mathbf E_x \Big[ &\sum_{i= 1}^{n} \int_{[0,{\tau}\wedge \tau^2_*)} e^{-rs} G^2(q^{ 1}_i) \Lambda_s^1  a_i dL_s^{q^{ 1}_i} + e^{-r \tau^2_*} \Lambda^1_{\tau^2_*} R^2(X_{\tau^2_*}) 1_{\tau^2_* < {\tau}}   \nonumber \\
& & + e^{-r{\tau }}  \Lambda_{\tau}^1 \bar J^2(X_{\tau })1_{\tau \leq\tau^2_*} \Big] \label{fortanaka33} \\
&  = \mathbf E_x \Big[ & \sum_{i = 1}^{n } \int_{[0, \tau\wedge \tau^2_*)}
((\bar J^2(q^{ 1}_i) a_i - \frac{1}{2}\Delta  (\bar J^2)'(q^{ 1}_i))
e^{-rs}  \Lambda^1_s  dL_s^{q^{ 1}_i} + e^{-r \tau^2_*} \Lambda^1_{\tau^2_*} R^2(X_{\tau^2_*}) 1_{\tau^2_* < {\tau}}   \nonumber \\
& & + e^{-r{\tau }}  \Lambda_{\tau}^1 \bar J^2(X_{\tau })1_{\tau \leq\tau^2_*} \Big] \label{exp2}
 \end{eqnarray}
\end{lemma}
As the proof will show,
Equation (\ref{fortanaka33}) follows from the strong Markov property and
Equation (\ref{exp2}) follows from an application of the It\^o-Meyer-Tanaka formula. Let us admit  for a while Lemma \ref{CN}. It follows that, for every  stopping time $\tau$,
\[
\mathbf E_x \left[\sum_{i= 1}^{n} \int_{[0,{\tau})} 1_{s<\tau^2_*} e^{-rs} G^2(q^{ 1}_i) a_i\Lambda_s^1   dL_s^{q^{ 1}_i} \right]= \mathbf E_x \left[\sum_{i = 1}^{n} \int_{[0, \tau)}1_{s<\tau^2_*}e^{-rs}
(\bar J^2(q^{ 1}_i) a_i - \frac{1}{2}\Delta \bar {J^2}'(q^{ 1}_i))
  \Lambda^1_s  dL_s^{q^{ 1}_i} \right].\]
Equivalently, for every  stopping time $\tau$, $\mathbf E_x[M_\tau] = \mathbf E_x[M_0] = 0$ where
\[ 
M_t \equiv \sum_{i= 1}^{n} \int_{[0,{\tau})}1_{s<\tau^2_*} e^{-rs} ( a_i [G^2(q^{ 1}_i) - \bar J^2(q^{ 1}_i)  ]+ \frac{1}{2}\Delta \bar {J^2}'(q^{ 1}_i))  \Lambda^1_s  dL_s^{q^{ 1}_i}.
\]
It follows that the process $(M_t)_{t \geq 0}$ is a martingale.\footnote{ See for instance Proposition 3.5 page 70  in Revuz and Yor (1999).}   Noting  that $(M_t)_{t \geq 0}$ is a continuous bounded variation process, we deduce  that, for every stopping time $\tau$,
\begin{equation} \label{1-29} 
M_\tau = M_0 = 0.
\end{equation}
Let us assume that for some $i$ with  $1 \leq i \leq n$, 
\begin{equation}\label{2-29}
a_{i} [G^2(q^{ 1}_{i}) - \bar J^2(q^{ 1}_{i})  ]+ \frac{1}{2}\Delta \bar {J^2}'(q^{ 1}_{i}) \neq 0.
\end{equation}
Let us consider  $X_0 = q^{ 1}_{i}$ and $T_{\varepsilon} := \inf \{t \geq 0 \, | \, X_t \notin (q^{ 1}_{i} - \varepsilon, q^{ 1}_{i} +\varepsilon) \}$ where $\varepsilon >  0$ is such that  $q^{ 1}_{i+ 1}<q^{ 1}_{i} -\varepsilon < q^{ 1}_{i} +\varepsilon < q^{ 1}_{i- 1}$ (with the convention that $q^{ 1}_{- 1} = \beta$ and $q^1_{n+1}=s^2$). From the  properties of the local time process we have that, for all $t>0$, $L_t^{q^{ 1}_{i}} >0$ $\mathbf{P}_{q^{ 1}_{i}}- a.s$.\footnote{ See for instance the proof of  Proposition 2.5 page 241 in  Revuz and Yor (1999).}    It then follows from (\ref{2-29}) that $M_{T_{\varepsilon}} \neq 0$, which contradicts (\ref{1-29}).  We conclude that $\bar J^2$ satisfies  (\ref{var14g'}) and that $\bar J^2$ is a solution to $\CV\CS^2$.

An analogous reasoning shows that if  $\bar J^1$ is the best reply value function to strategy $(\mu^2, S^2)$ then $\bar J^1$ is a solution to $\CV\CS^1$. \\

Let us  now prove Lemma \ref{CN}.

\noindent {\bf Proof of Lemma \ref{CN}}  
 From (v) of Proposition \ref{geneprop}, if $((\mu^1,S^1), (\mu^2, S^2))$ is a MPE, then $(0, (\alpha, s^2])$ is a best reply to $(\mu^1, S^1)$. Applying  the strong Markov property \eqref{MarkovpropertyJ} to 
the expected payoff of player 2 associated to the pair of Markov strategies
$((\mu^1,S^1), (0, (\alpha, s^2]))$ yields for every stopping time $\tau$ and every $x \in {\cal I}$
\begin{eqnarray*}
\bar J^2(x)  & = 
 \mathbf E_x \Big[ &\sum_{i= 1}^{n} \int_{[0,{\tau}\wedge \tau^2_*)} e^{-rs} G^2(q^{ 1}_i) \Lambda_s^1  a_i dL_s^{q^{ 1}_i} + e^{-r \tau_{s^{ 2} }} \Lambda^1_{\tau_{s^{ 2} }} R^2(X_{\tau^2_*}) 1_{\tau^2_* < {\tau}}   \nonumber \\
& & + e^{-r{\tau }}  \Lambda_{\tau}^1 \bar J^2(X_{\tau })1_{\tau \leq\tau^2_*} \Big]
\end{eqnarray*}
which shows  (\ref{fortanaka33}).

To prove (\ref{exp2}), we apply  the Ito-Meyer-Tanaka formula  to the process \break$e^{-r  (\tau \wedge \tau^2_* \wedge \tau_k)}  \Lambda_{\tau \wedge \tau^2_* \wedge \tau_k}^1 \bar J^2 ( X_{\tau \wedge \tau^2_* \wedge\tau_k})$ where, for each $k\in \NN$,  $\tau_k := \inf \{t \geq 0 \, : \, X_t \neq [\alpha_k , \beta_k ] \}$ with
  $([\alpha_k, \beta_k])_{k \in \NN}$  an increasing sequence of compacts intervals of ${\cal I}$ such that $\cup_{k \in \NN} [\alpha_k, \beta_k] = (\alpha, \beta)$.   Observe that $\tau_k  < \infty$ and that $X_t \in [\alpha_k, \beta_k]$ over $\{ t \leq \tau_k \}$, $\mathbf{P}_x$-almost surely for $x\in [\alpha_k,\beta_k]$.  Also, because $X$  is a regular diffusion, $\lim_{k \longrightarrow \infty} \tau_k = \infty$ and, hence, $\lim_{k \longrightarrow \infty} \tau \wedge \tau_k = \tau$ with  $\tau \in {\cal T}$.

We thus obtain from the Ito-Meyer-Tanaka formula 
\begin{eqnarray}
\bar J^2(x) &= &
e^{-r (\tau \wedge \tau^2_* \wedge\tau_k)} \Lambda^1_{\tau \wedge \tau^2_* \wedge\tau_k} \bar J^2(X_{\tau \wedge\tau^2_* \wedge \tau_k}) -
\int_{[0,{\tau \wedge\tau^2_* \wedge \tau_k} )}e^{-rs} \bar J^2(X_s)d\Lambda_s^1 \nonumber \\
& &+ \int_{[0,{\tau \wedge \tau^2_* \wedge\tau_k} )} e^{-rs} \Lambda_s^1 (r \bar J^2(X_s) - {\cal L} \bar J^2(X_s))  \Pi_{i =1}^{n} 1_{\{X_s \neq q^{ 1}_i\}}\, ds  \nonumber\\
& &- \int_{[0,{\tau \wedge \tau^2_* \wedge\tau_k})}  e^{-rs} \Lambda_s^1 \sigma(X_s) \bar {J^2}'(X_s) \Pi_{i =1}^{n} 1_{\{X_s \neq q^{ 1}_i\}}dW_s \\
& &  - \, \frac{1}{2} \sum_{i =1}^{n}\Delta \bar {J^2}' (q^{ 1}_i) \int_{[0,{\tau \wedge\tau^2_* \wedge \tau_k})}   e^{-rs} \Lambda_s^2  dL_s^{q^{ 1}_i}. \label{imt'}
\end{eqnarray}
After taking  expectations, we get,
\begin{eqnarray*}
\bar J^2(x) = \mathbf E_x [e^{-r{\tau \wedge\tau^2_* \wedge \tau_k}}  \Lambda_{\tau \wedge \tau^2_* \wedge\tau_k}^1 \bar J^2(X_{\tau \wedge\tau^2_* \wedge \tau_k})]  -  \mathbf E_x\Big[\int_{[0,{\tau \wedge \tau^2_* \wedge\tau_k})}e^{-rs} \bar J^2(X_s) d \Lambda_s^1\Big] \nonumber \\
- \mathbf E_x\Big[\frac{1}{2} \sum_{i =1}^{n} \Delta  ( \bar J^2)'(q^{ 1}_i) \int_{[0, {\tau \wedge\tau^2_* \wedge \tau_k})} e^{-rs} \Lambda_s^1 dL_s^{q^{ 1}_i}\Big], \end{eqnarray*}
where we have used the fact that $\bar J^2$ satisfies  (\ref{var11g}) and that, 
\begin{equation} \label{loc}
\mathbf E_x\Big[\int_{[0, \tau \wedge \tau^2_*  \wedge\tau_k)} e^{-rs} \Lambda_s^1 \sigma (X_s) \bar {J^2}' (X_s) \Pi_{i =1}^{n} 1_{\{X_s \neq q^{ 1}_i \}}d W_s\Big] = 0.
\end{equation}
To justify this latter equality, note that $\sigma$ is continuous on $I$, and  $\bar J^2\in { C \cal }^1 (I \setminus \{(q_i^1)_{1\leq i \leq n}\})$ with $|(\bar {J^2})' (x+)| < \infty $ and  $|(\bar{ J^2})' (x-)| < \infty $ for $x\in \{(q^{ 1}_i)_{1 \leq i \leq n}\}$. Thus,
there exists $C_k >0$ such that $|\sigma (X_t)  (\bar J^2)'(X_t) | \leq C_k$ over $\{t \leq   \tau^2_*  \wedge \tau_k\}$, $\mathbf{P}_x$-almost surely, which implies (\ref{loc}).
It follows that,
\begin{eqnarray*}
\bar J^2(x) & =&  \mathbf E_x[e^{-r{\tau \wedge \tau_k}}  \Lambda_{(\tau \wedge \tau_k)}^1 \bar J^2(X_{\tau \wedge \tau_k})1_{ {\tau \wedge \tau_k}\leq\tau^2_* }] + \mathbf E_x [e^{-r \tau^2_*} 
\Lambda^1_{\tau^2_*}  R^2 (X_{\tau^2_*}) 1_{\tau^2_* <{\tau \wedge \tau_k}} ]\\
& & +\, \mathbf E_x \Big [\sum_{i = 1}^{n } \int_{[0, \tau \wedge \tau^2_* \wedge \tau_k)} e^{-rs} \bar J^2(X_s) \Lambda^1_s a_i dL_s^{q^{ 1}_i} \Big ] -
\mathbf E_x \Big[ \frac{1}{2} \sum_{i =1}^{n} \Delta \bar {J^2}'(q^{ 1}_i) \int_{[0, {\tau \wedge\tau^2_* \wedge \tau_k})} e^{-rs} \Lambda_s^1 dL_s^{q^{ 1}_i}\Big]. \end{eqnarray*}
That is,
\begin{eqnarray}
\bar J^2(x) =  \mathbf E_x[e^{-r{\tau \wedge \tau_k}}  \Lambda_{(\tau \wedge \tau_k)^-}^1 \bar J^2(X_{\tau \wedge \tau_k})1_{ {\tau \wedge \tau_k}\leq\tau^2_* }] +
 \mathbf E_x [e^{-r \tau_{q^{ 1}_n }} 
\Lambda^1_{\tau_{q^{ 1}_n  }} R^2(X_{\tau^2_*}) 1_{\tau^2_*< {\tau \wedge \tau_k}} ] \nonumber\\
+ \,\mathbf E_x \Big [\sum_{i = 1}^{n } \int_{[0, \tau \wedge \tau^2_* \wedge \tau_k)}
((\bar J^2(q^{ 1}_i) a_i - \frac{1}{2}\Delta \bar {J^2}'(q^{ 1}_i))
e^{-rs}  \Lambda^1_s  dL_s^{q^{ 1}_i} \Big ] 
 \label{01121}
 \end{eqnarray}
We have that
\begin{eqnarray*}
\lim_{k \longrightarrow \infty} \mathbf E_x \Big [\int_{[0, \tau \wedge \tau^2_* \wedge \tau_k)}
((\bar J^2(q^{ 1}_i) a_i - \frac{1}{2}\Delta \bar {J^2}'(q^{ 1}_i))
e^{-rs}  \Lambda^1_s  dL_s^{q^{ 1}_i} \Big]\\ = \mathbf E_x \Big[\int_{[0, \tau \wedge  \tau^2_* )}
((\bar J^2(q^{ 1}_i) a_i - \frac{1}{2}\Delta \bar {J^2}'(q^{ 1}_i))
e^{-rs}  \Lambda^1_s  dL_s^{q^{ 1}_i} \Big]
\end{eqnarray*}
and
$$\lim_{k \longrightarrow \infty} \mathbf E_x [e^{-r \tau_{s^{ 2} }} 
\Lambda^1_{\tau_{s^{ 2}  }} R^2(X_{\tau^2_*}) 1_{\tau^2_* < {\tau \wedge \tau_k}} ]=\mathbf E_x [e^{-r \tau_{s^{ 2} }} 
\Lambda^1_{\tau_{s^{ 2}  }} R^2(X_{\tau^2_*}) 1_{\tau^2_*< {\tau }} ]$$
by the monotone convergence theorem.
Because from Proposition \ref{geneprop} $\bar J^2 \leq G^2$, it follows from  assumption {\bf A4} that  the sequence 
 $(e^{-r \tau \wedge \tau_k} \bar J^2(X_{\tau^- \wedge \tau_k})1_{ {\tau \wedge \tau_k}\leq\tau^2_*}1_{\tau < \infty})_{k\in \NN}$ is uniformly integrable. 
We have 
\[ \lim_{k \rightarrow  \infty} \mathbf E_x[ e^{- r( \tau \wedge \tau_k)} \Lambda^1_{(\tau \wedge \tau_k)} \bar J^2 (X_{\tau \wedge \tau_k})1_{ {\tau \wedge \tau_k}\leq\tau^2_*}1_{\tau < \infty}] =  \mathbf E_x [e^{- r \tau} \Lambda^1_{ \tau} \bar J^2 (X_{ \tau})1_{ {\tau }\leq\tau^2_*}1_{\tau < \infty}]\]
by Vitali's convergence theorem.

Finally, over $\{\tau = \infty \}$, we have $e^{- r( \tau \wedge \tau_k)} \Lambda^1_{(\tau \wedge \tau_k)} \bar J^2 (X_{\tau \wedge \tau_k})1_{ {\tau \wedge \tau_k}\leq\tau^2_*} =e^{- r \tau_k} \Lambda^1_{\tau_k} \bar J^2 (X_{ \tau_k})1_{ { \tau_k}\leq\tau^2_*}$. For $k$ large enough, $x \in (\alpha_k, \beta_k)$, thus 
\begin{align*}
\mathbf E_x[e^{- r \tau_k} \Lambda^1_{\tau_k} \bar J^2 (X_{ \tau_k})1_{ { \tau_k}\leq\tau^2_*}] &\leq \mathbf E_x [e^{- r \tau_k} \Lambda^1_{\tau_k} \bar J^2 (X_{ \tau_k}) 1_{\{X_{\tau_k}  = \alpha_k\}}] + E_x [e^{- r \tau_k} \Lambda^1_{\tau_k} \bar J^2 (X_{ \tau_k}) 1_{\{X_{\tau_k}  = \beta_k\}}]\\
&= \mathbf E_x [e^{- r \tau_{\alpha_k}} \Lambda^1_{\tau_{\alpha_k}} \bar J^2 (\alpha_k)  ] + E_x [e^{- r \tau_{ \beta_k}} \Lambda^1_{\tau_{\beta_k}} \bar J^2 (\beta_k) ]\\
&\leq \frac{\phi(x)}{\phi(\alpha_k)} G^2(\alpha_k) + \frac{\psi(x)}{\psi(\beta_k)}G^2(\beta_k), 
\end{align*}
We deduce from the growth properties (\ref{gp'}) that
\[\lim_{k \longrightarrow \infty}\mathbf E_x[e^{- r \tau_k} \Lambda^1_{\tau_k} \bar J^2 (X_{ \tau_k})1_{ { \tau_k}\leq\tau^2_*}] = 0.\]
Thus, 
letting $k$ go to $\infty$ in (\ref{01121}) yields
\begin{eqnarray*}
\bar J^2(x)  &=& \mathbf E_x \Big[  \sum_{i = 1}^{n } \int_{[0, \tau\wedge \tau^2_*)}
((\bar J^2(q^{ 1}_i) a_i - \frac{1}{2}\Delta  (\bar J^2)'(q^{ 1}_i))
e^{-rs}  \Lambda^1_s  dL_s^{q^{ 1}_i} + e^{-r \tau^2_*} \Lambda^1_{\tau^2_*} R^2(X_{\tau^2_*}) 1_{\tau^2_* < {\tau}}   \nonumber \\
& & + e^{-r{\tau }}  \Lambda_{\tau}^1 \bar J^2(X_{\tau })1_{\tau \leq\tau^2_*} \Big],
\end{eqnarray*}
which shows (\ref{exp2}).\\

To close the proof of Theorem \ref{CSmixed} it simply remains to establish 
 a verification Lemma that shows that, if a pair $(w^1,w^2)$  are solutions to the systems $\CV\CS^1$ and $\CV\CS^2$  then,
\[((\mu^1, S^1), (\mu^2, S^2)) = \left (\left (\sum_{i =1}^n a_i \delta_{q^{ 1}_i}, \emptyset \right), \left (\sum_{i = 1}^{n- 1} b_i \delta_{q^{ 2}_i}, (\alpha,s^{ 2}]\right )\right )\]
 is a mixed strategy MPE.  
 Specifically,
 Lemma \ref{verifg1} below shows that 
the strategy  $(\sum_{i=1}^{n} a_i \delta_{q^{ 1}_{i}}, \emptyset)$ is a PBR to player's 2 strategy  $(\sum_{i=1}^{n-1} b_i \delta_{q^{ 2}_{i}}, (\alpha,s^{ 2}])$ and  that 
the strategy    $(\sum_{i=1}^{n-1} b_i \delta_{q^{ 2}_{i}}, (\alpha,s^{ 2}])$
 is a PBR to player's 1 strategy $(\sum_{i=1}^{n} a_i \delta_{q^{ 1}_{i}}, \emptyset)$. This ends the proof of Theorem \ref{CSmixed}. 
\begin{lemma} \label{verifg1}
For $i=1,2$, let $w^i$ be a solution to $\CV\CS^i$ then, with $j\neq i$
\begin{eqnarray}
w^i(x) & \geq &  \sup_{\tau}J^i(x, \tau, (\mu^{j}, S^j)), \label{npbr11}\\
w^i(x) & = & J^i(x,(\mu^1, S^1), (\mu^2, S^2)). \label{npbr11'}
\end{eqnarray}
 \end{lemma}

  \noindent
 {\bf Proof of Lemma \ref{verifg1}.}
We prove Lemma \ref{verifg1} for $i= 2$ and $j= 1$.

At first, let us observe that (\ref{MarkovpropertyJ}) leads to
$$J^2(x, (\mu^1, S^1), \tau) = 
 \mathbf E_x[e^{-r \tau} R^2(X_{\tau}) \Lambda_{\tau}^1 \;+ \;\sum_{i =1}^{ n}\int_{[0, \tau)} e^{-rs} G^2 (X_s) \Lambda_s^1  a_i \, d L_s^{q^{ 1}_i} ].$$

Let  $w^2$ be a solution to ${\cal VS}^2$. 
 We apply the It\^o-Meyer-Tanaka formula to the process  $(e^{-r (\tau \wedge \tau_k)} \Lambda^1_{\tau \wedge \tau_k} w^2 (X_{\tau \wedge \tau_k}))_{t \geq 0}$ where $\tau_k := \inf \{t \geq 0 \, : \, X_t \neq [\alpha_k , \beta_k ] \}$ and where $([\alpha_k, \beta_k])_{k \in \NN}$ be an increasing sequence of compact intervals of ${\cal I}$ such that $\cup_{k \in \NN} [\alpha_k, \beta_k] = {\cal I}$ as in the proof of Lemma \ref{CN}.
We get
\begin{eqnarray}
 w^2(x) &= &
e^{-r (\tau  \wedge\tau_k)} \Lambda^1_{\tau  \wedge\tau_k}  w^2(X_{\tau \wedge \tau_k}) -
\int_{[0,{\tau \wedge \tau_k} )}e^{-rs}  w^2(X_s)d\Lambda_s^1 \nonumber \\
& &+ \int_{[0,{\tau \wedge\tau_k} )} e^{-rs} \Lambda_s^1 (r w^2(X_s) - {\cal L} w^2(X_s))  \Pi_{i =1}^{n} 1_{\{X_s \neq q^{ 1}_i\}}\, ds  \nonumber\\
& &- \int_{[0,{\tau \wedge\tau_k})}  e^{-rs} \Lambda_s^1 \sigma(X_s)  {w^2}'(X_s) \Pi_{i =1}^{n} 1_{\{X_s \neq q^{ 1}_i\}}dW_s \label{br}\\
& &  - \, \frac{1}{2} \sum_{i =1}^{n}\Delta  {w^2}' (q^{ 1}_i) \int_{[0,{\tau  \wedge \tau_k})}   e^{-rs} \Lambda_s^1  dL_s^{q^{ 1}_i}. \label{imt'}
\end{eqnarray}
From assumption {\bf A3} and (\ref{var11g'}), we have ${\cal L}w^2 - r w^2 = {\cal L}R^2 - r R^2 \leq 0$ on $(\alpha, s^{2})$.
It then follows from  (\ref{var11g}) that
  \begin{equation} \label{majhjb}
\mathbf E_x\left [\int_{[0,{\tau \wedge\tau_k} )} e^{-rs} \Lambda_s^1 (r w^2(X_s) - {\cal L} w^2(X_s))  \Pi_{i =1}^{n} 1_{\{X_s \neq q^{ 1}_i\}}\, ds \right ] \geq 0.
\end{equation}
From (\ref{var14g'}) and from the properties of the local time process $(L_t^{q^{ 1}_i})_{t \geq 0}$, we have  
\begin{eqnarray}
 - \, \frac{1}{2} \sum_{i =1}^{n}\Delta  {w^2}' (q^{ 1}_i) \int_{[0,{\tau \wedge \tau_k})}   e^{-rs} \Lambda_s^1  dL_s^{q^{ 1}_i}
=  \sum_{i = 1}^{n} a_i (G^2 (q^{ 1}_i ) - w^2 (q^{ 1}_i ))
\mathbf E_x[\int_{[0,{\tau  \wedge \tau_k})}   e^{-rs} \Lambda_s^1  dL_s^{q^{ 1}_i}
] \nonumber\\
= \mathbf E_x [ \sum_{i = 1}^{n}  \int_{[0,{\tau \wedge \tau_k})} e^{-rs} G^2(X_s) \Lambda^1_s  a_i d L_s^{q^{ 1}_i} ] - \mathbf E_x [  \sum_{i= 1}^{n}\int_{[0,{\tau \wedge \tau_k})} e^{-rs}\Lambda_s^1 w^2(X_s) a_i d L_s^{q^{ 1}_i}] \nonumber \\
= \mathbf E_x [ \sum_{i = 1}^{n}  \int_{[0,{\tau \wedge \tau_k})} e^{-rs} G^2(X_s) \Lambda^1_s  a_i d L_s^{q^{ 1}_i} ] +
\mathbf E_x[\int_{[0,{\tau \wedge \tau_k})} e^{-rs} w^2(X_s)d\Lambda_s^1]\nonumber \\
\label{tana}
\end{eqnarray}
We obtain from (\ref{imt'}), (\ref{majhjb}), (\ref{tana})
\begin{eqnarray*}
w^2(x)& \geq &\mathbf E_x [e^{-r(\tau\wedge \tau_k)}\Lambda_{\tau \wedge \tau_k}^1 w^2(X_{\tau \wedge \tau_k})] 
 + \mathbf E_x[\sum_{i = 1}^{n}\int_{[0,{\tau \wedge \tau_k})}  e^{-rs} G^2(X_s) \Lambda^1_s a_i d L_s^{q^{ 1}_i} ]\\
& \geq & \mathbf E_x [e^{-r(\tau\wedge \tau_k)}\Lambda_{\tau \wedge \tau_k}^1 R^2(X_{\tau \wedge \tau_k})] 
 + \mathbf E_x[\sum_{i = 1}^{n}\int_{[0,{\tau \wedge \tau_k})}  e^{-rs} G^2(X_s) \Lambda^1_s a_i d L_s^{q^{ 1}_i} ]
\end{eqnarray*}
where we have used that the stochastic integral in (\ref{br}) is a centered square integrable martingale as shown in Lemma \ref{CN}  and that $w^2(x)  \geq R^2(x) $ on ${\cal I}$. Using again the same arguments than in Lemma \ref{CN}, letting $k$ go to $\infty$ yields
$$w^2(x) \geq  \mathbf E_x [e^{-r(\tau)}\Lambda_{\tau }^1 R^2(X_{\tau })] 
 + \mathbf E_x[\sum_{i = 1}^{n}\int_{[0,{\tau })}  e^{-rs} G^2(X_s) \Lambda^1_s a_i d L_s^{q^{ 1}_i} ],$$ where we note that, from (\ref{forbpr}), the left-hand-side of the above equation corresponds to $J^2(x,(\mu^1, S^1), \tau)$.
Then,  taking the supremum over $\tau$  yields (\ref{npbr11}).

To establish (\ref{npbr11'}), we apply the Ito-Meyer-Tanaka formula to  $e^{-r \tau_k} \Lambda_{\tau_k}^1 \Lambda_{\tau_k}^2 w^2 ( X_{\tau_k})$ where   $w^2$ is a solution to ${\cal VS}^2$. We obtain that 
\begin{eqnarray}
w^2(x) = \mathbf E_x[e^{-r{\tau_k}} \Lambda_{\tau_k}^1 \Lambda_{\tau_k}^2 w^2(X_{\tau_k})] - \mathbf E_x[\int_{[0,{\tau_k})}e^{-rs} w^2(X_s) \Lambda_s^2 d \Lambda_s^1]  -  \mathbf E_x[\int_{[0,{\tau_k})}e^{-rs} w^2(X_s) \Lambda_{s}^1 d \Lambda_s^2] \nonumber \\
- \mathbf E_x[\frac{1}{2} \sum_{i =1}^{n} \Delta {w^2}'(q^{ 1}_i) \int_{[0, {\tau_k})} e^{-rs} \Lambda_s^1 \Lambda_s^2 dL_s^{q^{ 1}_i}].\nonumber\\ \label{601}
\end{eqnarray}
where, as in the proof of Lemma \ref{CN}, we have used that
$$\mathbf E_x\Big[\int_{[0, \tau_k)} e^{-rs} \Lambda_s^1 \Lambda_s^2 \sigma (X_s)  {w^2}' (X_s) \Pi_{i =1}^{n} 1_{\{X_s \neq q^{ 1}_i \}}d W_s\Big] = 0.$$
and
$$\mathbf E_x[ \int_{[0,{\tau  \wedge\tau_k} )} e^{-rs} \Lambda_s^1 \Lambda_s^2(r  w^2(X_s) - {\cal L} \bar w^2(X_s))  \Pi_{i =1}^{n} 1_{\{X_s \neq q^{ 1}_i\}}\, ds ] =0.$$ This latter equality follows from  (\ref{var11g})  and from the definition of $\Lambda_s^2$ which contains the indicator $1_{s <\tau^2_*}
$. 
Next, we have 
\begin{eqnarray}
- \mathbf E_x[\int_{[0,{\tau_k})}e^{-rs} w^2(X_s) \Lambda_{s}^1 d \Lambda_s^2] & = &
\mathbf E_x [\sum_{i = 1}^{n-1} \int_{[0,{\tau_k})} e^{-rs} w^2(X_s) \Lambda_{s}^1 \Lambda_s^2 b_i d L_s^{q^{ 2}_i}] \nonumber \\
&  &+\, \mathbf E_x[ e^{-r \tau^2_*} e^{-rs} w^2(X_{\tau_{ s^{ 2}}})  \Lambda_{\tau^2_*}^1  \Lambda_{\tau^2_*}^2 1_{\tau^2_* \leq {\tau_k}}] \nonumber\\ &= &
\mathbf E_x [\sum_{i = 1}^{n-1} \int_{[0,{\tau_k})} e^{-rs} R^2(X_s) \Lambda_{s}^1 \Lambda_s^2 b_i d L_s^{q^{ 2}_i}] \nonumber \\
&  &+\, \mathbf E_x[ e^{-r \tau^2_*} e^{-rs} R^2(X_{\tau_{ s^{ 2}}})  \Lambda_{\tau^2_*}^1  \Lambda_{\tau^2_*}^2 1_{\tau^2_* \leq {\tau_k}}]
\label{602}.
\end{eqnarray}
Using (\ref{var15n}), we get
\begin{eqnarray}
-\mathbf E_x [ \frac{1}{2} \sum_{i =1}^{n} \Delta {w^2}'(q^{ 1}_i) \int_{[0, {\tau_k})} e^{-rs} \Lambda_s^1 \Lambda_s^2 dL_s^{q^{ 1}_i}] 
& =& \mathbf E_x [\sum_{i = 1}^{n} \int_{[0,{\tau_k})} e^{-rs} G^2(q_i^{ 1}) \Lambda_s^1 \Lambda_s^2 a_i d L_s^{q^{ 1}_i}] \nonumber\\
&   & -\mathbf E_x [\sum_{i = 1}^{n} \int_{[0,{\tau_k})} e^{-rs} w^2(q_i^{ 1}) \Lambda_s^1 \Lambda_s^2 a_i d L_s^{q^{ 1}_i}] \nonumber\\
&= &\mathbf E_x [\sum_{i = 1}^{n} \int_{[0,{\tau_k})} e^{-rs}G^2(X_s) \Lambda_s^1 \Lambda_s^2 a_i d L_s^{q^{ 1}_i}] \nonumber \\
&  & + \mathbf E_x[\int_{[0,{\tau_k})}e^{-rs} w^2(X_s) \Lambda_s^2 d \Lambda_s^1]\label{603}.
\end{eqnarray}
We get from (\ref{601}), (\ref{602}), (\ref{603}) that,
\begin{eqnarray}
w(x) = \mathbf E_x[e^{-r\tau_k} \Lambda_{\tau_k}^1 \Lambda_{\tau_k}^2 w(X_{\tau_k})] + \mathbf E_x [\sum_{i = 1}^{n-1} \int_{[0,{\tau_k})} e^{-rs} R^2(X_s) \Lambda_{s}^1 \Lambda_s^2 b_i d L_s^{q^{ 2}_i}] \nonumber \\
 + \mathbf E_x[e^{-r \tau^2_*} R^2(X_{ s^{ 2}}) \Lambda_{\tau^2_*}^1  \Lambda_{\tau^2_*}^2 1_{\tau^2_* \leq {\tau_k}}]
+ \mathbf E_x [\sum_{i = 1}^{n} \int_{ [0,\tau_k)} e^{-rs} G^2(X_s) \Lambda_s^1 \Lambda_s^2 a_i d L_s^{q^{ 1}_i}]. \label{604}
\end{eqnarray}
 Letting $k$ go to $\infty$ as in Lemma \ref{CN}, yields
\begin{eqnarray*}
w^2(x) =  
\mathbf E_x[\sum_{i =1}^{ n-1}\int_{[0, \infty )} e^{-rs} R^2 (X_s) \Lambda_{s}^1 \Lambda_s^2 b_i \, d L_s^{q^{ 2}_i} \nonumber \\
  + \;  e^{-r \tau^2_*} R^2(X_{\tau_{ s^{ 2}}})  \Lambda_{\tau_{{s^{ 2}}}}^1  \Lambda_{\tau_{{s^{ 2}}}}^2 \;+ \;\sum_{i =1}^{ n}\int_{[0,{\infty})} e^{-rs} G^2 (X_s) \Lambda_s^1 \Lambda_{s}^2 a_i \, d L_s^{q^{ 1}_i} ].
\end{eqnarray*}
Finally, we observe   from (\ref{forbpr}) that  the right hand side of the latter equation corresponds to  $J^2(x, (\mu^1, S^1), (\mu^2, S^2))$, which ends the proof.

\end{document}